\documentclass[10pt]{mythesis}

\usepackage{braket}
\usepackage{amsmath2000,amsthm2000,amssymb}
\usepackage{latexsym}
\usepackage{proof}
\usepackage[matrix,frame,arrow]{xy}
\usepackage{fancyhdr}
\usepackage[english]{babel}
\usepackage{rotating}
\usepackage{epic,eepic}
\usepackage{newfile}
\usepackage{textcomp}
\usepackage{mathrsfs}

\DeclareRobustCommand{\inference}[2]{\raisebox{-1.3ex}{\infer{#1}{#2}}}
\newtheorem{thm}{Theorem}[chapter]
\newtheorem{cor}[thm]{Corollary} \newtheorem{lem}[thm]{Lemma}
\newtheorem{prop}[thm]{Proposition} \newtheorem{por}[thm]{Porism}
\theoremstyle{definition} \newtheorem{defin}[thm]{Definition}
\theoremstyle{remark} \newtheorem{rem}[thm]{Remark}

\newoutputstream{que} \openoutputfile{questions.out}{que}
\DeclareRobustCommand{\question}[2]{%
  \refstepcounter{thm}
\addtostream{que}{%
\vspace{0.8ex}\noindent{\bf Question \thethm}. {\it #1} \vspace{0.8ex} \par #2 \par}%
\vspace{0.8ex} \noindent{\bf Question \thethm}. {\it #1} \vspace{0.8ex}}

\DeclareRobustCommand{\conjecture}[2]{%
  \refstepcounter{thm}
\addtostream{que}{%
\vspace{0.8ex}\noindent{\bf Conjecture \thethm}. {\it #1} \vspace{0.8ex} \par #2 \par}%
\vspace{0.8ex} \noindent{\bf Conjecture \thethm}. {\it #1} \vspace{0.8ex}}

\DeclareRobustCommand{\pair}[1]{\left\langle #1 \right\rangle}
\DeclareRobustCommand{\imp}{\rightarrow}
\DeclareRobustCommand{\ekv}{\leftrightarrow}
\DeclareRobustCommand{\Vee}{\bigvee} 

\DeclareRobustCommand{\abs}[1]{\lvert #1 \rvert}

\DeclareMathOperator{\uniform}{\star}

\DeclareRobustCommand{\existsone}{\exists !}
\DeclareRobustCommand{\least}{\mu}

\DeclareRobustCommand{\xor}{\mathrel{\dot{\vee}}}

 \let\geq=\geqslant  \let\leq=\leqslant

\DeclareMathOperator{\len}{len}
\DeclareMathOperator{\conc}{\raisebox{0.8ex}{$\frown$}}
\DeclareRobustCommand{\eqdef}{\mathrel{=_{\rm df}}}
\DeclareRobustCommand{\num}[1]{n_{#1}}

\DeclareRobustCommand{\inferforall}[2]{\ldots \ #2 \ \ldots \ _{#1}}

\DeclareMathOperator{\Diag}{Diag} \DeclareMathOperator{\ElDiag}{ElDiag}

\DeclareRobustCommand{\F}{\mathscr{F}}
\DeclareRobustCommand{\Sat}{\ensuremath{\Sigma}}

\DeclareRobustCommand{\structure}[1]{\ensuremath{\mathfrak{#1}}}
\DeclareRobustCommand{\M}{\structure{M}}
\DeclareRobustCommand{\Mprop}{\ensuremath{\M_\mathrm{p}}}
\DeclareRobustCommand{\T}{\ensuremath{\mathfrak{T}}}

\DeclareMathOperator{\Sc}{\mathsf{S}} \DeclareRobustCommand{\v}{{\mathsf{v}}}
\DeclareRobustCommand{\c}{{\mathsf{c}}}
\DeclareRobustCommand{\Q}{{\mathop{\mathbf{\mathsf{Q}}}}}

\DeclareRobustCommand{\Nat}{\ensuremath{\mathbb{N}}}

\DeclareRobustCommand{\PA}{\ensuremath{\mathtt{PA}}}
\DeclareMathOperator{\Tr}{\mathtt{Tr}}

\DeclareMathOperator{\satcl}{{\mathtt{SatCl}}}
\DeclareMathOperator{\Pf}{{\mathtt{Pf}}}
\DeclareMathOperator{\form}{{\mathtt{Form}}}
\DeclareMathOperator{\term}{\mathtt{Term}}
\DeclareMathOperator{\sent}{{\mathtt{Sent}}}
\DeclareMathOperator{\clterm}{{\mathtt{ClTerm}}}

\DeclareMathOperator{\Th}{{\mathtt{Th}}}

\DeclareMathOperator{\FV}{{\mathtt{FV}}}
\DeclareMathOperator{\val}{{\mathtt{val}}}

\DeclareMathOperator{\Ax}{\mathtt{Ax}}

\DeclareRobustCommand{\Te}{{T}}

\DeclareRobustCommand{\apprx}[1]{\mathop\mathrm{app}(#1)}
\DeclareRobustCommand{\tmpl}[1]{\widehat{#1}}

\let\setminus=\setminus

\DeclareRobustCommand{\models}{\mathrel{\vDash}}
\DeclareRobustCommand{\nmodels}{\mathrel{\nvDash}}
\DeclareRobustCommand{\prf}[1][{}]{\mathrel{\vdash^{\!#1}}}

\DeclareRobustCommand{\mprf}[1][{}]{\mathrel{\vdash_{\!\!\!\M}^{\!  #1}}}
\DeclareRobustCommand{\mpropprf}[1][{}]{\mathrel{\vdash_{\!\!\!\Mprop}^{\!#1}\!\!\!}}

\DeclareRobustCommand{\propprf}{\mathrel{\vdash_{\!\!\!\raisebox{0.2ex}{$\scriptstyle
        \mathrm p$}}^{\!\!\!*}}}
\DeclareRobustCommand{\predprf}{\mathrel{\vdash^{\!\!\!*}}}
\DeclareRobustCommand{\nmprf}[1][{}]{\mathrel{\nvdash_{\!\!\!\M}^{\!#1}}}
\DeclareRobustCommand{\nmpropprf}[1][{}]{\mathrel{\nvdash_{\!\!\!\Mprop}^{\!#1}}}
\DeclareRobustCommand{\tprf}[1][{}]{\mathrel{\vdash_{\!\!\!\mathtt{T}}^{\!#1}}}
\DeclareRobustCommand{\ntprf}[1][{}]{\mathrel{\nvdash_{\!\!\!\mathtt{T}}^{\!#1}}}
\DeclareRobustCommand{\mmodels}{\mathrel{\models_{\!\!\!\raisebox{-0.4ex}{\scriptsize
        \M}}}}

\DeclareRobustCommand{\L}{\ensuremath{\mathscr{L}}}
\DeclareRobustCommand{\LM}{\ensuremath{^*\!\!\L_\M}}
\DeclareRobustCommand{\LMs}{\ensuremath{\L_\M}}
\DeclareRobustCommand{\LT}{\ensuremath{\L_{\mathtt{T}}}}
\DeclareRobustCommand{\LA}{\ensuremath{\L_A}}

\DeclareRobustCommand{\initialseg}[1]{I_{\mathord< #1}}

\newcommand{\forpropprf}[1]{\mathop{\mathtt{PropPrf}_{#1}}}
\newcommand{\forpredprf}[1]{\mathop{\mathtt{Prf}_{#1}}}

\newenvironment{case}{\left\{\begin{array}{ll}}{\end{array}\right.}

\def\godel#1{\,\godelalt{#1}\,}

\def\godelalt#1{\setbox0=\hbox{$#1$}
\dimen0=.3pt \dimen1=3.5pt 
\dimen2=1.5pt \dimen3=1.5pt
\dimen4=\dimen3 \advance\dimen4 by \ht0%
\raise\dimen4\leftcorner %
\dimen5=\dimen1 \advance\dimen5 by -\dimen2 %
\kern-\dimen5 %
#1 %
\kern-\dimen5 %
\raise\dimen4\rightcorner}

\def\hlineseg{\vbox to \dimen0{\hrule width\dimen1 depth0pt height\dimen0}}

\def\vlineseg{\hbox to \dimen0{\advance\dimen1 by -\dimen0 \vrule
    height\dimen0 depth\dimen1 width\dimen0}}

\def\leftcorner{\hbox{\vlineseg \kern-\dimen0 \hlineseg}}

\def\rightcorner{\hbox{\hlineseg \kern -\dimen0 \vlineseg}}

\begin{document}

\title{\Huge Satisfaction classes in nonstandard models of first-order
  arithmetic \\[3 cm]%
  \normalsize Revised version 1.1}
\author{Fredrik Engstr\"om} \date{} \year=2002 \month=8

\maketitle

\newpage \pagestyle{empty}

\vspace*{5cm} \noindent Satisfaction classes in nonstandard models
of first-order arithmetic\\
\noindent FREDRIK ENGSTR\"OM\\[5mm]
\noindent \copyright Fredrik S G Engstr\"om, 2002 \\[5mm]
\noindent ISSN 0347-2809/NO 2002:24 \\
\noindent Department of Mathematics\\
\noindent Chalmers University of Technology and G\"oteborg
University\\
\noindent 412 96 G\"oteborg \\
\noindent Sweden \\
\noindent Telephone +46 (0)31-772 1000

\vfill

\noindent Matematiskt centrum \\
\noindent G\"oteborg, Sweden 2002

\newpage

\vspace*{1 cm}

\noindent {\bf Abstract}

\vspace{5mm}

\noindent A satisfaction class is a set of nonstandard sentences
respecting Tarski's truth definition. We are mainly interested in full
satisfaction classes, i.e., satisfaction classes which decides all nonstandard
sentences. Kotlarski, Krajewski and Lachlan proved in 1981 that a countable
model of \PA\ admits a satisfaction class if and only if it is recursively
saturated. A proof of this fact is presented in detail in such a way that it
is adaptable to a language with function symbols. The idea that a satisfaction
class can only see finitely deep in a formula is extended to terms. The
definition gives rise to new notions of valuations of nonstandard terms; these
are investigated. The notion of a {\em free\/} satisfaction class is
introduced, it is a satisfaction class free of existential assumptions on
nonstandard terms.

It is well known that pathologies arise in some satisfaction classes. Ideas of
how to remove those are presented in the last chapter. This is done mainly by
adding inference rules to \M-logic. The consistency of many of these
extensions is left as an open question.

\vspace{3mm}

\noindent {\bf Keywords:} Satisfaction classes, Models of
arithmetic.

\vspace{3mm}

\noindent {\bf AMS 2000 Subject Classification:} 03C50, 03C62,
03H15

\vfill

\noindent {\bf  Acknowledgement}

\vspace{5mm}

\noindent My supervisor, Richard W. Kaye, deserves all my thanks
for letting me come to Birmingham under his inspiring supervision, and, of
course, for showing me `the way' the numerous times I got lost in the jungle
of logic. I would also like to thank Thierry Coquand, who is my supervisor in
G\"oteborg, and Jan Smith, who introduced me to logic, for their support which
made it possible for me to work in models of PA.

Also, I would like to thank all at Mathematical Sciences at G\"oteborg University and
Chalmers University of Technology, and the School of Mathematics and
Statistics at the University of Birmingham.  \vspace{3mm}

\noindent Thank you!

\tableofcontents
\chapter{Introduction}\label{chp:introduction}
\setcounter{page}{1} \pagestyle{fancy}

By the work of Skolem we know there are nonstandard models of Peano Arithmetic
(\PA) and we know how to arithmetise logic inside \PA\ due to G\"odel. It is
also easy to see that in any such nonstandard model there are nonstandard
elements which the model thinks are sentences. By Tarski's truth definition we
also know what it means for a standard sentence to be true. The obvious
question is:

\begin{quote}
  When is a nonstandard sentence true?
\end{quote}

Given a nonstandard model \M\ of \PA\ a satisfaction class is a
(non-definable) predicate which, in a special sense, is a truth definition for
nonstandard sentences, i.e., it respects Tarski's truth definition. More
formally, a satisfaction class is a set \Sat\ of standard and nonstandard
sentences, extending the elementary diagram of \M, such that
\begin{align*}
  \neg \varphi \in \Sat \quad &\text{iff} \quad \varphi \notin \Sat, \\
  \varphi \vee \psi \in \Sat \quad &\text{iff} \quad \text{$\varphi \in \Sat$
    or $\psi
    \in \Sat$,} \qquad \text{and}\\
  \mathop{\exists \v_i} \gamma \in \Sat \quad&\text{iff} \quad \text{there
    exists $a \in \M$ such that $\gamma[a/\v_i] \in \Sat$},
\end{align*}
for any nonstandard sentences $\varphi, \psi$ and $\exists \v_i \gamma$. It is
not obvious that we can construct such a set, in fact we cannot always do
this, depending on the saturation of \M.

This chapter is a short introduction to the subject with some motivation of
the study and a historical overview. Next chapter is intended as a review of
the prerequisites for this thesis; the results are given without proofs. In
Chapter~\ref{chp:satcl} we have rewritten the construction of satisfaction
classes in a new style, the language includes functions, as opposed to
\cite{Kotlarski.Krajewski.ea:81}, and the terms are treated as they should be
treated, i.e., satisfaction classes can only ``look'' finitely deep into them,
as opposed to \cite{Kaye:91*2}.

In Chapter~\ref{chp:wesatcl} we study some alternative definitions of
satisfaction classes which all are weaker in the sense that the sentence
$\exists x (t=x)$ does not have to be true. We call these satisfaction classes
`free' since they are free of existential assumptions on nonstandard terms.

The last chapter is devoted to pathological examples that arise in
satisfaction classes and how to remove them. We introduce satisfaction classes
closed under propositional proofs and even stronger notions.

All structures studied in this thesis are models of \PA\ in the language $\LA$
with symbols $\Set{\Sc,+,\cdot,0}$.

\section{Motivation}

From the philosophical point of view there is an obvious motivation for the
study of satisfaction classes. We know what truth is for standard formulas by
the truth definition of Tarski and with the arithmetisation of logic we have
the notion of nonstandard formulas. The obvious question is then, what does it
mean for a nonstandard formula to be true? The study of satisfaction classes
is an attempt to answer this question. In fact, satisfaction classes can be
seen as nonstandard models in the same way as complete consistent Henkin
theories represent their Henkin models. There are $2^{\aleph_0}$ satisfaction
classes in a countable recursively saturated model of \PA, which implies that
the structure of satisfaction classes is rich in some sense\footnote{In fact,
  the set of satisfaction classes is dense in the Stone space of the
  nonstandard Lindenbaum algebra, see \cite{Smith:84}.} and we can see the
study of satisfaction classes as a nonstandard model theory.

In \cite{Smith:84}, Smith characterises the recursively saturated models of
\PA\ in terms of satisfaction classes, he constructs a $\Sigma_1^1$ formula
characterising them. He also $\Delta_2^1$ characterises resplendent models in
terms of satisfaction classes.  This gives us some idea that the satisfaction
classes are more important, in a mathematical sense, then just as truth
definitions. For mathematicians there are a lot of questions to be answered by
this study. For example, one of my motivations for this thesis has been to
find characterisations of other interesting model theoretic properties, such
as saturation properties stronger than recursive saturation.

The study could also be seen as an example of how to work with ill founded
objects, there may even be applications to computer science.

The main reason for the existence of this thesis is to show that there are
many unstudied related notions of satisfaction classes.  Some of them arise
naturally when we add function symbols in the language and others when we try
to remove certain ``pathologies.''

\section{Historical background}

In 1963 Abraham Robinson published a paper, `On languages which are based on
nonstandard arithmetic' \cite{Robinson:63}, where he discusses syntax and
semantics for nonstandard languages. This is, as far as I know, the first time
nonstandard languages are defined and investigated explicitly. He does not use
the word satisfaction class, but he gives two different examples of semantics
for nonstandard languages. He calls them the internal and the external truth
definitions. The external one, defined by the help of Skolem operators, is
defined only for formulas with finite `Robinson-rank' (which is a complexity
measure on nonstandard formulas), therefore it is not a full satisfaction
class.

The internal truth definition is defined as $\Sat$ if
\[(\Nat,\Sat_0) \prec \pair{\M,\Sat},\footnote{In
  fact, Robinson's notion is a bit different, but for the purpose of this
  survey we can think of internal truth defined in this way.}\] where \Nat\ is
the standard model of arithmetic and $\Sat_0$ is the standard truth
definition, i.e.,
\[
\Sat_0=\Set{\varphi | \text{$\varphi$ is a standard sentence in \LA\ and $\Nat
    \models \varphi$}}.
\] 
He proves that the external and the internal truth definitions do not coincide
and leaves it, more or less, there.

Later Krajewski \cite{Krajewski:76} returns to the question of the semantics
of nonstandard languages. He defines satisfaction classes (even though his
definition is rather weak) and investigates some related notions. He also
proves that for some models of cardinality $\lambda$ there exists $2^\lambda$
full satisfaction classes. He does not mention the question of which models
admit satisfaction classes, he only proves that some specific models do.

In \cite{Kotlarski.Krajewski.ea:81} and \cite{Lachlan:81} this question is
answered. In the first paper it is shown that if a countable model is
recursively saturated then it admits a satisfaction class. This is done by
using a version of $\omega$-logic (the idea is due to Jeff Paris) and the
result is not surprising. What is more surprising is the result in the second
paper by Lachlan in which he shows that if a model (of any cardinality) admits
a satisfaction class then it is recursively saturated. He proves this using a
sort of overspill he gets from the satisfaction class (proving the result with
induction in the language $\LA \cup \Set{\Sat}$ is easy).

In \cite{Smith:84} Smith strengthens the results in the two papers above. He
shows that any resplendent model has a satisfaction class and is able to find
a $\Delta_2^1$ characterisation of resplendency. He also formulates Lachlan's
proof in a syntactical way which makes it possible to find a $\Sigma_1^1$
formula characterising recursive saturation.

Other important contributions to the study of satisfaction classes are
\cite{Kotlarski:85}, \cite{Kotlarski.Ratajczyk:90} and
\cite{Kotlarski.Ratajczyk:90*1}. These papers discuss the question of when the
structure $\pair{\M,\Sat}$ satisfies either induction over $\Sigma_k$ formulas
or full induction. In the first paper Kotlarski shows that the satisfaction
classes satisfying $\Delta_0$ induction are precisely those closed under
nonstandard proofs of first-order logic and including all nonstandard
instances of the axiom of induction. It is easy to see that if a model admits
such a satisfaction class then it satisfies the scheme of reflection:
\[
\text{`}\PA \prf \varphi\text{'} \imp \varphi.
\]

The definitions of satisfaction classes in the preceding works are all using a
relational language, i.e., the language $0,1,\sigma,\pi$ where $\sigma$ and
$\pi$ are ternary relational symbols, supposed to express addition and
multiplication. The theory \PA\ has to be extended by some axioms expressing
that these relational symbols are in fact functions. In \cite{Kaye:91*2} Kaye
investigates the case when \PA\ is expressed in the language with symbols:
$+$, $\cdot$, $<$, $0$ and $1$, where $+$ and $\cdot$ are binary function
symbols. Throughout this thesis we will be using the language $\LA$ which has
one unary function symbol $\Sc$, two binary function symbols $+$ and $\cdot$,
and one constant symbol $0$ (and the equality predicate $=$).

\section{Theorems, Propositions, Lemmas, Corollaries and Porisms}

We will use the term `theorem' sparsely, it is used to put extra emphasise on
an important result. Propositions are the results which have a value on their
own, and lemmas (or lemmata) are results which help us to prove propositions
or theorems.  Corollaries are simple consequences of propositions or theorems
(and in some rare occasions of lemmas), but the, somewhat, unusual term porism
is used for simple consequences of a {\em proof\/} of a proposition, theorem
or lemma. It could for example be a simple generalisation of a proposition
which you get by some minor modifications of the presented proof of the
proposition.  The term was used by Euclid, but he used it with a somewhat
different meaning, which is not entirely known, see \cite{Tweddle:00} for more
information.

\chapter{Prerequisites}

To be able to read this thesis the reader needs some background knowledge of
first-order logic, \cite{Mendelson:97} is more than enough, and also some
knowledge of first-order arithmetic (Peano Arithmetic), especially its model
theory. A good general reference for this is \cite{Kaye:91*2}. In this chapter
we will review some of the material, omitting the proofs of the results.

\section{Peano Arithmetic}

Throughout this thesis \M\ will be a structure in the language \LA\ with
symbols $\Sc$, $+$, $\cdot$ and $0$, where $\Sc$ is a unary function, $+$ and
$\cdot$ are binary function symbols and $0$ is a constant symbol. We will as
usual write terms and equalities in the more convenient way by using infix
notation, e.g., instead of writing $+\bigl(t_1,t_2\bigr)$ we will write
$t_1+t_2$.

We denote interpretations as usually; the interpretation of $\Sc$ in \M\ is
denoted by $\Sc^\M$, the interpretation of $+$ is denoted $+^\M$, and so on.

We will consider formulas built up from the logical connectives $\neg$ and
$\vee$ and the quantifier $\exists$. The other connectives and quantifiers are
considered to be abbreviations in the usual way, $\varphi \wedge \psi$ is {\em
  defined\/} to be $\neg\bigl(\neg \varphi \vee \neg\psi\bigr)$, $\varphi \imp
\psi$ is $\neg \varphi \vee \psi$, $\varphi \ekv \psi$ is $\bigl(\varphi \imp
\psi\bigr) \wedge \bigl(\psi \imp \varphi\bigr)$ and $\mathop{\forall \v_i}
\varphi$ is $\neg {\exists \v_i}\, \neg \varphi$. We also define exclusive or
$\varphi \xor \psi$ as $(\varphi \vee \psi) \wedge \neg (\varphi \wedge
\psi)$.  Later it will be important to note that all these abbreviations are
of constant depth, i.e., the depth only increases a constant number when
replacing the definiendum with the definiens. For example, the depth of
$\neg\bigl(\neg \varphi \vee \neg \psi\bigr)$ is always two more than that of
$\varphi \wedge \psi$.

The variables in the language are $\v_0,\v_1, \ldots$. Sometimes we will be a
bit sloppy in the notation and use $x,y,w,\ldots$ as names for variables.  If
$\varphi$ is a formula then
\[
\varphi[t_1,\ldots,t_k / \v_{i_1},\ldots,\v_{i_k}]
\]
will denote the formula you get by substituting all free occurrences of
$\v_{i_l}$ with the term $t_l$, we will always assume that $\v_{i_l}$ is free
for $t_l$. Sometimes, when it will not cause confusion, we will write
$\varphi(t)$ to mean the formula obtained from $\varphi$ by replacing all the
occurrences of the free variable {\em under consideration} by the term $t$,
i.e., $\varphi[t/\v_i]$ where $\v_i$ is the variable under consideration. In
short, we will adopt all the usual notations and abbreviations used in the
literature.

Let \LMs\ be the language $\LA \cup \Set{\c_a | a \in \M }$ ($0$ and
$\c_{0^\M}$ will be regarded as the same symbol; this is to simplify some
definitions below) and let $\ElDiag(\M)$ be the theory of the structure \M\ in
the language \LMs, i.e., all \LMs-formulas true in $\M^+$, where $\M^+$ is the
expanded model of $\M$ which interprets each symbol $\c_a$ as $a$. Mostly we
will not distinguish between the two structures \M\ and $\M^+$, hoping this
will not cause any confusion for the reader. Sometimes, mostly when dealing
with standard formulas and terms, we will identify the element $a \in \M$ with
the constant symbol $\c_a$.

$\text{Diag}(\M)$ is the set of all true standard atomic and negated atomic
formulas in the language \LMs, so $\text{Diag}(\M) \varsubsetneq \ElDiag(\M)$.

$\PA^-$ is the theory with the universal closures of
\begin{gather*}
  \Sc(x)=\Sc(y) \imp x=y,\\
  \Sc(x) \neq 0, \\
  x \neq 0 \imp \exists y \bigl(\Sc(y)=x\bigr),\\
  x+0=x,\\
  x+\Sc(y)=\Sc(x+y),\\
  x \cdot 0 = 0 \quad \text{and}\\
  x \cdot \Sc(y)=x\cdot y + x
\end{gather*}
as axioms. If we add the axiom scheme of induction:
\[
\mathop{\forall \bar{x}} \Bigl(\varphi(0,\bar{x}) \wedge {\forall y}
\bigl(\varphi(y,\bar{x}) \imp \varphi(\Sc(y),\bar{x})\bigr) \imp {\forall y}\,
\varphi(y,\bar{x})\Bigr)
\]
for all \LA-formulas $\varphi$, we get Peano Arithmetic or \PA\ for short.

The symbol \M\ will always be assumed to be a structure, in a language
extending \LA, satisfying \PA\ and not isomorphic to the standard model \Nat\ 
of \PA, i.e., \M\ will be assumed to be a nonstandard model of \PA.

The predicate $x<y$ is defined as 
\[
\exists z \bigl(z \neq 0 \wedge x+z=y\bigr).
\]
Once again, it will be important later that this definition is
of constant depth, i.e., $\exists z (z \neq 0 \wedge t+z =r)$ is of depth at
most four more than the depth of $t < r$.  We also define the function $x-1$
by the following equation
\[
x-1 \eqdef (\least z) \Bigl[\bigl(x=0 \wedge z=0\bigr) \vee \bigl(x \neq 0 \wedge
\Sc(z)=x\bigr) \Bigr],
\]
where $(\least x) \varphi(x)$ means `the least $x$ such that $\varphi(x)$.'

\eject

We will also identify the smallest initial segment of \M\ (i.e., the smallest
non\-empty subset of \M\ closed under successor and less than) with the
standard model \Nat, i.e., we assume that $\Nat \subseteq \M$. We will reserve
the symbol \Nat\ to denote the standard model and $\omega$ to denote the first
infinite ordinal, which is the same as the domain of $\Nat$.

If $\varphi(x)$ is a formula with a free variable we will write $\varphi(\M)$
for the set of elements in \M\ satisfying $\varphi(x)$, i.e.,
\[
\Set{a \in \M{} | {}\M \models \varphi(a)}.
\]

Given $a \in \M$ we define $\initialseg{a}$ to be the initial segment
\[
\Set{ x \in \M | x < a}.
\]

\section{Coding}

We will assume a notion of finite sets and a definable predicate $x \in y$
(i.e., an \LA-formula with two free variables) such that the universal
closures of
\begin{gather*}
  x \in y \imp x < y, \\
  \forall w \bigl(w \in x \ekv w \in y\bigr) \imp x=y, \\
  \exists z \forall y \bigl(y \in z \ekv y=x\bigr),\\
  \exists z \forall w \bigl(w \in z \ekv \bigl[w \in
  x \vee w \in y\bigr]\bigr),\\
  \exists z \forall w \bigl(w \in z \ekv \bigl[w \in
  x \wedge w \in y\bigr]\bigr), \quad \text{and} \\
  \exists z \forall w \bigl(w \in z \ekv \bigl[w \in x \wedge w \notin
  y\bigr]\bigr),
\end{gather*}
are all provable in \PA. The sets $z$ (which all are unique by the second
property) in the last four formulas will be denoted $\set{x}$, $x \cup y$, $x
\cap y$ and $x \setminus y$ respectively. We will also write $\Set{x,y}$,
$\Set{x,y,z}$, \dots\ for $\Set{x} \cup \Set{y}$, $\Set{x} \cup \Set{y} \cup
\Set{z}$, \dots\ respectively. The membership predicate can be defined by
using the exponentiation function, we will not go into the details of this
here; for a good reference see \cite{Hajek:98}.

We also assume a notion of finite sequences, either derived from the notion of
finite sets (see \cite{Hajek:98}) or by the Chinese remainder theorem (see
\cite{Kaye:91*2}). The predicate $(x)_y=z$ is assumed to be such that the
universal closures of the following formulas are provable in \PA,
\begin{gather*}
  \existsone z \, (x)_y=z,\\
  (x)_y \leq x,\\
  \exists y \, (y)_0=x, \quad \text{and}\\
  \exists w \bigl(\forall i \mathord< y \, (x)_i=(w)_i \wedge (w)_y=z\bigr).
\end{gather*}
We define the following provable recursive functions (and constant)
\begin{align*}
  \len(x)\eqdef & (x)_0,\\
  [x]_y \eqdef & (\least z)\bigl(y < \len(x) \wedge (x)_{\Sc(y)} =
  z\bigr) \vee \bigl(y \geq \len(x) \wedge z=0\bigr), \\
  x \conc y \eqdef & (\least z) \len(x)+\len(y)=\len(z) \\
  & \quad {} \wedge \forall i \mathord< \len(x) \, [x]_i=[z]_i
  \wedge \forall i \mathord< \len(y) \, [y]_i=[z]_{\len(x)+i},\\
  [] \eqdef & (\least x)\len(x)=0, \\
  [x] \eqdef & (\least y)\len(y)=1 \wedge [y]_0=x, \\
  [x_0,\ldots,x_{k}] \eqdef & [x_0] \conc \cdots
  \conc [x_{k}],\\
  x \upharpoonright y \eqdef & (\least z) \len(z)=y \wedge \forall i \mathord<
  y \, [z]_i=[x]_i,
\end{align*}
here $(\least x)\varphi(x)$ means `the least $x$ such that $\varphi(x)$.'

\section{Nonstandard languages}

We need a G\"odel numbering for the formulas and terms in the language \LMs.
We might define the G\"odel number for a formula $\varphi$, denoted
$\godel{\varphi}$, to be (a code for) the sequence of the G\"odel numbers of
the symbols in $\varphi$, thus
\[
\godel{\mathstrut\Sc(0)=\v_0}=\Bigl[\godel{\mathstrut\Sc},\godel{\mathstrut(\,},
\godel{\mathstrut0}, \godel{\mathstrut\,)}, \godel{\mathstrut0},
\godel{\mathstrut\v_0}\Bigr].
\]
The exact definition we use for G\"odel numbering is unimportant. But it will
{\em not}, except in some special occasions, be assumed to be defined in this
way; any numbering such that the properties below hold will work.

There are \LA-formulas $\form(x)$, $\sent(x)$, $\term(x)$ and $\clterm(x)$
coding, in \PA, the formulas, sentences, terms and closed terms of \LMs\ 
respectively. Let $\FV$ be the function (defined in \PA) such that
$\FV(\varphi)$ is (an element coding) the set of G\"odel numbers of the free
variables of $\varphi$. The precise construction of the formulas are not
important, but some of the properties of them are. Those are listed below.

The notation $\godel{\Sc(x)}$ will be taken, when appropriate, to mean the
{\em function} which takes a G\"odel number of a term $t$ and returns the
G\"odel number of the term $\Sc(t)$. By `when appropriate' we mean that in
some cases $\godel{\Sc(x)}$ means the G\"odel number of the term $\Sc(x)$, but
we will try to write $\godel{\Sc(\v_i)}$ in that case. Of course this also
applies to for example $\godel{\exists \v_i\, x}$ which is a function taking a
G\"odel number of a formula $\varphi$ and an $i$ and returning the G\"odel
number of the formula $\exists \v_i \,\varphi$. With the assumption that a
G\"odel number is a sequence of (G\"odel numbers of) symbols we get that
\[
\PA \prf \forall x \,
\godel{\mathstrut\Sc(x)}=\Bigl[\godel{\mathstrut\Sc},\godel{\mathstrut(\,}\bigr]\conc
x \conc \bigl[\godel{\mathstrut\,)}\Bigr].
\]

If $x$ is not a G\"odel number, or a G\"odel number of ``wrong type,'' then
the functions can be defined to take the value $0$.

In \PA\ we can define the substitution function that takes the G\"odel number
of a formula/term, a term and a variable and returns the G\"odel number of the
formula/term we get by substituting all occurrences of the given variable with
the given term. We denote this function $x[y/z]$.

That the universal closures of the following formulas are provable in \PA\ 
tells us that all elements that should satisfy $\term(x)$ does so.
\begin{gather*}
  \term(\godel{0}) \wedge \term(\godel{\c_x}) \wedge \term(\godel{\v_i}), \\
  \term(x) \imp \term(\godel{\Sc(x)}), \\
  \term(x) \wedge \term(y) \imp \term(\godel{x+y}), \quad \text{and} \\
  \term(x) \wedge \term(y) \imp \term(\godel{x \cdot y}).
\end{gather*}
The analogous formulas for $\form(x)$ are:
\begin{gather*}
  \term(x) \wedge \term(y) \imp \form(\godel{x=y}), \\
  \form(x) \imp \form(\godel{\neg x}), \\
  \form(x) \wedge \form(y) \imp \form(\godel{x \vee y}) \quad \text{and} \\
  \form(x) \imp \form(\godel{\exists \v_i \, x}).
\end{gather*}
The next properties tell us that nothing other than what is supposed to
satisfies $\term(x)$. This is the inductive property of terms:
\begin{multline*}
  \varphi(\godel{0}) \wedge \forall i \, \varphi(\godel{\v_i}) \wedge \forall
  x \, \varphi(\godel{\c_x}) \wedge \forall x,y
  \bigl[ \term(x) \wedge \term(y) \wedge \varphi(x) \\
  {}\wedge \varphi(y) \imp \varphi(\godel{\Sc(x)}) \wedge \varphi(\godel{x+y})
  \wedge \varphi(\godel{x \cdot y})\bigr] \imp \forall x \bigl(\term(x) \imp
  \varphi(x)\bigr)
\end{multline*}
for all \LMs-formulas $\varphi(x)$. The analogous property for $\form(x)$:
\begin{multline*}
  \forall x,y \bigl(\term(x) \wedge \term(y) \imp \varphi(\godel{x=y})\bigr)
  \wedge \forall x,y,i \bigl[ \form(x) \wedge \form(y) \\
  {}\wedge \varphi(x) \wedge \varphi(y) \imp \varphi(\godel{\neg x}) \wedge
  \varphi(\godel{x \vee y}) \wedge \varphi(\godel{\exists \v_i x})\bigr] \imp
  \forall x \bigl(\form(x) \imp \varphi(x)\bigr)
\end{multline*}
for all \LMs-formulas $\varphi(x)$.

There are also some similar properties for $\FV(x)$.
\begin{gather*}
  \FV(\godel{0})=0 \wedge \forall x \, \FV(\godel{\c_x})= 0, \\
  \mathop{\forall i } \FV(\godel{\v_i})= \Set{\godel{\v_i}}, \\
\begin{split}
  \forall x,y \bigl[ \term(x) \wedge \term(y) \imp \FV(\godel{\Sc(x)})&=\FV(x)
  \wedge \FV(\godel{x+y}) =\FV(\godel{x \cdot y}) \\ & {}=\FV(\godel{x=y})
  =\FV(x) \cup \FV(y)\bigr], \quad \text{and}
\end{split}\\
\begin{split}
  \forall x,y,i \bigl[\form(x) \wedge \form(y) \imp \FV(\godel{\neg x}) =
  \FV(x) &\wedge \FV(\godel{x \vee y}) = \FV(x) \cup \FV(y) \\ &{} \wedge
  \FV(\godel{\exists \v_i x}) = \FV(x) \setminus \Set{\godel{\v_i}} \bigr].
\end{split}
\end{gather*}
Last we have the defining properties of $\clterm(x)$ and $\sent(x)$:
\begin{align*}
  \clterm(x) &\ekv \term(x) \wedge \FV(x)=0 \quad \text{and}\\
  \sent(x) &\ekv \form(x) \wedge \FV(x)=0. \\
\end{align*}

Observe that these properties define the formulas $\term(x)$, $\form(x)$,
$\clterm(x)$ and $\sent(x)$ and the function $\FV(x)$ up to provable
equivalence in \PA.

Two very important properties, which follows from the properties above, of
\PA\ is the following, usually called the `unique readability property.' \PA\ 
proves the following two sentences:
\begin{multline*}
  \forall x \Bigl( \term(x) \imp \bigl[ \existsone i \, x=\godel{\v_i} \xor
  \existsone y \, x=\godel{\c_y} \xor \existsone y \bigl(\term(y) \wedge
  x=\godel{\Sc(y)}\bigr) \\ {}\xor \existsone y,z \bigl(\term(y) \wedge
  \term(z) \wedge x=\godel{y+z}\bigr) \\ {}\xor \existsone y,z \bigl(\term(y)
  \wedge \term(z) \wedge x=\godel{y \cdot z}\bigr)\bigr]\Bigr)
\end{multline*}
and
\begin{multline*}
  \forall x \Bigl(\form(x) \imp \bigl[ \existsone y,z \bigl(\term(y) \wedge
  \term(z) \wedge x=\godel{z=y}\bigr) \\ {} \xor \existsone y \bigl(\form(y)
  \wedge x=\godel{\neg y}\bigr) \xor \existsone y,z \bigl(\form(y) \wedge
  \form(z) \wedge x=\godel{y \vee z}\bigr) \\{}\xor \existsone y,i
  \bigl(\form(y) \wedge x=\godel{\exists \v_i\, y}\bigr)\bigr]\Bigr).
\end{multline*}

These properties make it possible to handle nonstandard languages.  Let \LM\ 
be the nonstandard language which corresponds to \LMs, i.e., the ``terms'' of
\LM\ are all $a \in \M$ such that $\M \models \term(a)$ and the ``formulas'' are
all $a \in \M \models \form(a)$, etc. The unique readability properties give
us the possibility to handle these ``terms'' and ``formulas'' in much the same way
as the standard ones, with the important exception that they need not be
well-founded, e.g., $\neg\neg\ldots\neg 0=1$, where the dots represent a {\em
  nonstandard\/} number of negation signs, is a \LM-formula. Therefore, and
this is very important, we do not have ``external'' induction on \LM-terms and
\LM-formulas.\footnote{We have induction {\em inside} the model, that is what
  the inductive property of $\form(x)$ and $\term(x)$ tells us.}


\section{Partial truth definitions} \label{sec:part.tr.def}

Due to Tarski's theorem on the undefinability of truth we cannot find an
\LMs-formula $\varphi$ such that
\[\PA \prf \varphi(\godel{\psi}) \ekv \psi\]
for all \LMs-sentences $\psi$. What we can do is the following.

\begin{thm}
  There is an \LA-definable function $\val$ such that \PA\ proves
\[
\val(\godel{\,t\,})=t
\]
for all closed \LMs-terms $t$, and
\begin{multline*}
  \forall x \bigl(\clterm(x) \imp
  \val(\godel{\Sc(x)})=\Sc(\val(x))\bigr)\\
  {} \wedge \forall x,y\bigl(\clterm(x) \wedge \clterm(y) \imp
  \val(\godel{x+y})=\val(x)+\val(y)\bigr) \\
  {}\wedge \forall x,y \bigl(\clterm(x) \wedge \clterm(y) \imp
  \val(\godel{x\cdot y})=\val(x)\cdot \val(y)\bigr).
\end{multline*}
\end{thm}

Let $\Delta_k$, $\Sigma_k$ and $\Pi_k$ be defined as usual, e.g., $\Sigma_1$
is the set of all formulas of the form $\exists \bar{x} \, \varphi(\bar{x})$
where $\varphi(\bar{x})$ is in $\Delta_0$.  There are \LA-formulas coding
these sets in \PA\ in the usual sense. Let us write, for example, $x \in
\Delta_0$ for the formula coding $\Delta_0$ applied to the variable $x$.

\begin{thm}
  There are \LA-formulas $\Tr_\Gamma(x)$, where $\Gamma$ is $\Delta_k,
  \Sigma_k$ or $\Pi_k$, such that \PA\ proves
\[
\Tr_\Gamma(\godel{\psi}) \ekv \psi
\]
for all sentences $\psi \in \Gamma$,
\begin{gather*}
  \forall x,y \bigl[\clterm(x) \wedge \clterm(y) \imp
  \bigl(\Tr_{\Delta_0}(\godel{x=y}) \ekv \val(x)=\val(y)\bigr)\bigr],\\
  \forall x \bigl[x \in \Delta_0 \wedge \sent(x) \imp
  \bigl(\Tr_{\Delta_0}(\godel{\neg x}) \ekv \neg
  \Tr_{\Delta_0}(x)\bigr)\bigr], \\
  \forall x,y \bigl[x,y \in \Delta_0 \wedge \sent(x) \wedge \sent(y) \imp
  \bigl(\Tr_{\Delta_0}(\godel{x \vee y}) \ekv \Tr_{\Delta_0}(x)
  \vee \Tr_{\Delta_0}(y)\bigr)\bigr],\\
  \forall x,i,y \bigl[ x \in \Delta_0 \wedge \sent(\godel{\exists \v_i \, x})
  \imp \bigl( \Tr_{\Delta_0}(\godel{\exists \v_i \mathord< \c_y \, x}) \ekv
  \exists z \mathord< y \, \Tr_{\Delta_0}(\godel{x[\c_z/\v_i]})\bigr)\bigr],
\end{gather*}
and also
\begin{gather*}
  \forall x,i \bigl[\godel{\exists \v_i \,x} \in \Gamma \wedge
  \sent(\godel{\exists \v_i \, x}) \imp \bigl(\Tr_\Gamma(\godel{\exists \v_i\,
    x}) \ekv \exists y\,
  \Tr_\Gamma(\godel{x[\c_y/\v_i]})\bigr)\bigr], \quad \text{and} \\
  \forall x,i \bigl[\godel{\forall \v_i\, x} \in \Gamma \wedge
  \sent(\godel{\forall \v_i \, x}) \imp \bigl(\Tr_\Gamma(\godel{\forall \v_i
    \, x}) \ekv \forall y \, \Tr_\Gamma(\godel{x[\c_y/\v_i]})\bigr)\bigr].
\end{gather*}
\end{thm}

For explicit constructions see \cite{Kaye:91*2}.

\section{Recursive saturation and resplendency}

In this section \M\ can be any structure in any recursive language $\L$. Let
$\LMs$ be the language which extends $\L$ with constant symbols naming all
elements in \M.

Given a theory $\Te$, a {\em type in $\Te$} is a countable set of formulas
\[
t(\bar{x}) = \Set{\varphi_i(\bar{x},\bar{a}) | i \in \omega}
\]
such that $\varphi_i(\bar{x},\bar{y})$ are \L-formulas with finitely many free
variables $\bar{x}$ and $\bar{y}$; $\bar{a}$ are parameters from \M; and $\Te +
t(\bar{\c})$, where $\bar{\c}$ are new constant symbols, is consistent. A type
over a model \M\ is a type in the elementary diagram of the model, i.e., in
$\ElDiag(\M)$. A type $t(\bar{x})$ is {\em recursive} if the set
\[
\Set{ \godel{\varphi( \bar{x},\bar{y})} | \text{there exists $\bar{a} \in \M$
    such that $\varphi( \bar{x},\bar{a}) \in t(\bar{x})$}} \subseteq \omega
\]
is recursive. If $t(\bar{x})$ is a type then it is {\em realized\/} in \M\ if
there are elements $\bar{m} \in \M$ such that $\M \models \varphi(\bar{m})$
for all $\varphi (\bar{x}) \in t(\bar{x})$.

\begin{defin}
  \M\ is {\em recursively saturated\/} if all recursive types in \M\ are
  realized.
\end{defin}

The following proposition says that all consistent theories have recursively
saturated models of any cardinality.

\begin{prop}
  For every \M\ there is an elementary extension $\structure{N}\succ \M$ that
  is recursively saturated and such that $\abs{\structure{N}} = \abs{\M}$.
\end{prop}

Now to a slightly different notion, that of resplendency. A $\Sigma_1^1$
formula is a second-order formula of the form $\exists X \,\varphi(X)$ where
$X$ is a set variable and $\varphi(X)$ is a first-order formula in the
language extended with the set variable $X$.

\begin{defin}
  \M\ is {\em resplendent\/} if for all $\Sigma_1^1$ \LMs-sentences $\Phi$,
  such that
\[
\ElDiag(\M) \cup \Set{\Phi}
\]
is consistent, we have $\M \models \Phi$.
\end{defin}

In other words \M\ is resplendent if as many as possible $\Sigma_1^1$
sentences are true and it is recursively saturated if all recursive types are
realized. Observe that both these notions apply to all structures, not only
models of \PA. The next theorem follows from work by Kleene \cite{Kleene:52}.

\begin{thm}
  If \M\ is resplendent then it is recursively saturated.
\end{thm}

There is a converse if the model is countable, this result is due to Barwise
and Schlipf and independently Ressayre.

\begin{thm}[\cite{Barwise.Schlipf:76}]\label{thm:sat.imp.res}
  If \M\ is countable and recursively saturated then it is resplendent.
\end{thm}

There are several model theoretic properties which are $\Sigma_1^1$, making
resplendent models easy to work with.  In this thesis we will express the
consistency of logics, which are definable in some model, by an $\Sigma_1^1$
formula:
\begin{multline*}
  \exists X \Bigl(\text{all axioms are in $X$}{} \wedge {}\text{$X$ is closed
    under the inference rules} \\ {}\wedge {}\exists x \bigl(\text{$x$ is a
    formula}{} \wedge x \notin X\bigr)\Bigr).
\end{multline*}

Hopefully, the reader is now ready to face satisfaction classes.

\chapter{Satisfaction classes}\label{chp:satcl}

When defining satisfaction classes we have two different approaches to choose
from. The historical way (used in \cite{Krajewski:76}) is to look at
nonstandard formulas of \LA\ and define a satisfaction class to be a set of
pairs; the first component being a nonstandard formula of \LA\ and the second
being a code for a sequence of elements in \M. The intention is that the
sequence satisfies the formula, i.e., if we substitute the free variable
$\v_i$ with the $i$th element of the sequence then the result is ``true'' in
\M. The disadvantage of this approach is that we need some machinery to handle
different ways of getting the ``same'' formula. We will give an example; let
$\epsilon_0$ be $\v_0\neq \v_0$ and $\epsilon_{i+1}$ be $\epsilon_i \vee
\epsilon_i$, also let $\epsilon_0'$ be $\v_0\neq 0$ and $\epsilon_{i+1}'$ be
$\epsilon_i' \vee \epsilon_i'$, then we clearly want $\pair{\epsilon_i,[0]}
\in \Sat$ iff $\pair{\epsilon_i',[0]} \in \Sat$.

\conjecture{There is a satisfaction class $\Sat$, in the sense of
  \cite{Krajewski:76} and \cite{Kotlarski.Krajewski.ea:81}, such that
  $\pair{\epsilon_a,[0]} \in \Sat$ and $\pair{\epsilon_a',[0]} \protect\notin
  \Sat$ for some $a \in \M$.}{Remember that $\epsilon_0$ is $\v_0\neq \v_0$ and
  $\epsilon_{i+1}$ is $\epsilon_i \vee \epsilon_i$; $\epsilon_0'$ is $\v_0\neq
  0$ and $\epsilon_{i+1}'$ is $\epsilon_i' \vee \epsilon_i'$. We think that by
  redefining \M-logic to work with pairs of \LA-formulas and elements of \M\ 
  it should be possible to prove the conjecture.}

We think that by redefining \M-logic to work with pairs of \LA-formulas and
elements of \M\ it should be possible to prove the conjecture by reproving the
results in this chapter.

We are going to define a satisfaction class to be a set of \LM-sentences; the
intention is that the sentences are those which are ``true'' in \M. The
disadvantage of this approach is that the name `satisfaction class' seems a
bit awkward, a better name would probably be `truth class', but for historical
reasons we will stick with it.

We will use the two notations $x \in \Sat$ and $\Sat(x)$ to mean the same
thing, i.e., in this case the symbol $\in$ has nothing to do with the coding
of finite sets. We hope this will not confuse the reader, but instead make the
formulas easier to read.

A comment on the word {\em class} should be made here. A class in the model
theory of arithmetic is a subset $C$ of the domain of the model such that for
all $a$, $C \cap \initialseg{a}$ is definable (remember that $
\initialseg{a}=\Set{x \in \M | x < a}$). It is just a mere and unfortunate
coincidence that the word is used in the term `satisfaction class.'

\begin{defin}\label{def:satcl}
  A satisfaction class \Sat\ is an external subset of \M\ satisfying the
  following conditions in \M:
\begin{gather}
\label{s1} x \in \Sat \imp \sent(x),\\
\label{s2} \godel{\c_a=\c_b} \in \Sat \ekv a=b,\\
\label{s3} \godel{t=t} \in \Sat,\\
\label{s4} \godel{t=r} \in \Sat \imp \godel{r=t} \in \Sat,\\
\label{s5} \godel{t=r} \in \Sat \wedge \godel{r=s} \in \Sat \imp \godel{t=s} \in \Sat,\\
\label{s6} \godel{\Sc(t)=\c_a} \in \Sat \ekv \exists
x\bigl(\godel{t=\c_x}
\in \Sat \wedge \Sc(x)=a\bigr),\\
\label{s7} \godel{t+r=\c_a} \in \Sat \ekv \exists x,y
\bigl(\godel{t=\c_x} \in \Sat
\wedge \godel{r=\c_y} \in \Sat \wedge x+y=a\bigr),\\
\label{s8} \godel{t\cdot r=\c_a} \in \Sat\ekv \exists
x,y\bigl(\godel{t=\c_x}
\in \Sat \wedge \godel{r=\c_y} \in \Sat \wedge x\cdot y=a\bigr), \\
\label{s9} \godel{t=r} \in \Sat \ekv \exists x \bigl(
\godel{t=\c_x} \in \Sat
\wedge \godel{r=\c_x} \in \Sat\bigr), \\
\label{s10} \godel{\neg \varphi} \in \Sat \ekv \varphi \notin \Sat, \\
\label{s11} \godel{\varphi \vee \psi} \in \Sat \ekv \bigl(\varphi
\in
\Sat \vee \psi \in \Sat\bigr) \quad \text{and}\\
\label{s12} \godel{\exists \v_i \, \gamma} \in \Sat \ekv \exists x
\bigl(\godel{\gamma[\c_x/\v_i]} \in \Sat\bigr),
\end{gather}
for all closed \LM-terms $t$, $r$ and $s$ and all \LM-sentences $\varphi$,
$\psi$ and $\exists \v_i \gamma$.
\end{defin}

The clauses (\ref{s2})-(\ref{s9}) take care of the atomic formulas and we
could replace them with the single clause
\[
t=r \in \Sat \ekv \M \models \val(t)=\val(r),
\]
but that would yield a stronger notion. We will look more closely at this
later in Chapter~\ref{chp:stsatcl}.

Observe that there is a first-order formula $\satcl(X)$ in the language $\LA
\cup \Set{X}$ such that for any subset $A$ of \M\, $A$ is a satisfaction class
iff $\M \models \satcl(A)$. The formula $\satcl(X)$ is the conjunction of the
universal closures (with the universal quantifiers restricted to closed terms
and sentences) of (\ref{s1})-(\ref{s12}), which all are first-order
properties.

\begin{prop}
  Let \Sat\ be a satisfaction class, then the following is true in the
  structure $\pair{\M,\Sat}$:
\begin{align*}
  \godel{\varphi \wedge \psi} \in \Sat &\ekv \varphi \in
  \Sat \wedge \psi \in \Sat, \\
  \godel{\varphi \imp \psi} \in \Sat &\ekv \bigl(\varphi \in
  \Sat \imp \psi \in \Sat\bigr), \\
  \godel{\varphi \ekv \psi} \in \Sat &\ekv \bigl(\varphi \in
  \Sat \ekv \psi \in \Sat\bigr), \quad \text{and} \\
  \godel{\forall \v_i \,\gamma} \in \Sat &\ekv \forall x
  \bigl(\godel{\gamma[\c_x/\v_i]} \in \Sat\bigr),
\end{align*}
for all \LM-sentences $\varphi$, $\psi$ and $\forall \v_i \,\gamma$.
\end{prop}

\begin{prop}\label{prop:thm.in.satcl}
  If \Sat\ is a satisfaction class, $\varphi$ an \LMs-sentence and $\M \models
  \varphi$ then $\varphi \in \Sat$, i.e., $\ElDiag(\M) \subseteq \Sat$.
\end{prop}
\begin{proof}
  We first prove the statement for sentences $\varphi$ of the form $t=\c_a$
  for some closed term $t$ and some $a \in \M$, by induction on the
  construction of $t$.
  
  If $t$ is $\c_b$ then $\M \models \c_b = \c_a$ so $\c_b=\c_a \in \Sat$ by
  (\ref{s2}).
  
  Suppose $t$ is $\Sc(t')$ for some closed term $t'$; then $\M \models
  t'=\c_b$ for some $b \in \M$ and by the induction hypothesis we have
  $t'=\c_b \in \Sat$. Also $\Sc^\M(b)=a$ so $\Sc(t')=\c_a \in \Sat$ by
  (\ref{s6}).
  
  Suppose $t$ is $t_1+t_2$, then there are $b$, $b' \in \M$ such that $\M
  \models t_1=\c_b$ and $\M \models t_2=\c_{b'}$. By the induction hypothesis
  we have $t_1=\c_b \in \Sat$ and $t_2=\c_{b'} \in \Sat$. We know that $b +^\M
  b'=a$, so $t_1+t_2 = \c_a \in \Sat$ by (\ref{s7}). The case when $t$ is $t_1
  \cdot t_2$ is treated in a similar way.
  
  Suppose now $\M \models t_1=t_2$ then there is a $b \in \M$ such that $\M
  \models t_1=b$ and $\M \models t_2=b$, so $t_1 = \c_b \in \Sat$ and $t_2 =
  \c_b \in \Sat$ which by (\ref{s4}) and (\ref{s5}) implies $t_1=t_2 \in \Sat$.
  
  If $\M \models t_1 \neq t_2$ then there are $b_1$, $b_2 \in \M$ such that
\[
\M \models t_1=b_1 \wedge t_2=b_2 \wedge b_1 \neq b_2,
\]
so $t_1=\c_{b_1} \in \Sat$ and $t_2=\c_{b_2} \in \Sat$. By (\ref{s2}) we also
have $\c_{b_1} \neq \c_{b_2} \in \Sat$ and by (\ref{s4}) and (\ref{s5}) we
have $t_1 \neq t_2 \in \Sat$.

We have proved that $\Diag(\M) \subseteq \Sat$. We should prove that if
$\varphi$ is an \LMs-sentence then
\[
\M \models \varphi \quad \Leftrightarrow \quad \varphi \in \Sat
\]
by induction on $\varphi$.
The case when $\varphi$ is atomic is proven. The induction step is easy and
left to the reader.
\end{proof}

\begin{rem}
  By Tarski's theorem of the undefinability of truth a satisfaction class
  \Sat\ is not definable. Moreover $\Sat \cap \initialseg{a}$ for $a > \omega$
  is not definable, since this would also yield a truth definition, therefore
  \Sat\ is not a class.
\end{rem}

Given a set of \LM-sentences $X$, define the binary relation $\sim_X$ (or just
$\sim$ if $X$ is understood from the context) on the set of closed \LM-terms,
$\clterm(\M)$, as follows
\[t_1 \sim_X t_2 \quad \text{iff}\quad t_1 = t_2 \in X.\]
If $\sim_X$ is an equivalence relation let $\clterm(\M) / \mathord\sim_X$ be the
set of equivalence classes $\overline{t}$ of $\clterm(\M)$ and define the functions
$\Sc^{\M_X}$, $+^{\M_X}$ and $\cdot^{\M_X}$ by
\begin{align*}
  \Sc^{\M_X}(\,\overline{t}\,)&\eqdef\overline{\Sc(t)}, \\
  \overline{t} +^{\M_X} \overline{r}&\eqdef\overline{t+r} \quad \text{and} \\
  \overline{t} \cdot^{\M_X} \overline{r}&\eqdef\overline{t\cdot r}
\end{align*}
if they all are well-defined. This defines an \LA-structure
\[
(\clterm(\M)/ \mathord\sim_X, \Sc^{\M_X},+^{\M_X},\cdot^{\M_X},\overline{0}).
\]
Let $\M_X$ denote this structure, when it is well-defined.

The following equivalent definition of a satisfaction class may be well worth
notice; it will play the main role of Chapter~\ref{chp:wesatcl}.

\begin{prop}\label{prop:alt.satcl.def}
  \Sat\ is a satisfaction class iff it is a unary predicate on \M\ such that
  $\sim_\Sat$ is an equivalence relation, $\M_\Sat$ is well-defined, the
  canonical map
\[f: \M \to \M_\Sat, \quad a \mapsto\overline{\c_a}\]
is an isomorphism and
\begin{gather}
  \tag{\ref{s1}} x \in \Sat \imp \sent(x),\\
  \tag{\ref{s10}}\godel{\neg \varphi} \in \Sat \ekv \varphi \notin \Sat, \\
  \tag{\ref{s11}}\godel{\varphi \vee \psi} \in \Sat \ekv \bigl(\varphi \in
  \Sat \vee \psi \in \Sat\bigr), \quad\text{and} \\
  \tag{\ref{s12}}\godel{\exists \v_i \,\gamma} \in \Sat \ekv \exists x
  \bigl(\godel{\gamma[\c_x/\v_i]} \in \Sat\bigr)
\end{gather}
holds for all \LM-sentences $\varphi$, $\psi$ and $\exists \v_i\, \gamma$.
\end{prop}
\begin{proof}
  Assume \Sat\ is a satisfaction class. By (\ref{s3})-(\ref{s5}) it is clear
  that $\sim_\Sat$ is an equivalence relation and by (\ref{s6})-(\ref{s8})
  $\M_\Sat$ is well-defined. It is also clear that the mapping $f$ is a
  bijection since if $a \neq b$ then $\c_a = \c_b \notin \Sat$ so $f(a) \neq
  f(b)$ and if $t \in \clterm(\M)$ then $t=t \in \Sat$, so by (\ref{s9}) there
  exists $m \in \M$ such that $t=\c_m \in \Sat$, i.e., $f(m)=\overline{t}$. We
  also have
\begin{gather*}
  f(\Sc^\M(a))=\overline{\c_{\Sc^\M(a)}}=\overline{\Sc(\c_a)}=
  \Sc^{\M_X}(\overline{\c_a})= \Sc^{\M_X}(f(a)), \\
  f(a+^\M b)=\overline{\c_{a+^\M b}}=\overline{\c_a+\c_b}=
  \overline{\c_a}+^{\M_X} \overline{c_b}= f(a)+^{\M_X}f(b) \quad\text{and} \\
  f(a\cdot^\M b)=\overline{\c_{a\cdot^\M b}}= \overline{\c_a\cdot \c_b}=
  \overline{\c_a}\cdot^{\M_X} \overline{c_b}=f(a)\cdot^{\M_X} f(b),
\end{gather*}
so $f$ is really an isomorphism.

Assume now that \Sat\ satisfies the conditions above.

We have to prove conditions (\ref{s2}) and (\ref{s6})-(\ref{s9}) in the
definition of satisfaction classes. The other follows immediately.

(\ref{s2}) It follows from the injectivity of $f$, since if $a \neq b \in \M$
then $f(a) \neq f(b)$ so $\overline{\c_a} \neq \overline{\c_b}$, i.e., $\c_a
\neq \c_b \in \Sat$.

(\ref{s6}) Assume $\Sc(t)=\c_a \in \Sat$, we have to find $m$ such that
$t=\c_m \in \Sat$ and $\Sc^\M (m)=a$. We have
\[
a = f^{-1}(\overline{\c_a}) = f^{-1}(\overline{\Sc(t)}) =
f^{-1}(\Sc^{\M_X}(\,\overline{t}\,)) = \Sc^\M (f^{-1}(\,\overline{t}\,)),
\]
so $m=f^{-1}(\,\overline{t}\,)$ works since $\overline{t}=f(m)$ implies
$t=\c_m \in \Sat$. On the other hand assume that there exists $m$ such that
$t=\c_m \in \Sat$ and $\Sc^\M(m)=a$ then
\[
\overline{\Sc(t)}=\Sc^{\M_X}(\,\overline{t}\,)=\Sc^{\M_X}(\overline{\c_m})=
\Sc^{\M_X}(f(m)) = f(\Sc^\M(m)) = \overline{\c_a}.
\]

(\ref{s7}) Assume $t+r=\c_a \in \Sat$ and let $m=f^{-1}(\,\overline{t}\,)$ and
$m'=f^{-1}(\overline{r})$ then $a=m+^\M m'$ and $t=\c_m \in \Sat$ and
$r=\c_{m'} \in \Sat$. If $m+^\M m'=a$, $t=\c_m \in \Sat$ and $r=\c_{m'} \in
\Sat$ then $f^{-1}(\,\overline{t}\,)=m$ and $f^{-1}(\overline{r})=m'$ so
\[
a=f^{-1}(\overline{t}+^{\M_X} \overline{r}) = f^{-1}(\overline{t+r}).
\]
In other words $f(a)=\overline{t+r}$, i.e., $t+r=\c_a \in \Sat$.  (\ref{s8})
Similar as above.

(\ref{s9}) Assume $t=r \in \Sat$ then
$f^{-1}(\,\overline{t}\,)=f^{-1}(\overline{r})=m$ so $t=\c_m$, $r=\c_m \in
\Sat$ and if $t=\c_m$, $r=\c_m \in \Sat$ then
$\overline{t}=\overline{\c_m}=\overline{r}$, so $t=r \in \Sat$.
\end{proof}

\question{Are satisfaction classes built up by two ``parts;'' one with closed
  terms and equality and one with ``the rest,'' in the following sense: Given
  a relation $\sim$ satisfying the conditions in
  Proposition~\ref{prop:alt.satcl.def} is there a satisfaction class \Sat\ 
  such that the relation $\sim_\Sat$ coincides with $\sim$?}{A partial answer
  is given in Proposition~\ref{prop:partial.answer}. The general question
  seems to be hard, since it is a question of whether a set of equalities is
  consistent or not in \M-logic. See also Question~\ref{que:chapter.five}.}

\label{que:satcl.two.parts}

A partial answer to this question is given in
Proposition~\ref{prop:partial.answer}.

\section{Inductive partial satisfaction classes}

This section is included as an introduction to satisfaction classes; the
results will not be used later and may therefore be skipped. For simplicity
the satisfaction classes in this section will all include the set
\[
\Set{t=r |\M \models \val(t)=\val(r)}.
\]
Therefore, in this section, $\Sat \subseteq \M$ is a satisfaction class iff
\[
\pair{\M,\Sat} \models \forall x \bigl(x \in \Sat \imp \sent(x)\bigr) \wedge
\forall x \Psi(\Sat,x),\] where $\Psi(\Sat,x)$ is the formula
\begin{multline*}
  \forall y,z \bigl[\clterm(y) \wedge \clterm(z) \wedge x=\godel{y=z} \imp
  \bigl(x \in \Sat \ekv
  \val(y)=\val(z)\bigr)\bigr] \\
  {} \wedge \forall y \bigl[ \sent(y) \wedge x = \godel{\neg y} \imp
  \bigl(x \in \Sat \ekv y \notin \Sat\bigr)\bigr]\\
  {} \wedge \forall y,z \bigl[\sent(y) \wedge \sent(z) \wedge x = \godel{y
    \vee z} \imp \bigl(x \in \Sat \ekv
  y \in \Sat \vee z \in \Sat\bigl)\bigr] \\
  {} \wedge \forall y, i \bigl[\sent(\exists \v_i \, y) \wedge
  x=\godel{\exists \v_i \, y} \imp \bigl(x \in \Sat \ekv \exists z\,
  \godel{y[\c_z/\v_i]} \in \Sat\bigr)\bigr].
\end{multline*}

We will assume that if $\psi$ is a subformula or subterm of $\varphi$ then
$\godel{\psi} \leq \godel{\varphi}$; under the assumption that the G\"odel
number of a formula is the sequence of G\"odel numbers of the symbols in the
formula this is true.

\begin{defin}
  A set $\Sat \subseteq \M$ is a {\em partial satisfaction class} if there
  exists $c \in \M \setminus \omega$ such that
\[
\pair{\M,\Sat} \models \forall x \bigl(x \in \Sat \imp \sent(x)\bigr) \wedge
\forall \varphi \mathord< c \, \forall x \,
\Psi(\Sat,\godel{\varphi[x/\v]}).\footnote{See
  Definition~\ref{def:substitution} for the definition of $\varphi[a/\v]$.}
\]
A partial satisfaction class is {\em inductive} if full induction in the
language $\LA \cup \Set{\Sat}$ holds, i.e., if
\[
\pair{\M, \Sat} \models \varphi(0) \wedge \forall x \bigl(\varphi(x) \imp
\varphi(\Sc(x))\bigr) \imp \forall x\, \varphi(x)
\]
for every $\LMs \cup \Set{\Sat}$-formula $\varphi(x)$.
\end{defin}

\begin{prop}
  If $\varphi$ is an \LMs-sentence and $\Sat$ a partial satisfaction class
  then
\[ 
\varphi \in \Sat \quad \text{iff} \quad \M \models \varphi.
\]
\end{prop}
\begin{proof}
  For atomic sentences $\varphi$ the proposition follows trivial from the
  definition of partial satisfaction classes. The rest is an easy induction on
  the construction of $\varphi$.
\end{proof}

\begin{thm}\label{thm:rec.imp.par}
  Every countable recursively saturated model \M\ of \PA\ admits an inductive
  partial satisfaction class.
\end{thm}
\begin{proof}
  Since \M\ is countable and recursively saturated it is resplendent.
  Therefore we only need to show that the theory
\begin{multline*}
  \ElDiag(\M) \\
{}+ \Set{\varphi(0) \wedge \forall x \bigl(\varphi(x) \imp %
\varphi(\Sc(x))\bigr) \imp \forall x \, \varphi(x)  |  \text{%
  $\varphi(x)$ an $\LMs \cup \Set{\Sat}$-formula}} \\
{}+  \forall x \bigl(x \in \Sat \imp \sent(x)\bigr) + \Set{\forall %
x \, \Psi(\Sat,\varphi[x/\v])  |  \text{$\varphi$ an %
  \LA-formula}}
\end{multline*}
is consistent. If this theory is consistent then clearly any $\Sat$ satisfying
it will be an inductive partial satisfaction class by a simple overspill
argument since
\[
\forall n \mathord\in \omega \, \pair{\M,\Sat} \models \forall x \mathord< n
\, \forall x\, \Psi(\Sat,\varphi[x/\v]).
\]

To prove the consistency take a finite subset of the theory, it will at most
involve a finite number of standard formulas $\varphi$ in the scheme
\[
\forall x \, \Psi (\Sat,\varphi[x/\v]).
\]
Since it is a finite number they will all be $\Sigma_k$ for some $k \in
\omega$, therefore we can define
\[
\Sat \eqdef \Set{a \in \M  |  \M \models \Tr_{\Sigma_k}(a)}.  
\qedhere
\]
\end{proof}

\begin{thm}\label{thm:par.imp.rec}
  If $\M \models \PA$ is a nonstandard model admitting an inductive partial
  satisfaction class $\Sat$ then \M\ is recursively saturated.
\end{thm}
\begin{proof}
  Suppose $p(x)$ is a recursive type in \M\ and $\bar{a} \in \M$ are the
  parameters of the type. The type, being recursive, is coded by some $b \in
  \M$. Since
\[
\forall k \mathord\in \omega \, \pair{\M,\Sat} \models \exists x \forall
\varphi \mathord< k \bigl(\varphi \in b \imp
\godel{\varphi\bigl[[x,\bar{a}]/\v\bigr]} \in \Sat\bigr)
\]
we can use overspill (it is here we are using that $\Sat$ is inductive) and
get
\[
\pair{\M,\Sat} \models \exists x \forall \varphi \mathord< c \bigl(\varphi \in
b \imp \godel{\varphi\bigl[[x,\bar{a}]/\v\bigr]} \in \Sat\bigr)
\]
for some $c \in \M \setminus \omega$. This gives us an $x \in \M$ realizing
$p(x)$.
\end{proof}

\begin{rem}
  We are going to do the same sort of argument when $\Sat$ does not satisfy
  induction. The argument there involves more work.
\end{rem}

Combining these two theorems we get the following.

\begin{thm}
  A nonstandard countable model of \PA\ is recursively saturated iff it admits
  an inductive partial satisfaction class.
\end{thm}

\section{Construction of satisfaction classes}
\label{sec:constructionofsatcl}

We will prove the following theorem.

\begin{thm}[\cite{Kotlarski.Krajewski.ea:81}]
  If \M\ is countable and recursively saturated then \M\ admits a satisfaction
  class.
\end{thm}

To prove it we will define a logic, called \M-logic, which we will prove to be
consistent (if \M\ is recursively saturated). Then we construct a maximally
consistent set of sentences in this logic, which turns out to be a
satisfaction class. The idea to use \M-logic is due to Jeff Paris.

\subsection{\M-logic}

As formulas in this logic we will consider all \LM-{\em sentences}. We will
consider a formal deduction system that derives actually finite sets of
formulas, usually denoted by upper case Greek letters. The intention is that
the set $\Gamma$ should be read as $\Vee \Gamma$, i.e., the disjunction of the
sentences in $\Gamma$. The symbol $\varphi$ will denote the singleton set
$\Set{\varphi}$ and $\Gamma, \varphi$ will denote $\Gamma \cup \Set{\varphi}$.
The axioms for the deductive system are:
\begin{gather}
  \tag{Axiom1}\label{Axiom:1}\varphi, \neg \varphi\\
  \tag{Axiom2}\label{Axiom:2}\c_a \neq \c_b \quad \text{if  $a \neq b$}\\
  \tag{Axiom3}\label{Axiom:3}t=t \\
  \tag{Axiom4}\label{Axiom:4}t\neq r,r=t\\
  \tag{Axiom5}\label{Axiom:5}t\neq r, r\neq s, t=s \\
  \tag{Axiom6}\label{Axiom:6}t\neq r,\Sc(t)=\Sc(r) \\
  \tag{Axiom7}\label{Axiom:7}t\neq t',r\neq r',t+r=t'+r' \\
  \tag{Axiom8}\label{Axiom:8}t\neq t',r\neq r',t\cdot r=t' \cdot r' \\
  \tag{Axiom9}\label{Axiom:9} \Sc(\c_a) = \c_{\Sc^\M(a)} \\
  \tag{Axiom10}\label{Axiom:10} \c_a + \c_b = \c_{a +^\M b} \\
  \tag{Axiom11}\label{Axiom:11} \c_a \cdot \c_b = \c_{a \cdot^\M b} \\
  \tag{Axiom12}\label{Axiom:12}\exists \v_0(t=\v_0)
\end{gather}
where $\varphi$ is an arbitrary \LM-sentence and $t$, $r$ and $s$ are
arbitrary closed \LM-terms.  The inference rules are the following:
\begin{equation}
  \tag{Weak}\label{weakening-rule}
  \inference{\Gamma,\varphi}{\Gamma}
\end{equation}
\begin{equation}
  \tag{$\vee$I1}\label{v1-rule}
  \inference{\Gamma,\varphi \vee \psi}{\Gamma,\varphi}
\end{equation}
\begin{equation}
  \tag{$\vee$I2}\label{v2-rule}
  \inference{\Gamma,\varphi \vee\psi}{\Gamma,\psi}
\end{equation}
\begin{equation}
  \tag{$\vee$I3}\label{v3-rule}
  \inference{\Gamma,\neg (\varphi \vee\psi)}{\Gamma,\neg \varphi &
  \Gamma,\neg \psi}
\end{equation}
\begin{equation}
  \tag{$\neg$I}\label{neg-rule}
  \inference{\Gamma,\neg\neg\varphi}{\Gamma,\varphi}
\end{equation}
\begin{equation}
  \tag{Cut}\label{cut-rule}
  \inference{\Gamma}{\Gamma, \varphi & \Gamma, \neg\varphi}
\end{equation}
\begin{equation}
  \tag{$\exists$I}\label{e-rule}
  \inference{\Gamma, \exists \v_i\,\varphi}{\Gamma,\varphi[\c_a/\v_i]}
\end{equation}
\begin{equation}
  \tag{\M-rule}\label{M-rule}
  \inference{\Gamma, \neg \exists \v_i\,\varphi}{\inferforall{a \in
  \M}{\Gamma, \neg \varphi
  [\c_a/\v_i]}}
\end{equation}
where $\Gamma$ is an arbitrary finite set of \LM-sentences, $\varphi$ is an
arbitrary \LM-sentence and the premises in \ref{M-rule} means $\Gamma, \neg
\varphi[\c_a/\v_i]$ for all $a \in \M$.

It might be worth noticing that it is only \ref{Axiom:2}, \ref{Axiom:9},
\ref{Axiom:10}, \ref{Axiom:11} and \ref{M-rule} that depends on the model \M.

At a first glance \ref{Axiom:12} seems to be an ugly duckling, but it turns
out to be of great importance. We will discuss its importance in Chapter
\ref{chp:wesatcl}.

Observe the connection with $\omega$-logic: the logic we get by replacing
\ref{M-rule} with
\[
\inference{\Gamma,\forall \v_i\, \varphi}{\inferforall{a \in
    \M}{\Gamma,\varphi[\c_a/\v_i]}}
\]
is essentially equivalent to \M-logic.

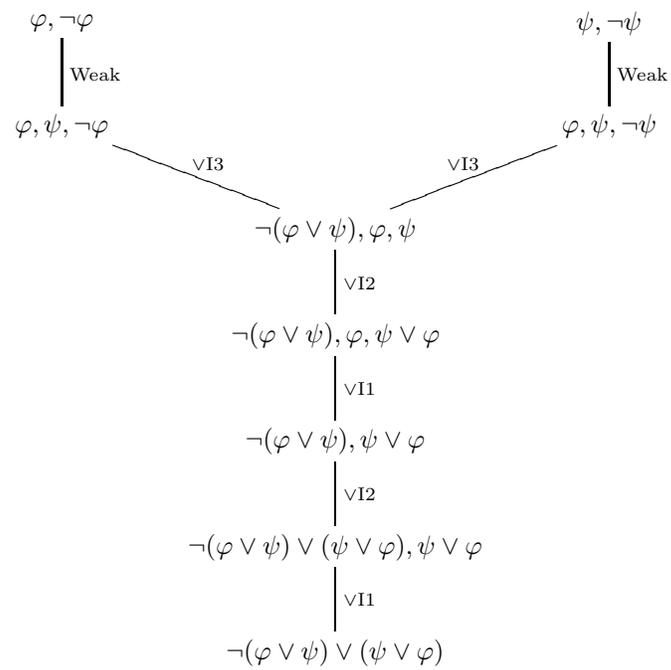
\begin{figure}
\[ \xymatrix{
  {\varphi,\neg\varphi}\ar@{-}[d]^{\text{\ref{weakening-rule}}} & &
  {\psi,\neg\psi}\ar@{-}^{\text{\ref{weakening-rule}}}[d]   \\
  {\varphi,\psi,\neg\varphi}\ar@{-}^{\text{\ref{v3-rule}}}[rd] & &
  {\varphi,\psi,\neg\psi}\ar@{-}_{\text{\ref{v3-rule}}}[ld] \\
  &{\neg(\varphi\vee \psi),\varphi,\psi}\ar@{-}^{\text{\ref{v2-rule}}}[d]&\\
  &{\neg(\varphi\vee \psi),\varphi,\psi \vee \varphi}\ar@{-}^{\text{\ref{v1-rule}}}[d]&\\
  &{\neg (\varphi \vee \psi),\psi \vee \varphi}\ar@{-}^{\text{\ref{v2-rule}}}[d]&\\
  &{\neg (\varphi \vee \psi)\vee (\psi \vee \varphi),\psi \vee
    \varphi}\ar@{-}^{\text{\ref{v1-rule}}}[d]  & \\
  &{\neg (\varphi \vee \psi) \vee (\psi \vee \varphi)} & } \]
\caption{The tree form of the proof of $(\varphi \vee \psi) \imp
  (\psi \vee \varphi)$.} \label{fig:tree}
\end{figure}

The definition of a proof follows the usual one. A proof is a tree where the
nodes are finite sets of sentences, the leaves are axioms and the root is the
conclusion and each edge from one node to another follows one of the inference
rules. For example;
\[
\infer[\text{\ref{v1-rule}}]{\neg (\varphi \vee \psi) \vee (\psi \vee
  \varphi)}{ \infer[\text{\ref{v2-rule}}]{\neg (\varphi \vee \psi) \vee (\psi
    \vee \varphi), \psi \vee \varphi}{ \infer[\text{\ref{v1-rule}}]{\neg
      (\varphi \vee \psi), \psi \vee \varphi}{
      \infer[\text{\ref{v2-rule}}]{\neg (\varphi \vee \psi), \varphi, \psi
        \vee \varphi}{ \infer[\text{\ref{v3-rule}}]{\neg (\varphi \vee \psi),
          \varphi, \psi}{ \infer[\text{\ref{weakening-rule}}]{\varphi, \psi,
            \neg \varphi}{\infer{\varphi, \neg \varphi}{\text{\ref{Axiom:2}}}}
          & \infer[\text{\ref{weakening-rule}}]{ \varphi,\psi, \neg
            \psi}{\infer{\psi, \neg \psi}{\text{\ref{Axiom:2}}}} }}}}}
\]
is a proof of the commutativity of $\vee$, i.e., of $(\varphi \vee \psi) \imp
(\psi \vee \varphi)$. Figure~\ref{fig:tree} shows the corresponding tree.

The notation $\mprf[p] \Gamma$ will mean that $p$ is a proof in \M-logic with
conclusion $\Gamma$ and $\mprf \Gamma$ means that there is a proof $p$ such
that $\mprf[p]\Gamma$. Please observe that the proofs are external objects and
might {\em not\/} be definable in \M, they are either finite or infinite
trees.

If $\Lambda$ is a (external, finite or infinite) set of \LM-sentences then
$\Lambda \mprf[p] \Gamma$ means that $p$ is a proof of $\Gamma$ in \M-logic
with the added axioms
\begin{equation}
\tag{Axiom$\Lambda$}\label{Axiom:L} \varphi \hspace{5mm} \text{if
} \varphi \in \Lambda.
\end{equation}
$\Lambda \mprf \Gamma$ means that there is a proof $p$ such that $\Lambda
\mprf[p] \Gamma$.

If $\Lambda \mprf[p] \Gamma$ then $\abs{p}$ denotes the height of $p$ which is
the height of $p$ looked on as a tree, or in other words;
\begin{defin}
  The {\em height}, $\abs{p}$, of a proof $p$ is defined as follows; if $p$ is
  an axiom then $\abs{p}=0$, if the last inference rule in $p$ is of the form
\[
\inference{\Gamma}{\deduce{\Gamma_1}{p_1} & \deduce{\Gamma_2}{p_2} & \ldots &
  \deduce{\Gamma_i}{p_i} & \ldots }
\] (there are either
one, two or infinitely many premises) then \[\abs{p}=\sup_i (\abs{p_i}+1).\]
\end{defin}

\begin{sidewaysfigure}
  \newcommand{\ab}{\c_{a_0}=\c_{a_0}} \newcommand{\abb}{\c_{a_1}=\c_{a_1}}
  \newcommand{\abbb}{\c_{a_2}=\c_{a_2}} \newcommand{\nab}{\neg
    \c_{a_0}\neq\c_{a_0}} \newcommand{\nabb}{\neg \c_{a_1}\neq\c_{a_1}}
  \newcommand{\nabbb}{\neg \c_{a_2}\neq\c_{a_2}}
\begin{equation*}
\infer{\neg\exists\v_1(\v_1 \neq \v_1)}
{\infer{\nab}{\infer{\ab}{\text{\ref{Axiom:3}}}} &
 \infer{\nabb}{
  \infer{\abb}{
    \infer{\abb,\varphi}{\infer{\abb}{\text{\ref{Axiom:3}}}} &
    \infer{\abb,\neg \varphi}{\infer{\abb}{\text{\ref{Axiom:3}}}}}} &
 \infer{\nabbb}{
   \infer{\abbb}{
     \infer{\abbb,\varphi}{
       \infer{\abbb}{
         \infer{\abbb,\varphi}{\infer{\abbb}{\text{\ref{Axiom:3}}}} &
         \infer{\abbb,\neg \varphi}{\infer{\abbb}{\text{\ref{Axiom:3}}}}}} &
     \infer{\abbb,\neg \varphi}{\infer{\abbb}{\text{\ref{Axiom:3}}}}}} &
     \dots}
\end{equation*}
\caption{An example of a proof of height $\omega$. Here
   $\{a_i\}_{i \in \omega}$ is an enumeration of the domain of \M\ and
  $\varphi$ is, for example, $\c_{a_0}=\c_{a_0}$.}
\label{fig:ett}
\end{sidewaysfigure}

The height $\abs{p}$ is an ordinal number in the ``real world.''  By the
\M-rule $\abs{p}$ might be infinite. For example Figure~\ref{fig:ett} shows a
proof of height $\omega$. As we will see later, the question of whether there
are sets $\Gamma$ provable in \M-logic but not provable by a proof of finite
height is equivalent to the question if \M-logic is consistent.

If $\alpha$ is an ordinal number then $\Lambda \mprf[\alpha] \Gamma$ means
that there is a proof $p$ such that $\abs{p}<\alpha$ and $\Lambda \mprf[p]
\Gamma$, e.g., $\mprf[\omega] \Gamma$ means that there is a finite height
proof of $\Gamma$.

As an alternative we may define
\[
\mprf[\alpha] \Gamma \quad \text{iff} \quad \Gamma \in I_\alpha
\]
where $I_1$ is the set of axioms of \M-logic,
\[
I_{\alpha +1}= \Set{ \Gamma | \text{$\Gamma$ is the result of applying one of
    the inference rules to sets in $I_\alpha$}}
\]
and
\[
I_\lambda = \bigcup_{\alpha < \lambda} I_\alpha
\]
for limit ordinals $\lambda$. Then we could have defined $\mprf \Gamma$ to
hold if there exists an ordinal $\alpha$ such that $\mprf[\alpha] \Gamma$.
This alternative definition avoids mentioning trees, simplifying things.

\begin{lem}\label{lem:deductive.or.neg}
  Let $\Lambda$ be an arbitrary set of \LM-sentences and $\Delta$ and $\Gamma$
  be finite sets of \LM-sentences. Then
\begin{enumerate}
\item $\Lambda, \Delta \mprf \Gamma \quad \text{iff} \quad \Lambda \mprf \neg
  \Delta, \Gamma \qquad \text{where $\neg \Delta = \Set{\neg \varphi
       | \varphi \in \Delta}$}$.
  
\item $\Lambda \mprf \Gamma, \varphi \vee \psi \quad \text{iff} \quad \Lambda
  \mprf \Gamma,\varphi,\psi$.
  
\item $\Lambda \mprf \Gamma,\varphi \quad\text{iff}\quad\Lambda \mprf \Gamma,
  \neg \neg \varphi$.
\end{enumerate}
\end{lem}
\begin{proof}
  1. If $\Delta=\emptyset$, the lemma is trivial. Suppose
  $\Delta=\Set{\varphi}$, then if $\Lambda \mprf \neg \varphi, \Gamma$ it is
  clear that $\Lambda, \varphi \mprf \neg \varphi, \Gamma$ and $\Lambda,
  \varphi \mprf \varphi, \Gamma$ (by \ref{Axiom:L} and \ref{weakening-rule})
  so by \ref{cut-rule} we get $\Lambda, \varphi \mprf \Gamma$.
  
  On the other hand assume that
\[
\Lambda, \varphi \mprf[p] \Gamma.
\]
We will do induction on $\abs{p}$; if $\abs{p} =0$ then $\Gamma$ is an axiom,
if $\Gamma \neq \Set{\varphi}$ it is clear that $\Lambda \mprf \Gamma,\neg
\varphi$ by \ref{weakening-rule} and if $\Gamma=\Set{\varphi}$ then $\Gamma,
\neg \varphi$ is \ref{Axiom:1}.

If $\abs{p}\geq 1$ then the last inference in $p$ is of the form
\[\inference{\Gamma}{\Gamma_0 & \Gamma_1 & \ldots} \]
and by the induction hypothesis we know that $\Lambda \mprf \Gamma_i,\neg
\varphi$, so by the inference
\[\inference{\Gamma,\neg \varphi}{\Gamma_0, \neg \varphi & \Gamma_1,
  \neg \varphi & \ldots} \]we get $\Lambda \mprf \Gamma,\neg \varphi$.

To prove the statement in the general case when $\abs{\Delta} > 1$ we iterate,
moving one formula at a time by induction on $\abs{\Delta}$.

2. Suppose $\Lambda \mprf \Gamma,\varphi,\psi$, then by \ref{v1-rule} $\Lambda
\mprf \Gamma,\varphi \vee \psi,\psi$ and by \ref{v2-rule} $\Lambda \mprf
\Gamma,\varphi \vee \psi$.

If $\Lambda \mprf \Gamma,\varphi \vee \psi$ then by \ref{Axiom:1}
\[
\Lambda \mprf \Gamma, \varphi, \psi,\neg \varphi \quad \text{and} \quad
\Lambda \mprf \Gamma, \varphi, \psi,\neg \psi,
\]
so by $\vee$3-rule we have $\Lambda \mprf \Gamma,\varphi,\psi,\neg (\varphi
\vee \psi)$ and then by \ref{weakening-rule} and \ref{cut-rule} we get
$\Lambda \mprf \Gamma,\varphi,\psi$.

3. One direction is \ref{neg-rule}. Suppose $\Lambda \mprf \Gamma, \neg \neg
\varphi$, then by \ref{weakening-rule} we have
\[
\Lambda \mprf \Gamma, \varphi, \neg \neg \varphi
\]
and by using \ref{Axiom:1} and \ref{weakening-rule} we get
\[
\Lambda \mprf \Gamma,\varphi, \neg \varphi,
\]
so, finally, \ref{cut-rule} gives us $\Lambda \mprf \Gamma, \varphi$.
\end{proof}

\begin{defin}
  If $\Gamma$ is a set of \LM-sentences, then it is said to be {\em consistent
    in \M-logic} if $\Gamma \nmprf \emptyset$.
\end{defin}

\begin{prop}\label{prop:consistency.equiv.mlogic}
  If $\Gamma$ is a set of \LM-sentences then the following statements are
  equivalent:
\begin{enumerate}
\item $\Gamma$ is consistent in \M-logic.
  
\item There exists a \LM-sentence $\varphi$ such that $\Gamma \nmprf \varphi$.
  
\item For all \LM-sentences $\varphi$ either $\Gamma \nmprf \varphi$ or
  $\Gamma \nmprf \neg \varphi$.
\end{enumerate}
\end{prop}
\begin{proof}
  $1 \Rightarrow 2$. Suppose $\Gamma$ proves all \LM-sentences. Let $\varphi$
  be any such, then $\Gamma \mprf \varphi$ and $\Gamma \mprf \neg \varphi$. By
  \ref{cut-rule} $\Gamma \mprf \emptyset$, i.e., $\Gamma$ is inconsistent.
  
  $2 \Rightarrow 3$. Suppose there is a \LM-sentence $\varphi$ such that
  $\Gamma \mprf \varphi$ and $\Gamma \mprf \neg \varphi$ then by
  \ref{cut-rule} $\Gamma \mprf \emptyset$ and so by \ref{weakening-rule}
  $\Gamma \mprf \psi$ for any \LM-sentence $\psi$.
  
  $3 \Rightarrow 1$. Suppose $\Gamma$ is inconsistent; by \ref{weakening-rule}
  $\Gamma$ proves every \LM-sentence, so there certainly exists a \LM-sentence
  $\varphi$ such that $\Gamma \mprf \varphi$ and $\Gamma \mprf \neg \varphi$.
\end{proof}

It follows from Lemma \ref{lem:deductive.or.neg} that if $\Gamma$ is
consistent (in \M-logic) and $\varphi$ is a \LM-sentence then $\Gamma,\varphi$
or $\Gamma, \neg \varphi$ is consistent (in \M-logic).

\begin{prop}\label{prop:mlogic.proves.thm}
  If $\varphi \in \ElDiag(\M)$ then $\mprf \varphi$.
\end{prop}
\begin{proof}
  The proof is by induction on the construction of $\varphi$.
  
  First we prove that if $t$ is any closed \LMs-term and
\[
\M \models t=\c_a \quad \text{then} \quad\mprf t=\c_a.
\]
This is done by induction on the construction of $t$. The base case, when $t$
is a constant, follows from \ref{Axiom:3}. We prove the case when $t$ is $r+s$
for some closed terms $r$ and $s$. Let $b$, $d \in \M$ be such that
\[
\M \models r=\c_b \wedge s=\c_d \wedge b+d=a.
\]
By the induction hypothesis $\mprf r=\c_b$ and $\mprf s=\c_d$, so by
\ref{Axiom:7} $\mprf t=\c_b+\c_d$ and by \ref{Axiom:10} and \ref{Axiom:5}
$\mprf t=\c_a$. The other cases, when $t$ is $\Sc(r)$ or $r \cdot s$ are
similar.

Now we prove, by induction on the construction of $\varphi$, that
\begin{align*}
  \text{if} \quad \M \models \varphi \quad&\text{then}\quad \mprf \varphi
  \qquad\text{and} \\
  \text{if} \quad\M \nmodels \varphi \quad&\text{then}\quad \mprf \neg
  \varphi.
\end{align*}
For the base case, when $\varphi$ is an atomic sentence, assume first that $\M
\models t=r$. Then, clearly, $\M \models t=\c_a \wedge r=\c_a$ for some $a \in
\M$. By the fact proved above $\mprf t=\c_a$ and $\mprf r=\c_a$, and by using
\ref{Axiom:3} and \ref{Axiom:4} by get $\mprf t=r$.

If $\M \models t \neq r$ then $\M \models t=\c_a \wedge r=\c_b$ for some $a
\neq b \in \M$. The formal deduction of $t \neq r$ from $t
= \c_a$; $r = \c_b$ and $\c_a \neq \c_b$ is written out in Figure
\ref{fig:mproof} as an illustration of a \M-logic proof.

\begin{sidewaysfigure}
\begin{equation*}
\infer{t \neq r}{
  \infer{t \neq r, \c_a = \c_b}{
    \infer{t \neq r, \c_a = \c_b, r \neq \c_b}{
      \infer{t \neq r, \c_a = \c_b, r \neq \c_b, \c_a \neq r}{
        \infer{\c_a = \c_b, r \neq \c_b, \c_a \neq r}{\text{\ref{Axiom:5}}}
      } &
      \infer{t \neq r, \c_a = \c_b, r \neq \c_b, \c_a = r}{
        \infer{t \neq r, \c_a = \c_b, \c_a = r}{
          \infer{t \neq r,\c_a = r}{
            \infer{t \neq r,\c_a = r, \c_a \neq t}{\text{\ref{Axiom:5}}} &
            \infer{t \neq r,\c_a = r, \c_a = t}{
              \infer{t \neq r, \c_a = t}{
                \infer{\c_a=t}{
                  \infer{\c_a=t,t\neq \c_a}{\text{\ref{Axiom:4}}} &
                  t=\c_a
                }
              }
            }
          }
        }
      }
    } &
    \infer{t \neq r, \c_a = \c_b, r=\c_b}{
      \infer{t \neq r,r=\c_b}{r=\c_b}
    }} &
  \infer{t \neq r, \c_a \neq \c_b}{
    \infer{\c_a \neq \c_b}{\text{\ref{Axiom:2}}}
  }
}
\end{equation*}
  \caption{A proof of $t \neq r$ from the hypothesis $t=\c_a$; $r=\c_b$ 
    and $a \neq b$.}
  \label{fig:mproof}
\end{sidewaysfigure}

For the inductive step we only handle the case when $\varphi$ is $\psi \vee
\sigma$, the others are similar and easy. Assume $\M \models \psi \vee \sigma$
then $\M \models \psi$ or $\M \models \sigma$ so by the induction hypothesis
we have $\mprf \psi$ or $\mprf \sigma$. In either case $\mprf \psi \vee
\sigma$ by \ref{v1-rule} or \ref{v2-rule}.

On the other hand if $\M \nmodels \psi \vee \sigma$ then $\M \nmodels \psi$
and $\M \nmodels \sigma$. Thus, by the induction hypothesis $\mprf \neg \psi$
and $\mprf \neg \sigma$. An application of \ref{v3-rule} yields $\mprf \neg
(\psi \vee \sigma)$.
\end{proof}

\subsection{Soundness and completeness}

If $\Gamma=\Set{\gamma_1,\gamma_2,\ldots,\gamma_k}$ is a finite set of
\LM-sentences, let $\Vee \Gamma$ be
\[
\gamma_1 \vee (\gamma_2 \vee (\ldots \vee \gamma_k))
\]
if $\Gamma$ is nonempty and $0 \neq 0$ if $\Gamma$ is empty.\footnote{We will
  later redefine $\Vee \Gamma$ so that it does not depend on the numbering of
  the set by choosing a canonical numbering.} If \Sat\ is a satisfaction class
then $\Vee \Gamma \in \Sat$ iff $\gamma_i \in \Sat$ for some $1\leq i \leq k$.
Please note that this is independent on the numbering of the set $\Gamma$.

\begin{prop}[Soundness of \M-logic]\label{prop:sound.mlogic}
  If $\mprf \Gamma$ and \Sat\ is a satisfaction class then $\Vee \Gamma \in
  \Sat$.
\end{prop}
\begin{proof}
  The proof is by induction on the height of the proof $p$ of $\Gamma$.  If
  $\abs{p}=0$ then $\Gamma$ is an axiom:
  
  \ref{Axiom:1}: If $\Gamma$ is $\psi,\neg\psi$, it is clear that $\psi \vee
  \neg \psi \in \Sat$ due to the definition of satisfaction classes.
  
  \ref{Axiom:2}: If $a\neq b$ then $\c_a \neq \c_b \in \Sat$ by (\ref{s2}).
  
  \ref{Axiom:3}: If $\Gamma$ is $t=t$ then $t=t \in \Sat$ by (\ref{s3}).
  
  \ref{Axiom:4}: This is (\ref{s4}).
  
  \ref{Axiom:5}: This is (\ref{s5}).
  
  \ref{Axiom:6}: If $\Gamma$ is $t\neq r,\Sc(t)=\Sc(r)$ we have to prove that
  $t \neq r \vee \Sc(t)=\Sc(r) \in \Sat$ or equivalent that if $t=r \in \Sat$
  then $\Sc(t)=\Sc(r) \in \Sat$. If $t=r \in \Sat$ let $a$ be such that
  $t=\c_a \in \Sat$ and $r=\c_a \in \Sat$. We have $\Sc(t)=\c_b \in \Sat$ and
  $\Sc(r)=\c_b \in \Sat$, where $b=\Sc^\M(a)$, by (\ref{s6}). Thus, by
  (\ref{s3})-(\ref{s5}), $\Sc(t)=\Sc(r) \in \Sat$.
  
  \ref{Axiom:7}: As above but use (\ref{s7}) instead of (\ref{s6}).
  
  \ref{Axiom:8}: As above but use (\ref{s8}) instead of (\ref{s6}).
  
  \ref{Axiom:9}, \ref{Axiom:10}, \ref{Axiom:11}: Directly from
  Proposition~\ref{prop:thm.in.satcl}.
  
  \ref{Axiom:12}: $\exists \v_0(t=\v_0) \in \Sat$ since by (\ref{s3}) we have
  $t=t \in \Sat$ and by (\ref{s9}) there exists $a \in \M$ such that $t=\c_a
  \in \Sat$ and so $\exists \v_0 (t=\v_0) \in \Sat$ by (\ref{s12}).
  
  For the inductive step we prove that all the inference rules are sound in
  the sense that if all disjunctions of the premises of an inference rule are
  in in a satisfaction class then so is the disjunction of the conclusion.
  This is easy to see for the first six rules. We prove it only for
  \ref{e-rule} and \ref{M-rule}.
  
  \ref{e-rule}: Suppose $\Vee \Lambda \vee \psi[\c_a/\v_i] \in \Sat$; then we
  can assume that
\[
\psi[\c_a/\v_i] \in \Sat
\]
since otherwise $\Vee \Lambda \in \Sat$ and $\Vee \Lambda \vee \exists \v_i
\,\psi \in \Sat$.  Thus, by (\ref{s1}), $\exists \v_i \,\psi \in \Sat$ and so
$\Vee \Lambda \vee \exists \v_i \,\psi \in \Sat$.

\ref{M-rule}: Suppose $\Vee \Lambda \vee \neg \psi[\c_a/\v_i] \in \Sat$ for
all $a \in \M$. We can assume
\[
\neg \psi[\c_a/\v_i] \in \Sat
\]
for all $a \in \M$. If $\neg \exists \v_i \,\psi \notin \Sat$ then $\exists
\v_i \,\psi \in \Sat$ and then $\psi[\c_a/\v_i] \in \Sat$ for some $a \in \M$
by (\ref{s1}), and we have a contradiction.
\end{proof}

\begin{rem}
  There is a stronger version of the proposition: If $\Lambda \mprf \Gamma$
  and $\Sat$ is a satisfaction class extending the set $\Lambda$, i.e.,
  $\Lambda \subseteq \Sat$, then $\Vee \Gamma \in \Sat$.
\end{rem}

Suppose \M\ is countable and the set $\Lambda$ is a set of \LM-sentences which
is consistent (e.g., if \M-logic is consistent itself we can choose $\Lambda$
to be the empty set). We will construct a satisfaction class including
$\Lambda$. The construction is very much as the completeness theorem for
first-order logic.

Let $\Set{\varphi_i}_{i=1}^\infty$ be an enumeration of all \LM-sentences and
let $\Gamma_0=\Lambda$. Define the sequence $\Set{\Gamma_i}_{i=0}^\infty$ of
consistent sets of \LM-sentences recursively as follows.

We know that if $\Gamma_i$ is consistent then either $\Gamma_i,\varphi_{i+1}$
or $\Gamma_i,\neg\varphi_{i+1}$ is consistent. Choose $\Gamma_{i+1}$ to be the
one which is consistent with the extra condition that if $\varphi_{i+1}$ is
$\exists \v_j \,\psi$ for some \LM-formula $\psi$ and we choose $\Gamma_{i+1}$
to be $\Gamma_i,\exists \v_j \,\psi$ then we also put $\psi[\c_a/\v_j]$ in
$\Gamma_{i+1}$ in such a way that this new set is consistent. We can always do
this, because otherwise we would have
\[
\Gamma_i,\exists \v_j \,\psi \mprf \neg \psi[\c_a/\v_j]
\]
for all $a \in \M$, but then
\[
\Gamma_i,\exists \v_j \,\psi \mprf \neg \exists \v_j \,\psi
\]
by \ref{M-rule} so by \ref{cut-rule} $\Gamma_i,\exists \v_i\,\psi \mprf
\emptyset$ contradicting the fact that $\Gamma_i,\exists\v_i\,\psi$ is
consistent.

Let \[\Gamma_\infty = \bigcup \Gamma_i.\]

\begin{lem}\label{lem:maximal.consistent.mlogic}
  $\Gamma_\infty$ is a maximally consistent set of \LM-sentences.
\end{lem}
\begin{proof}
  It is clearly maximal from the construction: for every sentence $\varphi$ we
  have $\varphi \in \Gamma_\infty$ or $\neg \varphi \in \Gamma_\infty$. To
  prove the consistency of $\Gamma_\infty$ we first notice that by the
  construction there is no sentence $\varphi$ such that $\varphi \in
  \Gamma_\infty$ and $\neg \varphi \in \Gamma_\infty$. Then we observe that
  for all rules with finitely many premises if the premises are in
  $\Gamma_\infty$ then the conclusion is too. This is easy to see since all
  the premises are in some $\Gamma_n$. Assume now that $\neg
  \varphi[\c_a/\v_i] \in \Gamma_\infty$ for all $a \in \M$ and that $\exists
  \v_0 \,\varphi \in \Gamma_\infty$, by the construction we have
  $\varphi[\c_a/\v_i] \in \Gamma_\infty$ for some $a$ which is a
  contradiction. We have proved that $\Gamma_\infty$ is closed under all the
  inference rules so since there is no \LM-sentence $\varphi \in
  \Gamma_\infty$ such that $\neg \varphi \in \Gamma_\infty$, $\Gamma_\infty$
  is consistent.
\end{proof}

\begin{prop}
  \Sat\ is a maximally consistent set of \LM-sentences iff \Sat\ is a
  satisfaction class.
\end{prop}
\begin{proof}
  If \Sat\ is a satisfaction class we can prove the consistency and the
  maximality very much as in Lemma~\ref{lem:maximal.consistent.mlogic}.  We
  prove the converse.
  
  Assume \Sat\ is a maximally consistent set of \LM-sentences. We will use
  Proposition~\ref{prop:alt.satcl.def} to prove that \Sat\ is a satisfaction
  class. First we have to prove that $\sim$ is an equivalence relation;
  \ref{Axiom:3}, \ref{Axiom:4} and \ref{Axiom:5} takes care of that.
  \ref{Axiom:6}, \ref{Axiom:7} and \ref{Axiom:8} takes care of the
  well-definability of $\M_\Sat$.  The canonical map $f: \M \to \M_\Sat$ is an
  embedding by \ref{Axiom:2}, \ref{Axiom:9}, \ref{Axiom:10} and
  \ref{Axiom:11}, $f$ is surjective by \ref{Axiom:12}.
\end{proof}

Adding up the results in this section we get the following:

\begin{prop}\label{prop:sound.complete.mlogic}
  If $\Lambda$ is a set of \LM-sentences then there is a satisfaction class
  $\Sat$ such that $\Lambda \subseteq \Sat$ iff $\Lambda \nmprf \emptyset$.
\end{prop}

\begin{rem}
  We can define $\Lambda \mmodels \Gamma$ for sets of \LM-sentences $\Gamma$
  and $\Lambda$, where $\Gamma$ is finite, to hold iff $\Vee \Gamma$ is
  included in every satisfaction class extending $\Lambda$. Then
  Proposition~\ref{prop:sound.complete.mlogic} can be reformulated as
\[
\Lambda \mprf \Gamma \quad\text{iff} \quad\Lambda \mmodels \Gamma.
\]
\end{rem}

Finally, we add a small remark for the confused reader.

\begin{rem}
  Even though by definition a consistent set $\Gamma$ is maximally consistent
  in a logic if for every sentence or formula $\varphi \notin \Gamma$, the set
  $\Gamma \cup \Set{\varphi}$ is inconsistent; in a logic where the deduction
  theorem holds a consistent set $\Gamma$ is maximally consistent iff for all
  sentences or formulas $\varphi$ either $\varphi \in \Gamma$ or $\neg \varphi
  \in \Gamma$.
\end{rem}

\subsection{The height of proofs}

The next result is about recursively saturated models. This is the only time
we use the recursive saturation in this chapter. It tells us that in
recursively saturated models we can prove everything provable in \M-logic with
proofs of finite height.

\begin{prop} \label{prop:finite.proofs.mlogic}
  If $\Lambda$ is an \LMs-definable set of \LM-sentences, $\Lambda \mprf
  \Gamma$ and \M\ is recursively saturated then $\Lambda \mprf[\omega]
  \Gamma$.
\end{prop}
\begin{proof}
  Let $\lambda(x)$ define $\Lambda$. We will recursively define formulas
  $\Pf_k(x)$, $k \geq 1$ such that
\[
\M \models \Pf_k(a) \quad \text{iff} \quad\Lambda \mprf[k] \Gamma,
\]
where $a$ is a code for the finite set $\Gamma$.

Let $\Pf_1(x)$ be
\begin{gather*}
  \forall z \bigl(z \in x \imp \sent(z)\bigr) \\
  {}\wedge \bigl[ \exists y \bigl( \sent(y) \wedge x = \Set{y} \wedge \lambda(y)\bigr)\\
  {}\vee \exists y \bigl( \sent(y) \wedge x = \Set{y, \godel{\neg y}}\bigr) \\
  {}\vee \exists y \bigl(x= \Set{\godel{\c_y=\c_y}}\bigr)\\
  {}\vee \exists y,z \bigl(x = \Set{\godel{\c_y \neq \c_z}} \wedge y \neq z\bigr)\\
  {}\vee \exists y \bigl(\clterm(y) \wedge x= \Set{\godel{y=y}}\bigr) \\
  {}\vee \exists y,z \bigl(\clterm(y) \wedge \clterm(z) \wedge
  x=\Set{\godel{y\neq z}, \godel{z=y}}\bigl)\\
\begin{split}
  {}\vee \exists y,z,w \bigl( \clterm(y) \wedge \clterm (w) \wedge &
  \clterm (z) \wedge  \\
  &x= \Set{\godel{y\neq z}, \godel{z\neq w}, \godel{y=w}}\bigr)
\end{split}\\
{}\vee \exists y,z \bigl(\clterm(y)\wedge x=
\Set{\godel{y\neq\c_z}, \godel{\Sc(y)= \Sc(\c_z)}}\bigr)\\
\begin{split}
  {}\vee\exists y,z,w,v\bigl(\clterm (y) \wedge & \clterm(z) \wedge\\
  & x= \Set{\godel{y \neq \c_w} ,\godel{z \neq
      \c_v},\godel{y+z=\c_w+\c_v}}\bigr)
\end{split}\\
\begin{split}
  {}\vee\exists y,z,w,v \bigl(\clterm(y) \wedge & \clterm(z) \wedge \\
  & x= \Set{\godel{y \neq \c_w} ,\godel{z \neq \c_v},\godel{y\cdot z=\c_w\cdot
      \c_v}}\bigr)
\end{split}\\
{}\vee \exists y \bigl(\clterm(y) \wedge x= \Set{\godel{\exists
    \v_0(y=\v_0)}}\bigr)\bigr]
\end{gather*}

And let $\Pf_{k+1}(x)$ ($k \geq 1$) be
\begin{gather*}
  \exists y,z \bigl(\sent(z) \wedge \Pf_k(y) \wedge x= y \cup \Set{z}\bigr)\\
  {} \vee \exists y,z,w \bigl(\sent(z) \wedge \sent(w) \wedge \Pf_k(y \cup
  \Set{z}) \wedge x= y \cup \Set{\godel{z \vee w}}\bigr) \\
  {}\vee \exists y,z,w \bigl(\sent(z) \wedge\sent(w) \wedge \Pf_k(y \cup
  \Set{w}) \wedge x= y \cup \Set{\godel{z \vee w}}\bigr) \\
\begin{split}
  {}\vee\exists y,z,w \bigl(\sent(z) \wedge\sent(w) \wedge\Pf_k(y \cup
  \Set{\godel{\neg z}}) \wedge &\Pf_k(y \cup \Set{\godel{\neg w}}) \wedge \\
  &x= y \cup \Set{\godel{\neg (z \vee w)}}\bigr)
\end{split}\\
\vee{} \exists y,z \bigl(\sent(z) \wedge \Pf_k(y \cup \Set{z})
\wedge x= y \cup \Set{\godel{\neg\neg z}}\bigr)\\
{} \vee \exists z \bigl( \sent(z) \wedge \Pf_k (x \cup \Set{z})
\wedge \Pf_k (x \cup \Set{\godel{\neg z}})\bigr) \\
{}\vee \exists y,z,w,i \bigl(\sent(\godel{\exists \v_i \,z}) \wedge \Pf_k (y
\cup \Set{\godel{z[\c_w/\v_i]}}) \wedge x= y \cup \Set{\godel{\exists \v_i \,z}}\bigr) \\
{}\vee \exists y,z,i (\sent(\godel{\neg \exists \v_i \,z}) \wedge \forall w
\Pf_k(y \cup \Set{\godel{\neg z[\c_w/\v_i]}}) \wedge x= y \cup
\Set{\godel{\neg \exists \v_i\, z}}).
\end{gather*}

Observe that if $a \in \M$ is code for a finite set of \LM-sentences $\Gamma$
and $\M \models \Pf_k(a)$ then $\abs{\Gamma} \leq k+2$.

Suppose the lemma is false and let
\[
A=\Set{ \Gamma | \text{$\Lambda \mprf \Gamma $ and $ \Lambda \nmprf[\omega]
    \Gamma$}}
\]
and
\[
B=\Set{ p | \text{$p$ is a proof from $\Lambda$ of $\Gamma$ for some $\Gamma
    \in A$}}.
\]
Let $p \in B$ be of smallest height, i.e., such that if $q \in B$ then
$\abs{p} \leq \abs{q}$. All the subproofs of $p$ must be of finite height, so
it is clear that $\abs{p}=\omega$ and that the last inference rule in $p$ is
\M-rule:
\[
\inference{\Delta,\neg\exists \v_i \,\varphi}{\inferforall{a \in
    \M}{\Delta,\neg \varphi[\c_a/\v_i]}}
\]
Let $a \in \M$ be a code for $\Delta$ and define
\[
q(x)=\Set{\neg \Pf_k\bigl(a \cup \Set{\godel{\neg\varphi[\c_x/\v_i]}}\bigr)
   |  k \geq 1}.
\]
It should be clear that $q(x)$ is a recursive type so, by recursive
saturation, it is realized by some $b \in \M$. We therefore have
\[
\nmprf[\omega] \Delta,\neg\varphi[\c_b/\v_i],
\]
which contradicts the fact that the height of $p$ is $\omega$.
\end{proof}

\begin{rem}
  In Proposition~\ref{prop:finite.proofs.mlogic} we do not need $\Lambda$ to
  be definable, we only need the expanded structure $\pair{\M, \Lambda}$ to be
  recursively saturated for the proof to work.
\end{rem}

\section{The inconsistency of \M-logic}\label{sec:inc}

In this section we prove that if \M\ admits a satisfaction class then \M\ is
recursively saturated (without any restriction on the cardinality of \M). We
will also mention a strengthening of this result by Smith.

Please do compare the next theorem (and proof) with
Theorem~\ref{thm:par.imp.rec}.

\begin{thm}[\cite{Lachlan:81}]\label{thm:lachlan} Let \M\ be an arbitrary
  nonstandard model of \PA\ admitting a satisfaction class then \M\ is
  recursively saturated.
\end{thm}
\begin{proof}
  Let \Sat\ be a satisfaction class on \M\ and assume \M\ is not recursively
  saturated. Let
\[\Set{\varphi_i(x)}_{i \in \omega}\] be a non-realized recursive
type. We can assume that
\begin{align*}
  \M \models &\forall x \bigl(\varphi_{i+1}(x) \imp \varphi_i(x)\bigr), \\
  \M \models &\exists x \bigl(\neg \varphi_{i+1}(x) \wedge
  \varphi_i(x)\bigr) \quad \text{and} \\
  \M \models &\forall x \,\varphi_0(x),
\end{align*}
for all $i \in \omega$. If not, we can replace the formulas $\varphi_i$ by
$\varphi_i'$ where
\begin{align*}
  \varphi_0'(x)\quad&\text{is} \quad x=x \qquad\text{and}\\
  \varphi_{i+1}'(x) \quad& \text{is} \quad\varphi_i'(x) \wedge \varphi_i(x)
  \wedge \exists y \mathord< x \,\varphi_i'(y).
\end{align*}
Let
\[
\alpha_i(x) \text{ be } \neg \varphi_{i+1}(x) \wedge \varphi_{i}(x)
\quad\text{for all $i \in \omega$}
\]
and
\[
A_i=\Set{a \in \M | \M \models \alpha_i(a)}.
\]
It should be clear that $\Set{A_i}_{i \in\omega}$ is a partition of \M.

We will recursively define a sequence $\Set{\beta_i(x)}_{i \leq \nu}$ of
formulas. Let $\beta_0(x)$ be $\alpha_0(x)$ and if $\beta_i(x)$ is defined let
$\theta_{ij}$ be
\[
\exists x \bigl(\beta_i(x) \wedge \alpha_j(x)\bigr)
\]
and $\beta_{i+1}(x)$ be
\begin{align*}
  \bigl(\neg \exists y \beta_i(y) \wedge \alpha_0(x)\bigr) \vee \Bigl(\exists
  y \,\beta_i(y)& \wedge \bigl(\theta_{i0} \wedge \alpha_1(x)\bigr) \vee
  \bigl[\neg \theta_{i0}\\
  &{}\wedge \bigl(\theta_{i1} \wedge \alpha_2(x)\bigr) \vee
  \bigl[\neg \theta_{i1}\\
  &{}\wedge \bigl(\theta_{i2} \wedge \alpha_3(x)\bigr) \vee \bigl[\neg \theta_{i2}\\
  & \qquad \qquad\vdots \\
  &{}\wedge \bigl(\theta_{ii} \wedge
  \alpha_{i+1}(x)\bigl)\cdots\bigr]\bigr]\bigl]\Bigl).
\end{align*}

The sequence $\beta_i(x)$ is recursive and therefore coded in $\M$ and so
extendable to a nonstandard $\nu \in \M \setminus \omega$, so that the
recursive definition holds for all $i<\nu$. Let
\[
B_i = \Set{a \in \M | \beta_i[\c_a/\v_j] \in \Sat}
\]
for all $i \leq \nu$, where $\v_j$ is the free variable in $\beta_i(x)$.

The idea of the definitions of the $\beta_i$s is that
\[
B_{i+1} =
\begin{case}
  A_0     &\quad\text{if  $B_i=\emptyset$}\\
  A_{k+1} &\quad\text{otherwise, where $k = (\least n \mathord\in \omega)\,
    B_i \cap A_n \neq \emptyset$.}
\end{case}
\]
Since a satisfaction class is able to ``look'' finitely deep into a formula we
can prove the following properties of the sequence $\Set{B_i}_{i \leq \nu}$:
\begin{equation}\label{eq:1}
B_i=A_k \quad \Rightarrow \quad B_{i+1}=A_{k+1} \qquad \text{for
all $i \leq \nu$ and $k \in \omega$}
\end{equation}
by the recursive definition of $\beta_i(x)$, and
\begin{equation}\label{eq:2}
\forall i \mathord\leq \nu \,\exists k \mathord\in \omega \, B_i =
A_k
\end{equation}
since $B_0=A_0$ and if $i > 0$ then either $B_{i-1}$ is the empty set, in
which case $B_i$ is $A_0$, or $B_{i-1}$ is not empty and then there is a least
$k \in \omega$ such that $B_{i-1}$ intersects $A_k$ (since $\Set{A_i}_{i \in
  \omega}$ is a partition of \M) and then $B_i$ is $A_{k+1}$.

Let us now finally define the (external) function $f : \initialseg{\nu} \to
\omega$ such that $f(i)=k$ iff $B_i=A_k$. By the property (\ref{eq:2}) this is
a total function and by property (\ref{eq:1}) $f(i+1)=f(i)+1$, so the sequence
\[
f(\nu) > f(\nu-1) > f( \nu -2) > \ldots
\]
is a strictly decreasing infinite sequence of natural numbers, which
contradicts the well-ordering of the natural numbers. Hence there could not be
a non-realized recursive type and therefore \M\ is recursively saturated.
\end{proof}

There is a somewhat stronger version of the theorem:

\begin{thm}[\cite{Smith:84}]
  If \M-logic is consistent then \M\ is recursively saturated.
\end{thm}
\begin{proof}
  The proof is by modifying the proof of Lachlan's result; defining the sets
  $B_i$ by provability in \M-logic instead of by satisfaction classes, see
  \cite{Smith:84} for the details.
\end{proof}

Please observe that when \M\ is countable \M\ admits a satisfaction class
precisely when \M-logic is consistent, so the result of Smith is a
strengthening of Lachlan's result only when \M\ is uncountable.

\section{The consistency of \M-logic}

\subsection{Template logic}

In Section~\ref{sec:constructionofsatcl} we proved that \M-logic is consistent
iff we can find a satisfaction class for the model (assuming \M\ is
countable). In this section we prove that for any recursively saturated model
of \PA\ \M-logic is consistent.

Together with the results in Section~\ref{sec:inc} above we get that a
countable nonstandard model of \PA\ admits a satisfaction class iff it is
recursively saturated.

The question is; if \M\ is recursively saturated how do we prove the
consistency of \M-logic? Proving consistency of a logic can be done, mainly,
in two different ways; the proof theoretic way, by a cut-elimination theorem,
or the model theoretic way, by a soundness theorem. We will use the model
theoretic approach and prove a soundness theorem. For this we will define a
new kind of logic. The idea is really easy; instead of studying nonstandard
formulas we are going to replace some subformulas and subterms in the
nonstandard formula with templates, and in this way study formulas of finite
depth.

We will call this logic template logic, it is first-order logic with template
symbols added. Each \LM-formula and \LM-term has a corresponding template
symbol. These symbols may be looked on as predicates and functions of
nonstandard finite arity. For example, the template symbol corresponding to the
\LM-formula
\[
\v_0 = \v_0 \vee \v_1 = \v_1 \vee \ldots \vee \v_a = \v_a
\]
could be treated as a predicate of arity $a$, even if $a > \omega$.

Let \LT\ be the language
\[\LMs \cup \Set{\tmpl{\varphi}| \text{$\varphi$ is a \LM-formula}}
\cup \Set{\tmpl{t}| \text{$t$ is a \LM-terms}}
\]
where $\tmpl{\varphi}$ is treated as 0-ary relational symbols (i.e.,
propositional variables) and $\tmpl{t}$ is treated as constant symbols when
building up terms and formulas. For technical reasons we include the variables
$\Set{\v_a | a \in \M}$ in the language \LT.

The formulas and terms of template logic are defined in the usual way with
template symbols being 0-ary. The free variables of a term is defined in a
nonstandard way as follows
\begin{align*}
  \FV(\c_a)&\eqdef\emptyset \\
  \FV(\v_i)&\eqdef\Set{\godel{\v_i}} \\
  \FV(\,\tmpl{t}\,)&\eqdef\Set{a \in \M | \M \models a \in \FV(t)} \\
  \FV(\Sc(r))&\eqdef\FV(r)\\
  \FV(r + s)&\eqdef\FV(r) \cup \FV(s) \\
  \FV(r \cdot s)&\eqdef\FV(r) \cup \FV(s)
\end{align*}
Observe that in the third clause $\FV$ has two different meanings, in the
second appearance it is the function definable in \PA. And of a formula as
\begin{align*}
  \FV(r=s)&\eqdef\FV(r) \cup \FV(s) \\
  \FV(\tmpl{\varphi})&\eqdef\Set{a \in \M | \M \models a \in \FV(\varphi)}\\
  \FV(\neg \gamma) &\eqdef \FV(\gamma) \\
  \FV(\gamma \vee \delta) &\eqdef \FV(\gamma) \cup \FV(\delta) \\
  \FV(\exists \v_i \,\gamma)&\eqdef\FV(\gamma) \setminus \Set{\godel{\v_i}}
\end{align*}
We call a formula $\gamma$ a {\em sentence} if $\FV(\gamma)=\emptyset$ and a
term $t$ {\em closed\/} if $\FV(t)=\emptyset$.

Substitution is also defined in a nonstandard way; for terms
\begin{align*}
  \c_b[\c_a/\v_i] &\quad\text{is}\quad \c_b\\
  \v_j[\c_a/\v_i] &\quad\text{is}\quad
\begin{cases}
  \c_a & \quad\text{if $i=j$} \\
  \v_j & \quad\text{otherwise}
\end{cases}\\
\tmpl{t}\,[\c_a/\v_i] &\quad\text{is}\quad \tmpl{t[\c_a/\v_i]}\\
\Sc(r)[\c_a/\v_i] &\quad\text{is}\quad \Sc(r[\c_a/\v_i])\\
(r + s) [\c_a/\v_i] &\quad\text{is}\quad r[\c_a/\v_i] +
s[\c_a/\v_i] \\
(r \cdot s) [\c_a/\v_i] &\quad\text{is}\quad r[\c_a/\v_i] \cdot s[\c_a/\v_i]
\end{align*}
and for formulas
\begin{align*}
  (r=s)[\c_a/\v_i] &\quad\text{is}\quad r[\c_a/\v_i] =s[\c_a/\v_i]  \\
  \tmpl{\varphi}[\c_a/\v_i] &\quad\text{is}\quad\tmpl{\varphi[\c_a/\v_i]}\\
  (\neg \gamma)[\c_a/\v_i] &\quad\text{is}\quad \neg(\gamma[\c_a/\v_i])\\
  (\gamma \vee \delta)[\c_a/\v_i] &\quad\text{is}\quad \gamma[\c_a/\v_i]
  \vee \delta[\c_a/\v_i]\\
  \exists \v_i \,\gamma [\c_a/\v_j] &\quad\text{is}\quad
\begin{cases} \exists \v_i (\gamma[\c_a/\v_j]) &\quad\text{if $i
    \neq j$}\\ \exists \v_i \,\gamma &\quad\text{otherwise.}
\end{cases}
\end{align*}

\begin{defin}\label{def:substitution}
  If $\varphi$ is a \LM-formula and $a \in \M$ then $\varphi[a/\v]$ is
\[\varphi[\c_{[a]_{i_1}-1}, \ldots, \c_{[a]_{i_k}-1}/ \v_{i_1},
\ldots, \v_{i_k}]\] where $\Set{i_1, \ldots, i_k} = \Set{i \in \M  |  [a]_i
  \neq 0}$. We define $t[a/\v]$, for terms $t$, in the obvious similar way.
\end{defin}

\begin{defin} The relation (between \LM-formulas) $\varphi \cong \psi$ holds
  if there exists a \LM-formula $\gamma$ and $a$, $b \in \M$ such that
  $\gamma[a/\v]$ is $\varphi$ and $\gamma[b/\v]$ is $\psi$.  The relation
  between \LM-terms $t \cong r$ holds if there exists a \LM-term $s$ and $a$, $b \in
  \M$ such that $s[a/\v]$ is $t$ and $s[b/\v]$ is $r$.
\end{defin}

\begin{prop}
  The relation $\cong$ is an equivalence relation on the set of \LM-formulas
  and \LM-terms.
\end{prop}
\begin{proof}
  The reflexive and symmetric properties are trivial, only the transitive
  property involves some work. We prove the proposition for terms, the case
  with formulas is similar.
  
  Let us assume that $t \cong r$ and $r \cong s$, and prove $t \cong s$. The
  proof is by induction on $r$ inside \PA; we use the inductive property of
  $\term(x)$. If $r$ is a constant or a variable then $t$ and $s$ are also
  either constants or variables.  In any case $t \cong s$ since any constant
  is related to any constant or variable and a variable is related to a
  variable iff they are equal.
  
  For the inductive step the case when $r$ is of the form $\Sc(r')$ is easy
  and the two cases when $r$ is $r_1+r_2$ or $r_1 \cdot r_2$ are similar,
  therefore we only handle the case when $r$ is $r_1 + r_2$.
  
  It is easy to see that $t$ and $s$ also are of this form, i.e., $t$ is $t_1
  + t_2$ and $s$ is $s_1 + s_2$. Let $p$, $q$, $a_t$, $a_r$, $b_r$, and $b_s$
  be such that
\[
t=p[a_t/\v], r= p[a_r/\v], r=q[b_r/\v] \quad\text{and}\quad s=q[b_s/\v]
\]
and let $p$ be $p_1 + p_2$ and $q$ be $q_1 + q_2$. Clearly,
\[
t_i=p_i[a_t/\v], r_i=p_i[a_r/\v], r_i=q_i[b_r/\v] \quad\text{and}\quad
s_i=q_i[b_s/\v]
\]
for $i=1$, $2$, therefore $t_i \cong r_i$ and $r_i \cong s_i$.  Thus, by the
induction hypothesis $t_i \cong s_i$. Let $o_i$, $a_{it}$ and $a_{is}$ be such
that
\[
t_i=o_i[a_{it}/\v] \quad\text{and}\quad s_i=o_i[a_{is}/\v].
\]
We have to ``unify'' $a_{it}$ and $a_{is}$ in such a way that the results,
$a_t$ and $a_s$, work for both $i$s, i.e., such that $t_i=o_i'[a_t/\v]$ and
$s_i=o_i'[a_s/\v]$, where $o_i$ are some new terms. Then clearly $t=(o_1' +
o_2')[a_t/\v]$ and $s=(o_1 + o_2)[a_s/\v]$ so $t \cong s$.

To make this happen let $o_1'$ be $o_1$ and let $i_0 \in \M$ be bigger than
all indices of free variables of $o_1$ and all $i$s such that
\[
[a_{1t}]_i\neq 0 \quad\text{or}\quad [a_{1s}]_i \neq 0.
\]
Let $o_2'$ be as $o_2$ except that if $\v_i$ occurs as a free variable in
$o_2$ and $[a_{2t}]_i \neq 0$ it is substituted by $\v_{i+i_0}$. Finally,
define $a_t$ and $a_s$ such that
\begin{align*}
  [a_t]_i &= [a_{1t}]_i \quad \text{and}\\
  [a_s]_i &= [a_{1s}]_i
\end{align*}
for all $i < i_0$ and
\begin{align*}
  [a_t]_{i+i_0} &= [a_{2t}]_i \quad \text{and}\\
  [a_s]_{i+i_0} &= [a_{2s}]_i
\end{align*}
for all $i \in \M$. Clearly, $t_i=o_i'[a_t/\v]$ and $s_i=o_i'[a_s/\v]$ for
$i=1$, $2$.
\end{proof}

So far, so good; we have a new logic to play with (even though we have not
defined the axioms and inference rules yet). But how does template logic
connect with \M-logic, the object of study? For the connection to work we need
a way to approximate a \LM-formula or a \LM-term with a template formula or
term.  Given a \LM-formula you can, by replacing some of the subformulas and
subterms with corresponding template symbols, e.g., replacing the subformula
$\varphi$ by $\tmpl{\varphi}$, make an \LT-formula. This is the idea behind
approximations.

If $\psi$ is an \LT-formula and $\delta$ a \LM-formula we define
$\F_{\delta}(\psi)$ to be $\psi$ with all occurrences of symbols
$\tmpl{\varphi}$, for $\varphi \cong \delta$, replaced by
\begin{align*}
  \tmpl{t_1}=\tmpl{t_2}& \quad\text{if $\varphi$ is $t_1=t_2$}\\
  \neg \tmpl{\gamma} & \quad\text{if $\varphi$ is $\neg \gamma$}\\
  \tmpl{\gamma}\vee \tmpl{\sigma} & \quad\text{if
    $\varphi$ is $\gamma \vee \sigma$}\\
  \exists \v_i \,\tmpl{\gamma}& \quad\text{if $\varphi$ is $\exists
    \v_i\,\gamma$}
\end{align*}
and we define $\F_{r}(\psi)$, where $r$ is a \LM-term, to be $\psi$ with all
occurrences of symbols $\tmpl{t}$, for $t \cong r$, replaced by
\begin{align*}
  \v_i & \quad\text{if $t$ is $\v_i$}\\
  \c_a & \quad\text{if $t$ is $\c_a$}\\
  \Sc(\tmpl{r}) & \quad\text{if $t$ is $\Sc(r)$}\\
  \tmpl{r} + \tmpl{s} & \quad\text{if $t$ is $r + s$} \\
  \tmpl{r}\cdot \tmpl{s} & \quad\text{if $t$ is $r\cdot s$.}
\end{align*}

Now, the definition of an approximation.

\begin{defin}
  An \LT-formula (or \LT-term) $\psi$ is an approximation of another
  \LT-formula (or \LT-term) $\delta$ if there exists \LM-formulas or -terms
  $\tau_1$, \dots, $\tau_k$ such that
\[
\delta=\F_{\tau_k} \circ \F_{\tau_{k-1}} \circ \ldots \circ \F_{\tau_1}(\psi).
\]
\end{defin}

If
\[
\F = \F_{\tau_k} \circ \F_{\tau_{k-1}} \circ \ldots \circ \F_{\tau_1}
\]
we call $\F$ an approximating function and say that the length of $\F$,
denoted $\abs{\F}$, is $k$. We say that $\psi$ is an approximation of a
\LM-formula $\varphi$ if $\psi$ is $\F(\tmpl{\varphi})$, for some
approximating function $\F$.

We define approximations of finite sets of formulas by letting
\[
\F(\Delta)=\Set{\F(\delta) | \delta \in \Delta}.
\]

Finally, for convenience, we define $\F(\varphi)$, where $\varphi$ is a
\LM-sentence, to be $\F(\tmpl{\varphi})$.

The formal proof system for template logic as just like the one for \M-logic,
but for completeness we write it down anyway. The axioms are
\begin{gather}
  \tag{Axiom1$_t$}\label{Axiom:1'}\gamma, \neg \gamma\\
  \tag{Axiom2$_t$}\label{Axiom:2'}\c_a \neq \c_b \quad\text{if  $a \neq b$}\\
  \tag{Axiom3$_t$}\label{Axiom:3'}t=t \\
  \tag{Axiom4$_t$}\label{Axiom:4'}t\neq r,r=t\\
  \tag{Axiom5$_t$}\label{Axiom:5'}t\neq r, r\neq s, t=s \\
  \tag{Axiom6$_t$}\label{Axiom:6'}t\neq r,\Sc(t)=\Sc(r) \\
  \tag{Axiom7$_t$}\label{Axiom:7'}t\neq t',r\neq r',t+r=t'+r' \\
  \tag{Axiom8$_t$}\label{Axiom:8'}t\neq t',r\neq r',t\cdot r=
  t' \cdot r' \\
  \tag{Axiom9$_t$}\label{Axiom:9'} \Sc(\c_a) = \c_{\Sc(a)} \\
  \tag{Axiom10$_t$}\label{Axiom:10'} \c_a + \c_b = \c_{a + b} \\
  \tag{Axiom11$_t$}\label{Axiom:11'} \c_a \cdot \c_b = \c_{a \cdot b} \\
  \tag{Axiom12$_t$}\label{Axiom:12'}\exists \v_0(t=\v_0)
\end{gather}
where $\gamma$ is an arbitrary \LT-sentence, and $t$, $r$ and $s$ are
arbitrary closed \LT-terms. The inference rules are the following:
\begin{equation}
  \tag{Weak$_t$}\label{weakening-rule'}
  \inference{\Gamma,\gamma}{\Gamma}
\end{equation}
\begin{equation}
  \tag{$\vee$I1$_t$}\label{v1-rule'}
  \inference{\Gamma,\gamma \vee \delta}{\Gamma,\gamma}
\end{equation}
\begin{equation}
  \tag{$\vee$I2$_t$}\label{v2-rule'}
  \inference{\Gamma,\gamma \vee \delta}{\Gamma,\delta}
\end{equation}
\begin{equation}
  \tag{$\vee$I3$_t$}\label{v3-rule'}
  \inference{\Gamma,\neg (\gamma \vee \delta)}{\Gamma,\neg \gamma &
  \Gamma,\neg \delta}\\
\end{equation}
\begin{equation}
  \tag{$\neg$I$_t$}\label{neg-rule'}
  \inference{\Gamma,\neg\neg\gamma}{\Gamma,\gamma}\\
\end{equation}
\begin{equation}
  \tag{Cut$_t$}\label{cut-rule'}
  \inference{\Gamma}{\Gamma, \gamma& \Gamma, \neg\gamma}\\
\end{equation}
\begin{equation}
  \tag{$\exists$I$_t$}\label{e-rule'}
  \inference{\Gamma, \exists \v_i\,\gamma}{\Gamma,\gamma[\c_a/\v_i]}\\
\end{equation}
\begin{equation}
  \tag{\M-rule$_t$}\label{M-rule'}
  \inference{\Gamma, \neg \exists \v_i\,\gamma}{\inferforall{a \in \M}{\Gamma,
  \neg \gamma[\c_a/\v_i]}}
\end{equation}
where $\Gamma$ is an arbitrary finite set of \LT-sentences, $\gamma$, $\delta$
and $\exists \v_i\, \gamma$ are arbitrary \LT-sentences, $t$ and $r$ are
arbitrary closed \LT-terms and $a \in \M$.

Similar to \M-logic $\tprf[p] \Delta$ means that $p$ is a proof of $\Delta$ in
template logic. All other definitions in \M-logic transform almost verbatim to
template logic, that is also the case for Lemma~\ref{lem:deductive.or.neg} and
Proposition~\ref{prop:consistency.equiv.mlogic}.

Please observe that if we are only studying standard formulas template logic
extends \M-logic, so any proof in \M-logic using only standard formulas and
terms is also a proof in template logic. On the other hand if we restrict
template logic to formulas and terms without template symbols, then \M-logic
is an extension of template logic.

\subsection{Some Technical Results}

We need some more information on how the approximating functions work for
later use.

\begin{lem}\label{lem:subst.and.approx.commute}
  If $\F$ is an approximating function and $\gamma$ is an \LT-formula then
  $\F(\gamma[\c_a/\v_i])=\F(\gamma)[\c_a/\v_i]$.
\end{lem}
\begin{proof}
  Observe first that it is enough to prove the lemma for approximating
  functions $\F=\F_\tau$; the general result follows by ``moving'' one
  $\F_\tau$ at a time.
  
  Since both substitution and approximating functions commute with the symbols
  $\neg$, $\vee$, $\exists$, $=$, $\Sc$, $+$ and $\cdot$ we only have to check
  the base cases, i.e., when $\gamma$ is a constant, variable or template
  symbol.
  
  If $\gamma$ is a constant or variable (or even more generally if $\gamma$
  does not contain any template symbols) we have
\[
\F_\tau(\gamma[\c_a/\v_i]) = \gamma[\c_a/\v_i] = \F_\tau(\gamma)[\c_a/\v_i].
\]

If $\gamma$ is a template symbol, say $\tmpl{\delta}$, it is clear that if
$\tau \ncong \delta$ then
\[
\F_{\tau}(\tmpl{\sigma})[\c_a/\v_i]=\tmpl{\sigma}[\c_a/
\v_i]=\tmpl{\sigma[\c_a/\v_i]}=\F_\tau(\tmpl{\sigma[\c_a/\v_i]}=
\F_\tau(\tmpl{\sigma}[\c_a/\v_i]).
\]
In the third equality we are using the fact that $\tau \ncong
\delta[\c_a/\v_i]$, this is easy to see since if $\tau \cong
\delta[\c_a/\v_i]$ then, since $\cong$ is a equivalence relation and
$\delta[\c_a/\v_i] \cong \delta$, $\tau \cong \delta$.

Suppose $\psi$ is $\neg\gamma$ and $\varphi \cong \psi$ then
\[\F_{\varphi}(\tmpl{\psi})[\c_a/\v_i]=
\neg \tmpl{\gamma}[\c_a/\v_i]= \neg \tmpl{\gamma[\c_a/\v_i]}=
\F_{\varphi}(\tmpl{\psi[\c_a/\v_i]})= \F_\varphi(\tmpl{\psi}[\c_a/\v_i]).\]
The case when $\psi$ is $\gamma \vee \delta$ is treated in a similar way.
Suppose that $\psi$ is $\exists \v_j\,\gamma$ then
\begin{multline*}
  \F_\varphi(\tmpl{\psi})[\c_a/\v_i]=(\exists \v_j \,\tmpl{\gamma})
  [\c_a/\v_i]= \\
\begin{cases}
  \exists \v_j \,\tmpl{\gamma[\c_a/\v_i]}= \F_\varphi (\tmpl{\psi
    [\c_a/\v_i]} )= \F_\varphi (\tmpl{\psi} [\c_a/\v_i]) & \quad\text{if $i\neq j$}\\
  \exists \v_j \,\tmpl{\gamma}= \F_\varphi(\tmpl{\psi}) = \F_\varphi
  (\tmpl{\psi}[\c_a/\v_i]) & \quad\text{otherwise.}
\end{cases}
\end{multline*}
We also have to check the term cases. If the term is a composite term then it
is handled just as the $\neg \gamma$ case. For the other cases we have
\begin{multline*}
  \F_{\v_i}(\tmpl{\v_i})[\c_a/\v_j]=\v_i[\c_a/\v_j]= \\
\begin{cases}
  \v_i=\F_{\v_i}(\tmpl{\v_i}[\c_a/\v_j]) &\quad\text{if $i \neq j$} \\
  \c_a=\F_{\v_i}(\tmpl{\c_a})=\F_{\v_i}(\tmpl{\v_i}[\c_a/\v_j])
  &\quad\text{otherwise.}
\end{cases}
\end{multline*}
and
\[
\F_{\v_j}(\tmpl{\c_b})[\c_b/\v_i]=\c_b=\F_{\v_j}(\tmpl{\c_b} [\c_a/\v_i]).
\qedhere
\]
\end{proof}

\begin{lem}\label{lem:prov.template.implies.prov.aprox}
  If $\Delta$ is a finite set of \LT-sentences, $\F$ is an approximating
  function and $\tprf[\alpha] \Delta$ then $\tprf[\alpha] \F (\Delta)$.
\end{lem}
\begin{proof}
  The proof is by induction on the length of the proof. It should be clear
  that if $\Delta$ is a template axiom then so is $\F(\Delta)$. It is also
  easy, but tedious, to check that the inference rules are not affected by
  $\F$. We will not do it here.
\end{proof}

\begin{rem}
  It should be clear that the lemma could be strengthen as to say that if
  $\Lambda$ is a set of \LT-sentences such that if $\lambda \in \Lambda$ then
  $\F(\lambda) \in \Lambda$ for any approximating function $\F$, $\Delta$ a
  finite set of \LT-sentences, $\F$ an approximating function and $\Lambda
  \tprf[\alpha] \Delta$ then $\Lambda \tprf[\alpha] \F(\Delta)$.
\end{rem}

\begin{lem}\label{lem:uniform.approx.exists}
  If $\F_0$, $\F_1$, \dots\ are approximating functions such that $\abs{\F_i}
  \leq k$ and $\Gamma$ a finite set of \LM-sentences then there is an
  approximating function $\F$ such that
\[
\abs{\F} \leq \bigl(2^k-1\bigr)\abs{\Gamma}
\]
and
\[
\F(\Gamma)=\F(\F_i(\Gamma))
\]
for all $i \in \omega$.
\end{lem}
\begin{proof}
  Let
\[
\F=\F_{\delta_n} \circ \F_{\delta_{n-1}} \circ \ldots \circ \F_{\delta_1}
\]
where $\delta_1$, \dots, $\delta_n$ are all subformulas and subterms occurring
in some formula in $\Gamma$ at depth $\leq k$ in ``the right order,'' i.e., if
$\delta_i$ is a subformula or subterm of $\delta_j$ then $j<i$. Clearly $n\leq
(2^k-1)\abs{\Gamma}$ and $\F(\Gamma)=\F(\F_i(\Gamma))$ for all $i$.
\end{proof}

\begin{defin}\label{def:normal}
  An approximating function \[\F=\F_{\delta_n} \circ \ldots \circ
  \F_{\delta_1}\] is said to be in {\em normal form} if $j < i$ for every
  $i,j$ such that $\delta_i$ is a subformula or subterm of $\delta_j$.
\end{defin}

Observe that if $\F$ is an approximating function and $\F'$ is $\F_{\delta_k}
\circ \ldots \circ \F_{\delta_0}$ where $\F_{\delta_0}$, \dots, $\F_{\delta_k}$
are the approximating functions in $\F$ ordered such that if $\delta_i$ is a
subformula or subterm of $\delta_j$ then $i<j$, then $\F'\circ\F (\Gamma) =
\F'(\Gamma)$ and $\F'$ is in normal form. Therefore if $\F(\Gamma)$ is
provable so is $\F'(\Gamma)$.

{\em From now on we will assume that all approximating functions are in normal
  form.}

If $\F_1, \ldots , \F_k$ are approximating functions we can form {\em the
  uniform union}
\[
\F_1 \uniform \ldots \uniform \F_k
\]
of them which is any normal form of $\F_1 \circ \ldots \circ \F_k$. This
definition is not unique, the reader may try to make it unique in a suitable
way.

Please observe that if $\F_1$ and $\F_2$ are approximating functions and
$\tprf \F(\Gamma)$ then $\tprf \F_1 \uniform \F \uniform \F_2 (\Gamma)$ since
$\F_1 \uniform \F \uniform \F_2 ( \F (\Gamma)) = \F_1 \uniform \F \uniform
\F_2 (\Gamma)$.

\begin{lem}\label{lem:approx.functions.commute}
  Let $\psi$ be any \LM-sentence and $\F$ an approximating function (in normal
  form) such that $\F_\psi$ is in $\F$, then
\begin{itemize}
\item if $\psi$ is $\neg \gamma$ then $\F(\neg \gamma) = \neg \F(\gamma)$,
  
\item if $\psi$ is $\gamma_1 \vee \gamma_2$ then $\F(\gamma_1 \vee \gamma_2)
  =\F(\gamma_1) \vee \F(\gamma_2)$ and
  
\item if $\psi$ is $\exists \v_i \,\gamma$ then $\F(\exists \v_i \,\gamma) =
  \exists \v_i \,\F(\gamma)$.
\end{itemize}
\end{lem}
\begin{proof}
  The proof is more or less trivial and left to the reader.
\end{proof}

\subsection{Semantics}

We will now start to look at the semantics of template logic and prove a
soundness theorem which implies the consistency of the logic. We end the
section by also proving a completeness theorem.

\begin{defin}
  An {\em \LT-structure} $\T$ is a pair $\pair{\T_t,\T_v}$ of a set $\T_t$ of
  \LM-sentences and a map $\T_v$ from the closed \LM-terms into \M.
\end{defin}

\begin{defin}
  If $\T$ is an \LT-structure define $\val_\T(t)$ for closed \LT-terms $t$
  inductively as follows:
\begin{align*}
  \val_\T(\tmpl{r})&\eqdef\T_v(r), \\
  \val_\T(\Sc(r))&\eqdef\Sc(\val_\T(r)), \\
  \val_\T(r + s)&\eqdef\val_\T(r) + \val_\T(s) \quad\text{and} \\
  \val_\T(r \cdot s)&\eqdef\val_\T(r) \cdot \val_\T(s).
\end{align*}
\end{defin}

\begin{defin}
  If \T\ is an \LT-structure then define the predicate $\T \models \varphi$ on
  \LT-sentences $\varphi$ inductively as follows:
\begin{align*}
  \T \models t=r&\quad\text{iff}\quad \val_\T(t)=\val_\T(r), \\
  \T \models \tmpl{\varphi} &\quad\text{iff}\quad \varphi \in\T_t, \\
  \T \models \neg \gamma &\quad\text{iff}\quad \T \nmodels \gamma, \\
  \T \models \gamma\vee \delta &\quad\text{iff}\quad \T \models \gamma
  \text{ or } \T \models \delta \quad\text{and}\\
  \T \models \exists \v_i \gamma &\quad\text{iff}\quad \text{there exists $a
    \in \M$ such that $\T \models \gamma[\c_a/\v_i]$}.
\end{align*}
\end{defin}

We will now prove a soundness property for template logic. Let $\Vee \Delta$
be
\[\delta_0 \vee (\delta_1 \vee (\ldots \vee \delta_k))\]
if $\Delta=\Set{\delta_0,\ldots,\delta_k}$ and $0 \neq 0$ if
$\Delta=\emptyset$.

\begin{prop}[Soundness of template logic]\label{prop:template.sound}
  Let $\Delta$ be a finite set of \LT-sentences, let $\gamma$ be $\Vee \Delta$
  and $\Lambda$ any set of \LT-sentences. If \T\ is an \LT-structure making
  all the sentences of $\Lambda$ true and $\Lambda \tprf \Delta$ then $\T
  \models \gamma$.
\end{prop}
\begin{proof}
  We have to check that all the axioms of template logic are true in all
  \LT-structures and that all the inference rules are sound, i.e., if the
  premises of a rule are true in some \LT-structure then the conclusion is
  also true in the same \LT-structure. The axioms are quite obvious true and
  the inference rules are also easy to check; we only prove that \ref{M-rule'}
  is sound. Suppose $\T \models \neg\psi[\c_a/\v_i]$ for all $a \in \M$ and
  $\T \models \exists \v_i \,\psi$, then there is an $a \in \M$ such that $\T
  \models \psi[\c_a/\v_i]$ which is a contradiction.
\end{proof}

\begin{defin}
  An \LT-sentence $\gamma$ is said to be {\em true} in \M\ if $\T \models
  \gamma$ for all \LT-structures \T.
\end{defin}

\begin{cor}\label{cor:template.sound}
  If $\tprf \Gamma$ then $\Vee \Gamma$ is true in \M.
\end{cor}

Now, we easily get the consistency of template logic.

\begin{cor}
  In any model \M\ we have $\ntprf 0 \neq 0$.
\end{cor}
\begin{proof}
  Use Corollary~\ref{cor:template.sound} and the fact that for any template
  structure $\T$ we have $\T \nmodels 0 \neq 0$.
\end{proof}

\begin{prop}[Completeness of template logic]\label{prop:template.complete}
  If $\psi$ is an \LT-sentence true in \M\ then $\tprf \psi$.
\end{prop}
\begin{proof}
  Suppose $\ntprf \psi$ then $\Set{\neg \psi}$ is a consistent set in template
  logic and could be extended to a maximally consistent set $\Lambda$ in the
  same way as when we constructed satisfaction classes from consistent sets in
  \M-logic.
  
  Define the \LT-structure \T\ by letting
\[
\T_t=\Set{\varphi | \text{$\varphi$ is a \LM-sentence and $\tmpl{\varphi} \in
    \Lambda$}}
\]
and
\[
\T_v(t)=a \quad\text{iff}\quad \tmpl{t}=\c_a \in \Lambda.
\]

We prove that
\[
\varphi \in \Lambda \quad\text{iff}\quad \T \models \varphi
\]
by induction on the construction of $\varphi$. To handle the case when
$\varphi$ is atomic we first prove that if $t$ is a closed \LT-term then
$t=\c_a \in \Lambda$ iff $\val_{\T}(t)=a$.

Suppose $t$ is $\tmpl{r}$, then the claim is trivially true from the
definition of $\T_v$. If $t$ is $\Sc(r)$ and $t=\c_a \in \Lambda$, by
\ref{Axiom:12'}, there is a $b \in \M$ such that $r=\c_b \in \Lambda$ and
$\Sc^\M(b) = a$. By the induction hypothesis $\val_{\T}(r)=b$ so
$\val_{\T}(t)=\Sc^\M(a)=b$. The case when the term is $r + s$ or $r \cdot s$
is handled in a similar way.

If $t=r \in \Lambda$ then by \ref{Axiom:12} (and \ref{Axiom:3}, \ref{Axiom:4}
and \ref{Axiom:5}) there is an $a \in \M$ such that $t=\c_a \in \Lambda$ and
$r=\c_a \in \Lambda$. By the fact proved above $\val_{\T}(t)=\val_{\T}(r)$, so
$\T \models t=r$.

On the other hand if $\T \models t=r$ then $\val_{\T}(t)=\val_{\T}(r)=a$ for
some $a \in \M$. Since, by the maximality of $\Lambda$, there are $b$, $d \in
\M$ such that $t= \c_b \in \Lambda$ and $r = \c_d \in \Lambda$ we have, by the
fact proved above, $b=d=a$. Therefore $t=r \in \Lambda$.

If $\varphi \vee \psi \in \Lambda$ then either $\varphi \in \Lambda$ or $\psi
\in \Lambda$ by the maximality of $\Lambda$ and so, by the induction
hypothesis, either $\T \models \varphi$ or $\T \models \psi$, either way $\T
\models \varphi \vee \psi$. On the other hand, if $\varphi \vee \psi \notin
\Lambda$ then neither $\varphi$ nor $\psi$ is in $\Lambda$ so by the induction
hypothesis $\T \nmodels \varphi$ and $\T \nmodels \psi$ which implies that $\T
\nmodels \varphi \vee \psi$.

If $\neg \varphi \in \Lambda$ then $\varphi \notin \Lambda$ by the consistency
of $\Lambda$ and by the induction hypothesis $\T \nmodels \varphi$, therefore
$\T \models \neg \varphi$. If $\neg \varphi \notin \Lambda$ then $\varphi \in
\Lambda$ so $\T \models \varphi$ and $\T \nmodels \neg \varphi$.

If $\exists \v_i \,\varphi \in \Lambda$ then there exists $a \in \M$ such that
$\varphi[\c_a/\v_i] \in \Lambda$ and so $\T \models \varphi[\c_a/\v_i]$ and
$\T \models \exists \v_i \,\varphi$. And if $\exists \v_i \,\varphi \notin
\Lambda$ then $\varphi[\c_a/\v_i] \notin \Lambda$ for all $a \in \M$. By the
induction hypothesis $\T \nmodels \varphi[\c_a/\v_i]$ for all $a \in \M$ so
$\T \nmodels \exists \v_i \,\varphi$.
\end{proof}

To sum up this section, the main result is that for an \LT-sentence $\gamma$
we have that $\gamma$ is true in \M\ iff $\tprf \gamma$.

\subsection{A link between \M-logic and template logic}

In this section we prove that finite provability in \M-logic implies
provability of some approximation in template logic. Since template logic is
consistent (Theorem~\ref{prop:template.sound}) this will imply that \M-logic
is consistent.

\begin{prop}\label{prop:prov.mlogic.implies.template}
  There is a (recursive) function $G : \omega \to \omega$ such that if
  $\mprf[n] \Gamma$ then there is an approximating function $\F$ such that
  $\abs{\F}\leq G(n)$ and $\tprf \F(\Gamma)$.
\end{prop}
\begin{proof}
  If $n=1$ then $\Gamma$ is an axiom. It is easy to see that an approximating
  function of length 9 is enough to make $\F(\Gamma)$ into an axiom of
  template logic.
  
  Suppose we have defined $G$ for all values $\leq n$ and take a proof of
  height $n+1$, by the induction hypothesis we get a proof of some
  approximation of the premises of the last inference of the proof, with the
  approximating functions of length $\leq G(n)$.  Suppose the last inference
  is
  
  \textbf{\ref{weakening-rule}:}
\[\inference{\Lambda,\varphi}{\Lambda}\]
By the induction hypothesis we have an approximating function $\F$ such that
$\tprf \F(\Lambda)$, by \ref{weakening-rule'} we get a proof of
$\F(\Lambda),\F(\varphi)$ which is the same as $\F(\Lambda, \varphi)$,
therefore $G(n+1)=G(n)$ is enough for this case.

\textbf{\ref{v1-rule} (and \ref{v2-rule}):}
\[
\inference{\Lambda,\varphi \vee \psi}{\Lambda,\varphi}
\]
By the induction hypothesis we get an approximating function $\F_0$ such that
\[
\tprf \F_0(\Lambda,\varphi).
\]
Let $\F = \F_0 \uniform \F_{\varphi \vee \psi}$ then $\F(\Lambda,\varphi)$ is
provable and by \ref{v1-rule'}
\[\F(\Lambda),\F(\varphi) \vee \F (\psi)=\F(\Lambda,\varphi \vee
\psi)\] is provable. Therefore $G(n+1)=G(n)+1$ is enough for this case.

\textbf{\ref{v3-rule}:}
\[\inference{\Lambda,\neg (\varphi \vee \psi)}{\Lambda,\neg \varphi &
  \Lambda,\neg \psi}\] By the induction hypothesis we get approximating
functions $\F_0$ and $\F_1$ such that $\F_0(\Lambda,\neg \varphi)$ and
$\F_1(\Lambda,\neg \psi)$ are both provable. Let
\[
\F=\F_0 \uniform \F_1 \uniform \F_{\varphi \vee \psi} \uniform \F_{\neg
  (\varphi \vee \psi)} \uniform \F_{\neg \varphi} \uniform \F_{\neg \psi},
\]
by Lemma~\ref{lem:prov.template.implies.prov.aprox} $\F(\Lambda,\neg \varphi)$
and $\F(\Lambda,\neg \psi)$ are both provable. Since
\[
\F(\Lambda,\neg \varphi)= \F(\Lambda), \neg \F(\varphi) \text{ and }
\F(\Lambda,\neg \psi) =\F(\Lambda), \neg \F(\psi)
\]
we have by \ref{v3-rule'} a proof of
\[
\F(\Lambda), \neg (\F(\varphi) \vee \F(\psi)) = \F(\Lambda, \neg (\varphi \vee
\psi)).
\]
Therefore $G(n+1)=2G(n)+4$ is enough for this case.

\textbf{\ref{neg-rule}:}
\[\inference{\Lambda,\neg\neg\varphi}{\Lambda,\varphi}\]
By the induction hypothesis we get an approximating function $\F_0$ such that
\[
\F_0(\Lambda,\varphi)
\]
is provable. Let
\[
\F=\F_0\uniform \F_{\neg\varphi}\uniform \F_{\neg\neg\varphi},
\]
then
\[
\F(\Lambda,\varphi)=\F(\Lambda),\F(\varphi)
\]
is provable, so \ref{neg-rule'} gives us a proof of
\[
\F(\Lambda),\neg\neg \F(\varphi)=\F(\Lambda),\F(\neg\neg\varphi) =\F(\Lambda,
\varphi).
\]
Thus $G(n+1)=G(n)+2$ is enough for this case.

\textbf{\ref{cut-rule}:}
\[\inference{\Gamma}{\Gamma, \varphi & \Gamma, \neg \varphi}\]
By the induction hypothesis we get approximating functions $\F_0$ and $\F_1$
such that $\F_0(\Lambda,\varphi)$ and $\F_1(\Lambda,\neg \varphi)$ are
provable. Let
\[\F=\F_0 \uniform \F_1 \uniform \F_{\neg \varphi},\]
then
\begin{gather*} \F(\Gamma,\varphi)=\F(\Gamma),\F(\varphi) \quad \text{and}\\
  \F(\Gamma,\neg \varphi)=\F(\Gamma),\neg\F(\varphi)
\end{gather*}
are both provable. By \ref{cut-rule'} we get a proof of $\F(\Gamma)$. Thus
$G(n+1)=2G(n)+1$ is enough for this case.

\textbf{\ref{e-rule}:}
\[\inference{\Lambda, \exists \v_i \,\varphi}{\Lambda,\varphi[\c_a/\v_i]}\]
By the induction hypothesis we get an approximating function $\F_0$ such that
\[\F_0(\Lambda, \varphi [\c_a/\v_i])\] is provable. Let \[\F=\F_0 \uniform
\F_{\exists \v_i \,\varphi}\] then
\[
\F(\Lambda,\varphi[\c_a/\v_i])= \F(\Lambda),\F(\varphi)[\c_a/\v_i]
\]
is provable (the equality is Lemma~\ref{lem:subst.and.approx.commute}) so by
\ref{e-rule'} we get a proof of
\[
\F(\Lambda),\exists \v_i \,\F(\varphi) = \F(\Lambda),\F(\exists \v_i
\,\varphi)=\F(\Lambda, \exists \v_i \,\varphi).
\]
Thus $G(n+1)=G(n)+1$ is enough for this case.

\textbf{\ref{M-rule}:}
\[
\inference{\Lambda, \neg \exists \v_i\, \varphi}{\inferforall{a \in
    \M}{\Lambda, \neg \varphi [\c_a/\v_i]}}
\]
By the induction hypothesis we have approximating functions $\F_a$ such that
\[
\F_a(\Lambda,\neg \varphi[\c_a/\v_i])=\F_a(\Lambda), \F_a(\neg
\varphi)[\c_a/\v_i]
\]
are provable for all $a \in \M$. Let $\F'$ be as in Lemma
\ref{lem:uniform.approx.exists} and
\[
\F=\F' \uniform \F_{\exists \v_i\, \varphi} \uniform \F_{\neg \exists \v_i
  \,\varphi},
\]
then
\[
\F(\Lambda),\F(\neg \varphi)[\c_a/\v_i] = \F(\Lambda),\neg
\F(\varphi)[\c_a/\v_i]
\]
are all provable by Lemma~\ref{lem:prov.template.implies.prov.aprox}. By
\ref{M-rule'} we get a proof of
\[
\F(\Lambda),\neg \exists \v_i \,\F(\varphi) = \F(\Lambda),\F(\neg \exists
\v_i\, \varphi) = \F(\Lambda, \neg \exists \v_i\, \varphi).
\]
Thus $G(n+1)=(n+2)(2^{G(n)}-1)+2$ is enough since if $\mprf[n] \Gamma$ then
$\abs{\Gamma} \leq n+2$ and therefore $\abs{\F'} \leq (n+2)(2^{G(n)}-1)$ by
Lemma~\ref{lem:prov.template.implies.prov.aprox}.

Thus, if we define $G$ recursively by
\begin{align*}
  G(1)&=9\\
  G(n+1)&=(n+2)(2^{G(n)}-1)+2 \quad \text{for $n \geq 1$},
\end{align*}
then $G$ satisfies the proposition.
\end{proof}

\begin{defin}
  If $\Lambda$ is a set of \LM-sentences then $\apprx{\Lambda}$ is a set of
  \LT-sentences defined as
\[
\apprx{\Lambda}= \Set{\F(\lambda)| \text{$\lambda \in \Lambda$ and $\F$ is an
    approximating function}}.
\]
\end{defin}

\begin{por}\label{por:prov.mlogic.implies.template}
  If $\Lambda$ is a set of \LM-sentences and $\Lambda \mprf[\omega] \Gamma$
  then $ \apprx{\Lambda} \tprf \F(\Gamma)$ for some approximating function
  $\F$.
\end{por}

In fact we could strengthen the porism: Let $\Delta$ be any set of
\LT-sentences closed under approximating functions, i.e., if $\psi \in \Delta$
and $\F$ is an approximating function then $\F(\psi) \in \Delta$. Let
$\Lambda$ be a set of \LM-sentences and $k$ a natural numbers such that for
every $\varphi \in \Lambda$ there exists an approximating function $\F$ such
that $\abs{\F} \leq k$ and $\F(\varphi) \in \Delta$. If $\Lambda \mprf[\omega]
\Gamma$ then there exists $\F$ such that $\Delta \tprf \F(\Gamma)$.

\subsection{The consistency of \M-logic}

\begin{prop}
  If \M\ is recursively saturated and $\mprf \Gamma$ then there is an
  approximating function $\F$ such that $\tprf \F(\Gamma)$.
\end{prop}
\begin{proof}
  Combine Propositions~\ref{prop:finite.proofs.mlogic} and
  \ref{prop:prov.mlogic.implies.template}.
\end{proof}

\begin{prop}\label{prop:template.structure.implies.consistent}
  Suppose \M\ is recursively saturated, $\Lambda$ is a definable set and there
  is a template structure \T\ such that
\[
\T \models \lambda, \quad\text{for all $\lambda \in \apprx{\Lambda}$}
\]
then $\Lambda$ is consistent in \M-logic.
\end{prop}
\begin{proof}
  By Proposition~\ref{prop:finite.proofs.mlogic} if $\Lambda \mprf \emptyset$
  then $\Lambda \mprf[\omega] \emptyset$. By
  Porism~\ref{por:prov.mlogic.implies.template} there is a proof in template
  logic of $\emptyset$ from the set $\apprx{\Lambda}$, i.e., $\apprx{\Lambda}
  \tprf \emptyset$. By Proposition~\ref{prop:template.sound} we then have $\T
  \models 0 \neq 0$ which clearly is a contradiction.
\end{proof}

\begin{cor}\label{cor:mlogic.consistent}
  If \M\ is recursively saturated then \M-logic is consistent, i.e., $\nmprf
  \emptyset$.
\end{cor}

In fact, we can prove something a bit stronger by carefully examine the proofs
of Proposition~\ref{prop:finite.proofs.mlogic},
\ref{prop:prov.mlogic.implies.template} and~\ref{prop:template.sound}:

\begin{thm}\label{thm:str.cons.crit}
  Let $\Lambda$ be a set of \LM-sentences such that the expanded structure
  $\pair{\M,\Lambda}$ is recursively saturated, let $\Delta$ be a set of
  \LT-sentences closed under approximating functions and $n$ a natural number
  such that if $\varphi \in \Lambda$ then there exists an approximating
  function $\F$ such that $\abs{\F} \leq n$ and $\F(\varphi) \in \Delta$.
  Moreover let
\[
\Delta_k = \Set{ \F(\varphi) | \abs{\F} \leq k \wedge \varphi \in \Lambda
  \wedge \F(\varphi) \in \Delta}.
\]
If there exists template structures $\T_k$ such that
\[
\T_k \models \delta \quad\text{for all $\delta \in \Delta_k$}
\]
then $\Lambda$ is consistent in \M-logic.
\end{thm}

The next result is a, sort of, negative result. Usually it is expressed as
\M-logic admits full pathology.

Let $\delta_0$ be $0\neq 0$ and by induction define $\delta_{a+1}$ to be
$\delta_a \vee \delta_a$ for all $a \in \M$.

\begin{prop}
  If \M\ is recursively saturated and countable and $a \in \M \setminus
  \omega$ then \M\ admits a satisfaction class \Sat\ such that $\delta_a \in
  \Sat$.
\end{prop}
\begin{proof}
  The approximations of $\delta_a$ are
\begin{gather*}
  \tmpl{\delta_a}, \\
  \tmpl{\delta_{a-1}} \vee \tmpl{\delta_{a-1}}, \\
  (\tmpl{\delta_{a-2}} \vee \tmpl{\delta_{a-2}}) \vee (\tmpl{\delta_{a-2}}
  \vee \tmpl{\delta_{a-2}}),
\end{gather*}
and so on. The template structure $\T=\pair{\T_v,\T_t}$ with $\T_v(t)=0$ for
all closed terms $t$ and $\T_t=\Set{\delta_{a-k}: k \in \omega}$ makes all
these approximations true. Applying
Proposition~\ref{prop:template.structure.implies.consistent} gives us the
proposition.
\end{proof}

\begin{prop}
  If $t$ is a term with no constants or multiplications at finite depth, e.g.,
\[
\Sc(\Sc(\ldots \Sc(0) \ldots ))
\]
with a nonstandard number of successor symbols, and $a \in \M \setminus
\omega$ then there exists a satisfaction class $\Sat$ such that $t=\c_a \in
\Sat$.
\end{prop}
\begin{proof}
  We have to prove that there is a template structure making all the
  approximations of $t= \c_a$ true. Define $\T_v(t)=a$ and by induction define
  $\T_v$ on all closed terms occurring in $t$ at constant depth. If
  $\T_v(r)=b$ and $r$ is $\Sc(s)$ then define $\T_v(s)=b-1$ and if $r$ is
  $s+s'$ let
\[
\T_v(s)=\T_v(s')=b/2
\]
if $b$ is even and
\[
\T_v(s)=\frac{b+1}{2},\quad \T_v(s')=\frac{b-1}{2}
\]
if $b$ is odd. Let $\T_v(\c_a)=a$ and for all closed terms $r$ not occurring
in $t$ at finite depth let $\T_v(r)=0$. Finally, let $\T_t=\Set{t=\c_a}$.

It should now be clear that $\val_\T(\F(t))=a$, and so $\T \models \F(t=\c_a)$
for any approximating function $\F$. Applying
Proposition~\ref{prop:template.structure.implies.consistent} gives us the
result.
\end{proof}

The proposition is false if we allow multiplication at finite depth in $t$,
since if $t$ is, for example,
\[
\Sc(\Sc(r)) \cdot \Sc(\Sc(s))
\]
and $a \in \M$ is prime, i.e.,
\[
\M \models \Sc(0) < a \wedge \forall x,y \bigl(x \cdot y = a \imp x=\Sc(0)
\vee y=\Sc(0)\bigr),
\]
then $\mprf t \neq \c_a$.

The next result is a partial answer to Question~\ref{que:satcl.two.parts}.

\begin{prop}\label{prop:partial.answer}
  Let $\sim$ be an \LMs-definable equivalence relation on $\clterm(\M)$ and \[ E =
  \Set{ t=r | t \sim r}. \] If $\M_E$ is well-defined and the canonical map
\[
f : \M \to \M_E,\quad a \mapsto \overline{\c_a}
\]
is an isomorphism then there is a satisfaction class $\Sat$ such that
\[
t \sim r \quad\text{iff}\quad t=r \in \Sat
\]
for all closed \LM-terms $t$ and $r$.
\end{prop}
\begin{proof}
  We prove that the \LMs-definable set $E$ is consistent in \M-logic. We do
  this by defining a template structure $\T= \pair{\T_t,\T_v}$ such that all
  \LT-sentences in $\apprx{E}$ are true in $\T$.
  
  Let
\[ \T_v(t) = f^{-1}(\,\overline{t}\,) \] and $\T_t= E$. We claim that
\[ \T \models \F(t=r) \]
for closed \LM-terms $t$ and $r$ such that $t \sim r$ and any approximating
function $\F$. To see this we first observe that for any closed \LM-term $t$
and any approximating function $\F$ we have
\[ \T \models \tmpl{t} = \F(t). \]
This is proved by induction on $\abs{\F}$. For $\abs{\F}=0$ it is trivial and
for the induction step all we have to do is to observe that
\begin{multline*}
  \val_\T(\Sc(\,\tmpl{t}\,)) = \Sc^\M(\val_\T(\,\tmpl{t}\,)) =
  \Sc^\M(f^{-1}(\,\overline{t}\,)) = f^{-1}(\Sc^{\M_E}(\,\overline{t}\,)) = \\
  f^{-1}(\overline{\Sc(t)}) = \val_\T(\tmpl{\Sc(t)})
\end{multline*}
and similar for $+$ and $\cdot$. This means that when we substitute, for
example, $\Sc(\,\tmpl{t}\,)$ for $\tmpl{\Sc(t)}$ the value of $\val_\T$ does
not change. Therefore
\[
\val_\T(\F(t)) = \val_\T(\,\tmpl{t}\,)
\]
for any closed \LM-term $t$ and any approximating function $\F$.  This proves
the claim.
\end{proof}

\subsection{Some auxiliary results}

We will end this chapter by proving converse results of
Proposition~\ref{prop:template.structure.implies.consistent} and
\ref{prop:prov.mlogic.implies.template}.

\begin{prop}\label{prop:consistent.mlogic.implies.template.structure}
  If $\Lambda$ is a consistent set in \M-logic of \LM-sentences then there is
  a template structure \T\ making all sentences $\lambda \in \apprx{\Lambda}$
  true.
\end{prop}
\begin{proof}
  Since $\Lambda$ is consistent it is contained in some satisfaction class
  \Sat. Define \T\ to interpret the approximation symbols exactly as \Sat\ 
  sees them, i.e., define
\[\T_t= \Sat\] and
\[\T_v(\,\tmpl{t}\,)=a \quad\text{iff}\quad t=\c_a \in \Sat\]
It is clear that $\T \models \tmpl{\varphi}$ if $\varphi \in \Lambda$,
therefore, all we have to prove is that
\[ \T \models \psi \quad\Rightarrow \quad\T \models \F_\tau(\psi),\] for any
\LT-sentences $\psi$ and any template symbols $\tau$.

If $\tau$ is $\tmpl{\neg \varphi}$ then by observing that
\[\T \models \tmpl{\neg \varphi} \quad\text{iff}\quad
\T \models \neg \tmpl\varphi\] we see that replacing all occurrences of
$\tmpl{\neg \varphi}$ by $\neg \tmpl\varphi$ does not change the truth value
of the sentence. In the same way we have
\begin{align*}
  \T \models \tmpl{\varphi\vee \psi} &\quad\text{iff}\quad
  \T \models \tmpl{\varphi}\vee \tmpl{\psi} \qquad\text{and}\\
  \T \models \tmpl{\exists \v_i \varphi} &\quad\text{iff}\quad \T \models
  \exists \v_i \tmpl{\varphi}.
\end{align*}
Therefore, substituting $\tmpl{\varphi} \vee \tmpl{\psi}$ for $\tmpl{\varphi
  \vee \psi}$ in a sentence does not change the truth value. The same holds
for $\exists \v_i \tmpl{\varphi}$ and $\tmpl{\exists \v_i \varphi}$.

We also have that if $\tau$ is $\tmpl{\Sc(t)}$ then
\[ 
\T \models \tmpl{\Sc(t)}=\Sc(\,\tmpl{t}\,)
\] 
so substituting $\Sc(\,\tmpl{t}\,)$ for $\tmpl{\Sc(t)}$does not change the
truth value either. The same is true for addition and multiplication since
\begin{align*}
  \T &\models \tmpl{t + r}=\tmpl{t}+\tmpl{r} \quad\text{and}\\
  \T &\models \tmpl{t \cdot r}=\tmpl{t}\cdot \tmpl{r}.
\end{align*}
We have to check two more cases, the following observations will handle those:
\begin{align*}
  \T &\models \tmpl{t=r}\quad\text{iff}\quad \T \models
  \tmpl{t}=\tmpl{r} \quad \text{and} \\
  \T &\models \tmpl{\c_a}=\c_a. 
\qedhere
\end{align*}
\end{proof}

We are now in a good position to prove a converse to
Proposition~\ref{prop:prov.mlogic.implies.template}.

\begin{prop}\label{prop:prov.template.implies.mlogic}
  If $\Gamma$ is a finite set of \LM-sentences, $\F$ is any approximating
  function and $\tprf \F(\Gamma)$ then $\mprf \Gamma$.
\end{prop}
\begin{proof}
  Suppose that $\nmprf \Gamma$, then $\neg \Gamma$ is consistent in \M-logic.
  By Proposition~\ref{prop:consistent.mlogic.implies.template.structure} there
  is a template structure \T\ satisfying $\apprx{\neg \Gamma}$, but if $\tprf
  \F(\Gamma)$ then
\[
\T \models \Vee \F(\Gamma),
\]
by Proposition~\ref{prop:template.sound}. Clearly, then there is a $\gamma \in
\Gamma$ such that
\[
\T \models \F(\gamma),
\]
therefore, $\T \nmodels \neg \F(\gamma)$ and so
\[
\T \nmodels \F \uniform \F_{\neg \gamma}(\neg \gamma).
\]
Since
\[
\F \uniform \F_{\neg \gamma}(\neg \gamma) \in \apprx{\neg \Gamma}
\]
this yields a contradiction.
\end{proof}

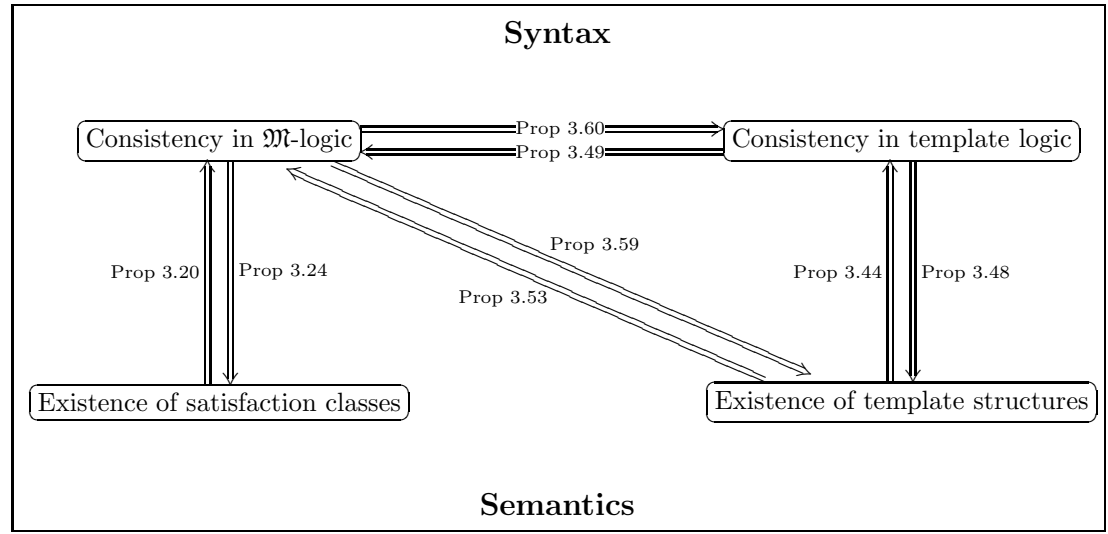
\begin{sidewaysfigure}
\begin{center}
  \fbox{ $\xymatrix{
      & \text{\large\bf Syntax} & \\
      *+[F-:<3pt>]{\text{Consistency in \M-logic}}
      \ar@<1ex>@{=>}[rr]|{\text{Prop \ref{prop:prov.template.implies.mlogic}}}
      \ar@<1ex>@{=>}[ddd]^{\text{Prop \ref{prop:sound.complete.mlogic}}}
      \ar@<0ex>@{=>}[]+<10ex,-2ex>;[dddrr]+<-8ex,2.5ex>^{\raisebox{4pt}{\text{%
            \scriptsize Prop
            \ref{prop:consistent.mlogic.implies.template.structure}}}} & &
        *+[F-:<3pt>]{\text{Consistency in template logic}}
        \ar@<1ex>@{=>}[ll]|{\text{\strut Prop
            \ref{prop:prov.mlogic.implies.template}}}
        \ar@<1ex>@{=>}[ddd]^{\text{\strut Prop \ref{prop:template.complete}}} \\
        \hbox{} & & \\ & & \hbox{}
        \\
        *+[F-:<3pt>]{\text{Existence of satisfaction classes}}
        \ar@<1ex>@{=>}[uuu]^{\text{\strut Prop \ref {prop:sound.mlogic}}}& &
        *+[F-:<3pt>]{\txt{Existence of template structures}}
        \ar@<0ex>@{=>}[]+<-12ex,2ex>;[lluuu]+<6ex,-2.5ex>^{\raisebox{-8pt}{\text{%
              \scriptsize Prop 
              \ref{prop:template.structure.implies.consistent}}}}
        \ar@<1ex>@{=>}[uuu]^{\text{\strut Prop \ref {prop:template.sound}}}
        \\
        & \text{\large\bf Semantics}& }$}
\end{center}
\caption{A diagram showing the main results of Chapter~\ref{chp:satcl}.}\label{fig:tva}
\end{sidewaysfigure}

To sum up we illustrate the chapter by Figure~\ref{fig:tva}.

\chapter{Weaker satisfaction classes}\label{chp:wesatcl}

In this chapter we will study sets which fail to be satisfaction classes, but
just merely; \ref{Axiom:12} might be false in it, i.e., there might be terms
$t$ such that
\[
\neg \exists \v_0 (t = \v_0)
\]
is in the set.

\section{Free \M-logic}

We will try to answer the question:
\begin{quote}
  What happens if we remove \ref{Axiom:12} from the axioms of \M-logic?
\end{quote}
Let us call a set of \LM-sentences satisfying the alternative definition of
satisfaction class given in Proposition~\ref{prop:alt.satcl.def} but with the
word `isomorphism' changed to `isomorphic embedding' for a {\em free}
satisfaction class. And \M-logic without \ref{Axiom:12} for {\em free}
\M-logic.\footnote{See Remark~\ref{rem:free} for an explanation of the name.}

\begin{prop}\label{prop:free.satcl.equiv.free.mlogic}
  Free \M-logic corresponds to free satisfaction classes in the same way as
  \M-logic corresponds to satisfaction classes, i.e., every maximally
  consistent set of sentences in free \M-logic is a free satisfaction class
  and every free satisfaction class is a maximally consistent set of
  sentences.
\end{prop}
\begin{proof}
  Assume $X$ is a maximally consistent set in free \M-logic. We have to prove
  that $X$ satisfies the definition of a free satisfaction class, the only
  nontrivial parts are to prove that $\sim_X$ is an equivalence relation, that
  $\M_X$ is well-defined and that the canonical map $f : \M \to \M_X$ is an
  isomorphic embedding.
  
  Clearly, \ref{Axiom:3}, \ref{Axiom:4} and \ref{Axiom:5} implies that
  $\sim_X$ is an equivalence relation, furthermore \ref{Axiom:6},
  \ref{Axiom:7} and \ref{Axiom:8} implies that $\M_X$ is well-defined. That
  $f$ is a homomorphism follows from \ref{Axiom:9}, \ref{Axiom:10} and
  \ref{Axiom:11}. Finally, the injectivity follows from \ref{Axiom:2}.
  
  For the converse, assume that \Sat\ is a free satisfaction class.  \Sat\ is
  clearly maximally consistent if it is consistent since for every
  \LM-sentence $\varphi$ either $\varphi \in \Sat$ or $\neg \varphi \in \Sat$.
  Thus, all we have to prove is that $\Sat$ is consistent in free \M-logic,
  i.e., we have to check that if $\Delta$ is an axiom then $\Vee \Delta \in
  \Sat$ and that $\Sat$ is closed under all inference rules; the consistency
  then follows from the fact that there are \LM-sentences not in $\Sat$.
  
  From the fact that $\sim_\Sat$ is an equivalence relation it is easy to see
  that the disjunctions of \ref{Axiom:3}, \ref{Axiom:4} and \ref{Axiom:5} all
  are in \Sat. For \ref{Axiom:6}; assume that $t=r \in \Sat$, by the
  well-definition of $\M_\Sat$
\[
\overline{\Sc(t)}=\Sc^{\M_\Sat}(\,\overline{t}\,)=\Sc^{\M_\Sat}(\overline{r})=
\overline{\Sc(r)},
\]
thus, $\Sc(t)=\Sc(r) \in \Sat$. Similar for \ref{Axiom:7} and \ref{Axiom:8}.
Clearly, the disjunction of \ref{Axiom:1} is in \Sat, and if $a \neq b$ then
\[ \overline{\c_a} = f(a) \neq f(b) = \overline{\c_b} \]
so $\c_a = \c_b \notin \Sat$. Therefore, by the maximality of $\Sat$ we have
$\c_a \neq \c_b \in \Sat$.

Furthermore,
\[
\overline{\c_{a+^\M b}} = f(a +^\M b) = f(a) +^{\M_\Sat} f(b) =
\overline{\c_a} +^{\M_\Sat} \overline{\c_b} = \overline{\c_a + \c_b}.
\]
Thus $\c_{a+^\M b}=\c_a + \c_b \in \Sat$ and similarly for \ref{Axiom:9} and
\ref{Axiom:11}.

That \Sat\ is closed under the inference rules is proved by using properties
(\ref{s10}), (\ref{s11}) and (\ref{s12}) of \Sat, it is left to the reader.
\end{proof}

A natural question now arises:

\begin{quote}
  Are there free satisfaction classes which are not satisfaction classes?
\end{quote}
The answer is yes as we now will prove. We prove that if $a \in \M \setminus
\omega$ then we can find a free satisfaction class \Sat\ such that $\neg
\exists \v_0 (\num{a} = \v_0) \in \Sat$, where $\num{a}$ is the closed term
defined inductively as follows:
\[
\num{a} =
\begin{cases}
  0 & \quad\text{if $a=0$} \\
  \Sc(\num{b}) & \quad\text{if $a=\Sc^\M(b)$}.
\end{cases}
\]

The consistency criteria we worked out in Chapter~\ref{chp:satcl} is too acute
to handle this since in every \LT-structure $\T$ we have
\[
\T \models \exists \v_0 (\F(\num{a}) = \v_0)
\]
for any approximating function $\F$. The solution to this problem is to
redefine and make \LT-structures more general.

\begin{defin}
  A {\em free} \LT-structure, $\T$, is a pair, $\pair{\T_t,\T_v}$, of a set, $\T_t$,
  of \LM-sentences and a map, $\T_v$, from the closed \LM-terms into some
  \LA-structure $\structure{N} \supseteq \M$.
\end{defin}

Truth in a free \LT-structure is defined in the obvious way.

We state the consistency criteria in one of its simplest forms, but it should
be evident that it could be strengthen as in Theorem~\ref{thm:str.cons.crit}.

\begin{prop}
  Let $\varphi$ be a \LM-sentence and $\F_0$ an approximating function such
  that there exists a free \LT-structure making $\F \circ \F_0(\varphi)$ true
  for any approximating function $\F$.  Then $\neg \varphi$ cannot be proved in
  free \M-logic.
\end{prop}
\begin{proof}
  Only some small modifications to the proofs of
  Proposition~\ref{prop:template.sound}
  and~\ref{prop:prov.mlogic.implies.template} is needed. The details are left
  to the reader.
\end{proof}

Let $\structure{N} \varsupsetneq \M$ and $b \in \structure{N} \setminus \M$.
Define
\begin{align*}
  \T_v(\num{a-k})&\eqdef b-k \quad\text{for all $k \in \omega$,}\\
  \T_t &\eqdef \emptyset,\\
  \T &\eqdef \pair{\T_t, \T_v}  \qquad\text{and}\\
  \F_0 &\eqdef \F_{\v_0} \circ \F_{\num{a} = \v_0} \circ \F_{\exists \v_0
    (\num{a} = \v_0)} \circ \F_{\neg \exists \v_0 (\num{a} = \v_0)}.
\end{align*}
Clearly
\[
\T \models \F \circ \F_0 (\neg \exists \v_0 (\num{a} = \v_0))
\]
for every approximating function $\F$, since $\val_\T(\F(\num{a}))=b$ for
every $\F$. Therefore, by the proposition, $\neg \exists \v_0 (\num{a} =
\v_0)$ is consistent in free \M-logic. By the usual construction we can find a
maximally consistent set \Sat\ including $\neg \exists \v_0 (\num{a} = \v_0)$
which by Proposition~\ref{prop:free.satcl.equiv.free.mlogic} is a free
satisfaction class. \Sat\ is not a satisfaction class since \ref{Axiom:12} is
not true.

\begin{rem}
Even though a free satisfaction class \Sat\ is a satisfaction class if for
every closed \LM-term $t$ there exists $a \in \M$ such that $t=\c_a \in \Sat$
and we get an equivalent definition of satisfaction class if we replace
(\ref{s9}) in Definition \ref{def:satcl} by 
\[
\exists x \, \godel{t=\c_x} \in \Sat
\]
for every closed \LM-term $t$, we do {\em not\/} get an equivalent definition
of a free satisfaction class by removing (\ref{s9}) in that definition. This
is a consequence of the fact that the \LM-sentence
\begin{equation}\label{eq:sent}
\Sc(\c_a)=\Sc(\num{b}) \wedge \c_a \neq \num{b} 
\end{equation}
is consistent in free \M-logic for any $a,b \in \M$ such that $b > \omega$,
and therefore included in some free satisfaction class. But (\ref{s6}) implies
that for any satisfaction class $\Sat$ in which \ref{s9} may fail the sentence
(\ref{eq:sent}) is false since if $\Sc(\c_a)=\Sc(\num{b}) \in \Sat$ then
$\c_{\Sc^\M(a)} = \Sc(\num{b}) \in \Sat$, thus (\ref{s6}) implies that $\c_a =
\num{b} \in \Sat$. To prove that the sentence (\ref{eq:sent}) is consistent in
free \M-logic it is enough to construct an \LA-structure $\structure{N}
\varsupsetneq \M $ such that for some $d \in \structure{N}$ we have
$\Sc^\structure{N}(d) = \Sc^\M(a)$ and $d \neq a$. Then define an
\LT-structure mapping $\num{b-k}$ to $d-k$ for every $k \in \omega$.
\end{rem}

\begin{rem} \label{rem:free}
  We are using the term {\em free} since free \M-logic is a sort of (positive)
  free logic, see \cite{Lambert:01}. In fact the part of free \M-logic where
  we only consider sentences of the form
  $\sigma[t_1,\ldots,t_k/\v_{i_1},\ldots,\v_{i_k}]$, where $\sigma$ is a
  \LMs-formula and $t_1$,~\dots,~$t_k$ are closed \LM-terms, is a (positive)
  free logic, with the existential predicate defined as
\[
\mathrm{E}! \, t \ekv_\text{def} \exists \v_0 (t=\v_0).
\]
In fact, for any $\varphi(x)$ of this form and any closed \LM-terms $t$ and
$r$ we can prove (in this restricted free \M-logic)
\[
\varphi(t) \wedge t=r \imp \varphi(r).
\]
Therefore, it is easy to see that
\[
\varphi(t) \wedge \mathrm{E}! \, t \imp \exists x \varphi(x)
\]
also is provable (in the same logic). See \cite{Lambert:01} or
\cite{Bencivenga:99} for more information on free logics.
\end{rem}

\question{Which \LA-structures $\structure{N} \supseteq \M$ are $\M_\Sat$ for
  some free satisfaction class \Sat?}{It is a natural question to ask.  It
  might be the case that the structures $\M_\Sat$ have very specific
  properties, analogue to the models arising in the arithmetised completeness
  theorem.}

\question{Are there free satisfaction classes \Sat, which are not satisfaction
  classes, such that the canonical map $f: \M \to \M_\Sat$ is an elementary
  embedding?}{If this is true, is there a corresponding extension of free \M-logic?}

\chapter{Stronger satisfaction classes}
\label{chp:stsatcl}

As we have seen in Chapter~\ref{chp:satcl} some ``pathological'' examples
arise in the study of satisfaction classes. For example, we can make the
sentences $\delta_a$ and $\num{a}=\c_b$ true in a satisfaction class (if $a$,
$b$ are nonstandard).\footnote{Remember that $\delta_0$ is $0 \neq 0$,
  $\delta_{a+1}$ is $\delta_a \vee \delta_a$, $\num{0}$ is $0$ and
  $\num{\Sc^\M(a)}$ is $\Sc(\num{a})$.} The main question we will try to
answer in this chapter is:
\begin{quote}
  What do we need to remove such ``pathological'' examples?
\end{quote}

To answer this question we have concentrated on extensions of \M-logic. Any
maximally consistent set of sentences in any of the extensions we will study
in this chapter is a satisfaction class removing pathologies of a certain
kind.

It should be remarked that this chapter is included to emphasise the vast
amount of open questions in this area. There is a lot of work to be done, I
have just scratched the surface. The big question of consistency of the
extensions is a hard question. We know the answer to some of them but not to
all, but a small remark is in order here:

\begin{rem}
  Since every logic we will study have axioms and inference rules definable in
  \PA, consistency is a $\Sigma_1^1$ statement:
\begin{multline*}
  \exists X \exists \varphi \bigl( \sent(\varphi) \wedge \varphi \notin X
  \wedge \text{$X$ includes all axioms} \\ {}\wedge \text{$X$ is closed under
    all inference rules}\bigr).
\end{multline*}
Therefore, in countable models, consistency of one of these extensions of
\M-logic could only depend on $\Th(\M)$ and not on other model theoretic
properties of \M, such as saturation properties, since if \M-logic is
consistent then \M\ is recursively saturated by Theorem~\ref{thm:lachlan} and so
resplendent by Theorem~\ref{thm:sat.imp.res}.
\end{rem}

\begin{figure}
\begin{center}\begin{tabular}{|c|c|c|c|}
\hline {\bf Pathology} & {\bf Solution} & {\bf Consistent?}  & {\bf Section} \\
\hline $\num{a}=\c_b$ & Axioms: $\Tr_\text{At}$ & Yes & \ref{sec:par.tru.def}\\
\hline $\delta_a$ & Axioms:  $\Tr_{\Delta_0}$  &Yes & \ref{sec:par.tru.def}\\
\hline $\epsilon_a^\varphi$ & Rule:  \ref{prop-rule}  & ? & \ref{sec:prop}\\
\hline $\mathop{\exists \v_0 \v_1 \ldots \v_a} 0 \neq 0$ & Axioms:
$\Tr_{\Sigma_1}$
            & Yes & \ref{sec:par.tru.def}\\
\hline & Rules: \ref{einf-rule}, \ref{Minf-rule}  & ? & \ref{sec:infinite}\\
\cline{2-4} \raisebox{1.5ex}[0mm][0mm]{$\bigl(\mathop{\exists \v_0
\v_1 \ldots
            \v_a} \varphi\bigr) \ekv \neg \varphi$} &Rule:  \ref{pred-rule}
            & ? & \ref{sec:predicate}\\
\hline & Rule: \ref{skolem-rule} & ? & \ref{sec:robinson}\\
\cline{2-4} \raisebox{1.5ex}[0mm][0mm]{$\bigl(\mathop{\exists \v_0
\forall \v_1
            \ldots \exists\v_{2a}} \varphi\bigr) \ekv \neg \varphi$} &
            Rule: \ref{pred-rule} & ? & \ref{sec:predicate}\\
\hline
\end{tabular}\end{center}
\caption{A summary of the pathological examples we will study in
  this chapter and their ``solutions.''}\label{fig:pat}
\end{figure}

Figure~\ref{fig:pat} is a summary of some of the pathologies and their
``solutions.'' The sentences $\epsilon_a^\varphi$ are defined as follows:
\begin{align*}
  \epsilon_0^\varphi \quad &\text{is} \quad \neg (\varphi \vee \neg
  \varphi)\qquad \text{and} \\
  \epsilon_{a+1}^\varphi \quad &\text{is} \quad \epsilon_a \vee \epsilon_a.
\end{align*}
In the table, $a$ is nonstandard and $\varphi$ is a sentence of high
complexity (it is not $\Sigma_k$ for any $k \in \omega$). The `Solution' to a
pathology tells us what we need to add to \M-logic (either axioms of inference
rules) to remove the pathology, i.e., to be able to prove the negation of the
pathology. The `Consistent?' column tells us if this logic is consistent or
not; a `Yes' means that in any recursively saturated model the logic is
consistent and a `?' means that we do not know the answer. The column named
`Section' is a reference for where to read more about the solution, it is the
section number in this chapter.

\question{How should Figure~\ref{fig:pat} be completed?}{The question marks
  are all in the column `Consistent?', thus this is a question of proving
  consistency of extensions of \M-logic. The only tool we have to do so is
  Theorem~\ref{thm:str.cons.crit} but it will not help us in this situation.
  The only plausible approach we have found is to alter the definition of
  template logic and template structure, but every attempt of this has ended
  with tears. We think this question is very hard.}
\label{que:chapter.five}

\section{Partial Truth Definitions}\label{sec:par.tru.def}

In \cite{Kaye:91*2} a satisfaction class is defined to extend the set
\[\Tr_\text{At} = \Set{t=r | \M \models \val(t)=\val(r)}.\]
By Proposition~\ref{prop:template.structure.implies.consistent} this set is
consistent since we can define a template structure $\T=\pair{\T_t,\T_v}$,
with
\begin{align*}
  \T_t&=\Set{t=r | \M \models \val(t)=\val(r)} \quad \text{and} \\
  \T_v&= \val,
\end{align*}
where $\val$ is the valuation function definable in \PA. This is a template
structure making all \LT-sentences in $\apprx{\Tr_\text{At}}$ true. In the
same manner we can find template structures making all sentences in
$\apprx{\Tr_{\Sigma_k}}$ true.

\begin{prop}
  The sets $\Tr_\text{At}, \Tr_{\Delta_k}, \Tr_{\Sigma_k}$ and $\Tr_{\Pi_k}$
  are all consistent.
\end{prop}
\begin{proof}
  It is clearly enough to prove that $\Tr_{\Sigma_k}$ is consistent for any $k
  \in \omega$. The sets are definable so, by
  Proposition~\ref{prop:template.structure.implies.consistent}, it suffices to
  find a template structure making the \LT-sentences in
  $\apprx{\Tr_{\Sigma_k}}$ true.
  
  Let $\T=\pair{\T_v,\T_t}$ be such that $\T_v(t)=a$ iff $\M \models
  \val(t)=\c_a$ and let $\T_t=\Tr_{\Sigma_k}$. By an easy induction it is easy
  to see that this structure satisfies the condition. The induction is left to
  the reader, but we remark that the properties of $\val$ and $\Tr_{\Sigma_k}$
  in Section~\ref{sec:part.tr.def} are used heavily.
\end{proof}

\section{Closure under propositional logic}
\label{sec:prop}

Satisfaction classes closed under nonstandard propositional proofs, in the
sense that if \M\ thinks $\varphi$ is provable in propositional logic from
sentences in \Sat\ then $\varphi \in \Sat$, is the next object of study.
Firstly, we have to define what it means for a model to think something is
provable in propositional logic, i.e., we need some formula expressing
propositional provability.

\begin{defin}
  If $\sigma(x)$ is any formula then $\forpropprf{\sigma}(y)$ is defined to be
  the formula
\begin{multline*}
  \exists x \Bigl( [x]_{\len(x)-1}=y \wedge \forall i {<} \len(x) \bigl[
  \sent([x]_i) \wedge \bigl[\Ax([x]_i) \vee \sigma([x]_i) \\ \vee \exists j,k
  {<} i \exists z \bigl([x]_j= \godel{z \imp [x]_i} \wedge
  [x]_k=z\bigr)\bigr]\bigr]\Bigr),
\end{multline*}
where $\Ax(x)$ is a formula defining the axioms of propositional
logic.\footnote{The axioms could be chosen in a variety of ways.  Use your
  favourite axiomatisation.}
\end{defin}

The formula $\forpropprf{\sigma}(y)$ says that there exists a sequence of
sentences such that every element in the sequence is a \LM-sentence $\varphi$
and either an axiom, satisfying $\sigma$, or a result of applying Modus Ponens
to other sentences occurring in the sequence prior to $\varphi$.

This section is about satisfaction classes closed under this relation, in the
sense that
\[
\M \models \forall x \bigl(\forpropprf{x \in \Sat}(x) \imp x\in \Sat\bigr).
\]

For simplicity we will write $\Lambda \propprf \varphi$ to mean $\M \models
\forpropprf{x \in \Lambda} (\varphi)$, but please do remember that all
propositional proofs are ``inside'' the model \M.

It is important to observe that we have a sort of compactness theorem even for
nonstandard proofs:

\begin{prop}\label{prop:prop.logic.compact}
  Let $\M^+$ be any expansion of \M, $\sigma(x)$ a formula in the language of
  $\M^+$ and $\sigma \propprf \varphi$, then there exists an \LMs-formula
  $\delta(x)$ such that
\[
\M^+ \models \forall x \bigl(\delta(x) \imp \sigma(x)\bigr) \wedge \exists x
\forall y \bigl(\delta(y) \imp y < x\bigr) \wedge
\forpropprf{\delta}(\varphi).
\]
\end{prop}
\begin{proof}
  Let $p$ be a proof of $\varphi$ from $\sigma$, i.e., a witness for the
  existential sentence $\forpropprf{\sigma}(\varphi)$. Define $\delta(x)$ to
  be
\[
\exists i \mathord< \len(p) \Bigl([p]_i=x \wedge\neg \Ax(x) \wedge \forall j,k
\mathord< \len(p) \forall y \bigl([p]_j=\godel{y \imp x} \imp [p]_k \neq
y\bigr)\Bigr),
\]
saying that $x$ is a sentence in the proof $p$ but it is not a result of Modus
Ponens, neither an axiom of propositional logic.  Clearly, this \LMs-formula
has the desired property.
\end{proof}

Define $\Vee \Gamma$, for any finite set of \LM-sentences $\Gamma$, in the
following way: let $\Vee \Gamma$ be the \LM-sentence
\[
\gamma_0 \vee \bigl(\gamma_2 \vee \bigl( \gamma_3 \vee \bigl(\ldots
\gamma_{k-1}\bigr)\bigr)\bigr)
\]
where $\gamma_i < \gamma_{i+1}$ for all $i<k-1$ and $\gamma_0$, \dots,
$\gamma_{k-1}$ enumerates the set $\Gamma$.

We will define a new logic, extending \M-logic, that corresponds to
satisfaction classes closed under propositional proofs. We call it
\Mprop-logic and it is similar to \M-logic, with the important difference that
\ref{weakening-rule}, \ref{v1-rule}, \ref{v2-rule}, \ref{v3-rule},
\ref{neg-rule} and \ref{cut-rule} are replaced by the single rule:

\begin{equation}
\tag{Prop}\label{prop-rule} \inference{\Gamma}{\inferforall{i <
a}{\Lambda,\varphi_i}} \qquad \text{if $\Set{\Vee \Lambda \vee
\varphi_0 ,\ldots, \Vee \Lambda \vee \varphi_{a-1}} \propprf \Vee
\Gamma$,}
\end{equation}
where $a \in \M$ and $\varphi_i$, for $i < a$, are \LM-sentences.

This means that the inference rules of \Mprop-logic are \ref{prop-rule},
\ref{e-rule} and \ref{M-rule}. Please observe that we still restrict the sets
of \LM-sentences we derive to be actually finite.

It is easy to see that \ref{weakening-rule}, \ref{v1-rule}, \ref{v2-rule},
\ref{v3-rule}, \ref{neg-rule} and \ref{cut-rule} are all derivable from the
single rule \ref{prop-rule}, therefore, $\Mprop$-logic is an extension of
\M-logic.

Let $\mpropprf$ denote provability in \Mprop-logic and let $\mpropprf[\alpha]$
mean provability with a proof of height less than $\alpha$, analogous to the
$\mprf$ and $\mprf[\alpha]$ relations.

We will show that \Mprop-logic is actually ``finite'' in the sense that if
something is provable then it is provable by a finite height proof (as
\M-logic also is). This will be true for all logics studied in this chapter,
but since all proofs follow the same line, we will only prove it for
\Mprop-logic.

We define formulas $\Pf_k'(x)$ as follows. Let $\Pf_1'(x)$ be the formula
defining the axioms of \Mprop-logic, i.e., $\M \models \Pf_1'(a)$ iff $a$ is
$\Vee \Gamma$ for $\Gamma$ an axiom of \Mprop-logic. Furthermore, let
\begin{align*}
  \Pf_{k+1}'(x)=& \forpropprf{\Pf_{k}'}(x) \\
  & {}\vee \exists y,z,i \bigl(\sent(y) \wedge \sent(\godel{\mathop{\exists
      \v_i} z}) \wedge
  x=\godel{\mathop{\exists \v_i} z \vee y} \\
  & \quad {}\wedge
  \mathop{\exists y} \Pf_{k}'(\godel{z[\c_y/\v_i] \vee y})\bigr) \\
  & {}\vee \exists y,z,i \bigl(\sent(y) \wedge \sent (\godel{\mathop{\exists
      \v_i} z}) \wedge x=\godel{\neg
    \exists \v_i \, z \vee y}\\
  & \quad {} \wedge \forall y \, \Pf_{k}'(\godel{\neg z[\c_y/\v_i] \vee
    y})\bigr),
\end{align*}
for $k > 1$.

These formulas ``code'' proofs of finite height in \Mprop-logic in the
following way:

\begin{lem} If $\Gamma$ is a finite set of \LM-sentences, $\varphi$
  is $\Vee \Gamma$ and $k>0$, then
\begin{align*}
  \mpropprf[k]\Gamma \quad &\Longrightarrow \quad \M \models
  \Pf_k'(\varphi) \qquad \text{and} \\
  \M \models \Pf_k'(\varphi) \quad &\Longrightarrow \quad \mpropprf[3k-2]
  \varphi.
\end{align*}
\end{lem}
\begin{proof}
  The proof is by induction on $k$. For the base case, assume that
  $\mpropprf[1] \Gamma$ then $\Gamma$ is an axiom of \Mprop-logic and $\M
  \models \Pf_1'(\varphi)$. On the other hand, if $\M \models \Pf_1'(\varphi)$
  then, let $\Lambda$ be the axiom such that $\Vee \Lambda$ is $\varphi$,
  since
\[
\Vee \Lambda \propprf \varphi,
\]
we have $\mpropprf[2] \varphi$.

For the induction step assume that the lemma holds for $k \leq n$.  We will
only prove the \ref{prop-rule} cases, the \ref{e-rule} and the \ref{M-rule}
cases are left to the reader.

Assume that $\mpropprf[n+1] \Gamma$ and that the last inference in the proof
of $\Gamma$ is \ref{prop-rule}, then there is a finite set $\Lambda$ and
\LM-sentences $\varphi_0,\ldots, \varphi_{a-1}$ such that
\[
\Set{ \Vee \Lambda \vee \varphi_0, \ldots, \Vee \Lambda \vee \varphi_{a-1}}
\propprf \Vee \Gamma,
\]
and $\mpropprf[n] \Lambda, \varphi_i$ for all $i < a$. By the induction
hypothesis,
\[
\M \models \Pf_{n}'(\Vee \Lambda, \varphi_i),
\]
so, clearly,
\[
\Pf_n' \propprf \Vee \Lambda \vee \varphi_i,
\]
for all $i < a$. Thus, $\Pf_{n}' \propprf \varphi$, and so $\M \models
\Pf_{n+1}'(\varphi)$.

On the other hand, if $\M \models \Pf_{n+1}'(\varphi)$ and
\[
\M \models \Pf_{n}' \propprf \varphi,
\]
then there exists, by Proposition~\ref{prop:prop.logic.compact}, \LM-sentences
$\varphi_0, \ldots, \varphi_{b-1}$ such that
\[
\M \models \Pf_{n}'(\varphi_i)
\]
for all $i < b$ and
\[
\Set{\varphi_i}_{i < b} \propprf\varphi.
\]
Therefore, by the induction hypothesis,
\[
\mpropprf[3n-2]\; \varphi_i
\]
for all $i < b$, and so, by using \ref{prop-rule},
\[
\mpropprf[3n-1]\; \varphi. 
\qedhere
\]
\end{proof}

The lemma tells us that
\[
\mpropprf[\omega] \Gamma \quad \text{iff} \quad \text{there exists $k \in
  \omega$ such that $\M \models \Pf_k'(\Vee \Gamma)$},
\]
for any finite set $\Gamma$ of \LM-sentences.

\begin{lem}
  If \M\ is recursively saturated and $\mpropprf \Gamma$ then
  $\mpropprf[\omega] \Gamma$.
\end{lem}
\begin{proof}
  Assume the lemma is false and let $\Gamma$ be such that
\[
\mpropprf[\omega+1] \Gamma \quad \text{but} \quad \nmpropprf[\omega] \Gamma,
\]
we can find such $\Gamma$ as in Lemma~\ref{prop:finite.proofs.mlogic}. Either
the last rule in the proof of height $\omega$ of $\Gamma$ is \ref{M-rule} or
\ref{prop-rule}.  For the first case the last inference is
\[
\infer{\Delta, \neg \exists \v_i \, \psi}{\inferforall{a \in \M}{\Delta,\neg
    \psi[\c_a/\v_i]}}.
\]
Define the type
\[p(x)=\Set{\neg \Pf_k'(\godel{\neg\psi[\c_x/\v_i] \vee \varphi}) | k \in
  \omega},\] where $\varphi$ is $\Vee \Delta$.

For the second case the last inference is
\[
\inference{\Gamma}{\inferforall{i < a}{\Lambda, \varphi_i}},
\]
where
\[
\Set{\psi \vee \varphi_0, \ldots, \psi \vee \varphi_{a-1}} \propprf \varphi,
\]
$\psi$ is $\Vee \Lambda$, $\varphi$ is $\Vee \Gamma$ and $a$ is nonstandard. By
Proposition~\ref{prop:prop.logic.compact} there is a \LMs-definable subset
$\Delta$ of $\Set{\psi \vee \varphi_i}_{i<a}$ such that
\[
\Delta \propprf \varphi
\]
Let $b \in \M$ enumerate $\Delta$ such that
\[
\Set{[b]_i | i < \len b} = \Delta
\]
and define the type
\[
p(x)=\Set{\neg \Pf_k'(\godel{\psi \vee (b)_x}) | k \in \omega} \cup \Set{x <
  \len b}.
\]
In either case $p(x)$ is a non-realized recursive type contradicting the
recursive saturation of \M.
\end{proof}

The Lemma tells us that to prove the consistency of \Mprop-logic we only need
to prove that
\[ \mathop{\forall k \mathord\in \omega}  \M \models \neg \Pf'_k(\godel{0 \neq 0}). \]
We have, however, not succeeding in doing so.\footnote{In a preliminary draft
  of this thesis a proof of the consistency was presented, but it turned out
  to be erroneous.}

\begin{rem}
  It is easy to see that there are models \M\ where \Mprop-logic is
  consistent. Let $\Nat$ be the standard model of \PA\ and let $\Sat_0$ be the
  standard satisfaction class on \Nat, i.e.,
\[
\Sat_0 = \Set{\godel{\varphi} | \varphi \in \ElDiag(\Nat)}.
\]
Clearly, $\Sat_0$ is closed under propositional proofs; thus,
\[
\pair{\Nat,\Sat_0} \models \satcl(\Sat_0) \wedge \forall x \bigl(\forpropprf{x
  \in \Sat_0}(x) \imp x \in \Sat_0\bigr).
\]
Therefore, if
\[
\pair{\Nat,\Sat_0} \prec \pair{\M,\Sat}
\]
then $\Sat$ is a satisfaction class closed under propositional proofs.
\end{rem}

Let us now instead prove that \Mprop-logic corresponds to satisfaction classes
closed under propositional proofs.

\begin{prop}\label{prop:propsatcl.equiv.mproplogic}
  Satisfaction classes closed under propositional proofs are exactly the
  maximally consistent sets in \Mprop-logic.
\end{prop}
\begin{proof}
  Let $\Sat$ be a maximally consistent set. It is easy to see that it is a
  satisfaction class, just as we did in Chapter~\ref{chp:satcl}. We check that
  it is closed under propositional logic. Suppose $\Sat \propprf \varphi$,
  then by the maximality either $\varphi$ or $\neg \varphi$ is in \Sat. If
  $\neg \varphi \in \Sat$ then $\Sat \mpropprf \emptyset$, therefore, since
  \Sat\ is consistent, $\varphi \in \Sat$.
  
  If \Sat\ is a satisfaction class closed under propositional proofs then it
  is consistent in \Mprop-logic, since it is closed under \ref{prop-rule},
  \ref{e-rule} and \ref{M-rule}. It is maximally consistent since for every
  $\varphi$ either $\varphi \in \Sat$ or $\neg \varphi \in \Sat$.
\end{proof}

\section{Infinite \ref{e-rule} and \ref{M-rule}}\label{sec:infinite}

The inference rule \ref{prop-rule} handles the propositional connectives in a
satisfying way. In a first try to handle quantifiers we add two infinite
versions of \ref{e-rule} and \ref{M-rule}:

\begin{equation}
\tag{I$\exists^\infty$}\label{einf-rule} \inference{\Gamma,
\mathop{\exists \v_{i_0},  \v_{i_1},\ldots, \v_{i_b}}
\varphi}{\Gamma,\varphi[\c_{(a)_0}, \c_{(a)_1}, \ldots,
\c_{(a)_{b}}/\v_{i_0}, \v_{i_1}, \ldots, \v_{i_b}]}
\end{equation}
and
\begin{equation}
\tag{$\M^\infty$-rule}\label{Minf-rule} \inference{\Gamma,\neg
{\exists \v_{i_0}, \ldots, \v_{i_b}}\, \varphi}{\inferforall{a \in
\M}{\Gamma,\neg \varphi[\c_{(a)_0}, \c_{(a)_1},\ldots,
\c_{(a)_b}/\v_{i_0},\v_{i_1}, \ldots, \v_{i_b}]}}
\end{equation}
where $b \in \M$ may be nonstandard.

\begin{defin}[\cite{Krajewski:76}] A satisfaction class \Sat\ is
  {\em $\exists$-complete} if
\begin{multline*}
  \mathop{\exists \v_{i_0}, \ldots, \v_{i_b}} \varphi \in \Sat
  \quad\text{iff} \quad \text{there exists $a \in \M$ such that}\\
  \varphi[\c_{(a)_0}, \ldots, \c_{(a)_b}/\v_{i_0}, \ldots, \v_{i_b}] \in \Sat.
\end{multline*}
\end{defin}

The following proposition should now be easy to prove.

\begin{prop}
  A set of \LM-sentences is a maximally consistent set in \M-logic with
  \ref{einf-rule} and \ref{Minf-rule} added iff it is a $\exists$-complete
  satisfaction class.
\end{prop}
\begin{proof}
  Left to the reader.
\end{proof}

%

\section{Skolem operators}\label{sec:robinson}

It seems that we could extend the infinite quantifier rules even more. To be
able to state this extended quantifier rule we need the notion of Skolem
operators.

Let $\Q$ be a \LM-definable sequence of quantifiers, i.e., $\Q \in \M$ such
that
\[
\forall j \mathord< \len(\Q) \exists i \bigl([\Q]_j=\godel{\exists \v_i} \vee
[\Q]_j=\godel{\forall \v_i}\bigr).
\]
Define a function $f_\Q : \M \mapsto \M$ such that the $f_\Q(j)$th
$\forall$-quantifier in $\Q$ is the first $\forall$-quantifier preceding (to
the left of) the $j$th $\exists$-quantifier in $\Q$.  More formally; there is
a \PA-definable function $F(x,y)$ such that the following is provable in \PA:
\begin{align*}
  F(x,y) &= 0 \quad \text{if $x$ is not a sequence of quantifiers.}\\
  F([],y) &=0 \\
  F(x,0) &= 0 \\
  F([\godel{\exists \v_i}] \conc x,\Sc(y)) &= F(x,y) \\
  F([\godel{\forall \v_i}] \conc x,y) &= \Sc(F(x,y)) \quad \text{if $y \neq
    0$}
\end{align*}
Define $f_\Q(x) = F(\Q,x)$ and let $\Q_\exists$ be the number of
$\exists$-quantifiers in $\Q$, i.e., $\Q_\exists = H(\Q)$ where
\begin{align*}
  H([]) &=0, \\
  H([\godel{\exists \v_i}] \conc x) & = \Sc(H(x)) \quad \text{and}
  \\
  H([\godel{\forall \v_i}] \conc x) & = H(x).
\end{align*}

\begin{defin}
  A function $\Psi : \M \mapsto \M$ is a {\em Skolem operator} for the
  sequence of quantifiers $\Q$ if
\[
\mathop{\forall x,y \forall j \mathord< \Q_\exists} \Bigl( \forall i \mathord<
f_\Q(\Sc(j)) \bigl[(x)_i=(y)_i\bigr] \imp (\Psi(x))_j = (\Psi(y))_j \Bigr).
\]
\end{defin}

If $\varphi$ is a \LM-formula and $\Q$ is a sequence of quantifiers then let
$\Q\varphi$ denote the formula we get by preceding $\varphi$ with the
quantifiers $\Q$.\footnote{The exact definition of $\Q\varphi$ depends on the
  G\"odel numbering.}

Define two \PA-definable functions $G^\forall, G^\exists : \M \mapsto \M$ such
that the following is provable in \PA:
\begin{align*}
  G^\exists([],y) = G^\forall ([],y)&=0 \\
  G^\forall([\godel{\forall \v_i}]\conc x,0) &= i \\
  G^\exists([\godel{\exists \v_i}]\conc x,0) &= i \\
  G^\forall ([\godel{\exists \v_i}] \conc x,y) &= G^\forall(x,y) \\
  G^\exists ([\godel{\exists \v_i}] \conc x,\Sc(y)) &= G^\exists(x,y) \\
  G^\forall ([\godel{\forall \v_i}] \conc x,\Sc(y))&=G^\forall(x,y) \\
  G^\exists ([\godel{\forall \v_i}] \conc x,y)&=G^\exists(x,y).
\end{align*}

Also, define $g^\forall_\Q(x) = G^\forall(\Q,x)$ and $g^\exists_\Q(x) =
G^\exists(\Q,x)$. Informally, if $q$ is either $\forall$ or $\exists$ then
$g^q_\Q(x)=i$ if the $(x+1)$st $q$-quantifier in $\Q$ bounds the variable
$\v_i$.

If $\Q\varphi$ is a sentence and $\Psi$ is a Skolem operator for $\Q$ then
$\varphi[\Psi,a]$ will denote the \LM-sentence we get by substituting $\v_i$,
where
\[
i = (\mathrm{M} x) \bigl[ g^\forall_\Q(x)=l \bigr],
\]
by $(a)_l$ and $\v_j$, where
\[
j = (\mathrm{M} x) \bigl[ g^\exists_\Q(x)=l \bigr],
\]
by $(\Psi(a))_l$, in $\varphi$. Here $(\mathrm{M} x)\varphi(x)$ means `the
greatest $x$ such that $\varphi(x)$,' we are using it instead of $(\mu x)$ to
take care of situations like
\[
\forall \v_0 \exists \v_1 \forall \v_0 (\v_0 = \v_1)
\]
where it is the second $\forall$ quantifier bounding $\v_0$, not the first.

If $\Psi$ is a Skolem operator such that the function $\Psi : \M \mapsto \M$
is \LMs-definable then we say that $\Psi$ is a {\em definable Skolem
  operator}.

\begin{defin}[\cite{Krajewski:76}]
  A satisfaction class $\Sat$ is {\em complete with respect to definable
    Skolem operators} if
\begin{multline*}
  \Q \varphi \in \Sat \quad \text{iff} \quad \text{there exists a definable
    Skolem operator $\Psi$ for $\Q$ such that} \\
  \text{$\varphi[\Psi,a] \in \Sat$ for all $a \in \M$}.
\end{multline*}
\end{defin}

\begin{rem} In the standard model $\Nat$ of \PA\ the only satisfaction
  class,
\[
\Sat_0=\Set{ \godel{\varphi} | \varphi \in \ElDiag(\Nat)},
\]
is complete with respect to definable Skolem operators. This follows from the
more general fact that in any $\M \models \PA$ and \LMs-sentence $\Q\varphi$
we have
\begin{multline*}
  \M \models \Q\varphi \quad \text{iff} \quad \text{there exists a
    definable Skolem operator for $\Q$ such that} \\
  \text{$\M \models \varphi[\Psi,a]$ for all $a \in \M$.}
\end{multline*}
\end{rem}

Let us define the corresponding rule:
\begin{equation}
\tag{Skolem-rule}\label{skolem-rule} \inference{\Gamma,\Q
\varphi}{\inferforall{a \in \M}{\Gamma,\varphi[\Psi,a]}}
\end{equation}
where $\Psi$ is a definable Skolem operator for the sequence of quantifiers
$\Q$.

\begin{prop}
  A maximally consistent set in \M-logic with \ref{skolem-rule} is a
  satisfaction class complete with respect to definable Skolem operators.
\end{prop}
\begin{proof}
  Left to the reader.
\end{proof}

\section{Closure under predicate logic}\label{sec:predicate}

Let $\forpredprf{\sigma}(y)$ denote the formula expressing that there is a
(nonstandard) predicate logic proof of $y$ from hypothesis satisfying
$\sigma(x)$. The exact definition of the formula is left to the reader to
figure out. We will write $\Lambda \predprf \varphi$ for 
\[
\M \models \forpredprf{x \in \Lambda} (\varphi).
\]

We will add the rule
\begin{equation}
\tag{Pred}\label{pred-rule} \inference{\Gamma}{\inferforall{i <
a}{\Lambda, \varphi_i}} \qquad \text{if $\Set{\Vee \Lambda \vee
\varphi_0 ,\ldots, \Vee \Lambda \vee \varphi_{a-1}} \predprf \Vee
\Gamma$}.
\end{equation}
to \M-logic.

\begin{prop}
  Satisfaction classes closed under nonstandard first-order provability are
  precisely the maximally consistent sets in \M-logic with \ref{pred-rule}.
\end{prop}
\begin{proof}
  Left to the reader.
\end{proof}

Let $\PA^*$ be the set of all standard and nonstandard instances of the axioms
of \PA. 

\begin{thm}[\cite{Kotlarski:85}]
  If \Sat\ is a satisfaction class such that
\[
\pair{\M,\Sat} \models \forall x \bigl(\sent(x) \wedge \forpredprf{x \in
  \PA^*}(x) \imp x \in \Sat\bigr)
\]
then $\pair{\M,\Sat}$ satisfies $\Delta_0$-induction, i.e.,
\[
\pair{\M,\Sat} \models \varphi(0,\Sat) \wedge \forall x \bigl(\varphi(x,\Sat)
\imp \varphi(\Sc(x),\Sat)\bigr) \imp {\forall x}\, \varphi(x,\Sat)
\]
for every $\LMs \cup \Set{\Sat}$-formula $\varphi(x,\Sat)$ which is
$\Delta_0$.
\end{thm}

\begin{rem}
  The arithmetical part of any $\pair{\M,\Sat}$ satisfying the condition in
  the theorem is stronger than \PA; for example, the consistency of \PA\ is
  provable in such a model, since otherwise $0\neq 0 \in \Sat$.
\end{rem}

\begin{rem}
  The previous remark shows that some countable recursively saturated models
  of \PA\ does {\em not} admit satisfaction classes closed under nonstandard
  provability in \PA.  It might still be the case that any countable
  recursively saturated model of \PA\ admits a satisfaction class closed under
  nonstandard first-order provability since it may fail to include $\PA^*$.
\end{rem}

\chapter{Conclusion and further work}

In the first and second chapter we gave some background information, including
a short historical survey of the study of nonstandard truth. In
Chapter~\ref{chp:satcl} we introduced a new definition of satisfaction class
in a language with function symbols. We also discussed the drawback of
defining a satisfaction class as a set of pairs of formulas and elements. The
main part of the chapter led up to Theorem~\ref{thm:str.cons.crit}, some
applications were presented in the end of the chapter, such as the existence
of satisfaction classes making the sentence
$\Sc(\Sc(\ldots(\Sc(0))\ldots))=\c_a$ true for any nonstandard number of $\Sc$
symbols and any nonstandard $a$. Our definitions of \M-logic and template
logic is rather different from other authors and we think our notions is
easier to work with.

In the chapter that followed we introduced {\em free} satisfaction class, it is
a weaker notion than satisfaction classes and in some sense it is a more
natural notion, e.g., free \M-logic is more natural than \M-logic. We proved
one characterisation of free satisfaction classes in terms of free \M-logic.

Chapter~\ref{chp:stsatcl} presented some ideas of how to remove pathologies.
One famous pathology is $0 \neq 0 \vee \ldots \vee 0 \neq 0$ which can be made true
for any nonstandard number of repetitions. We highlighted some other
pathologies and gave ideas of how to remove those. We stated the question of
whether there are satisfaction classes closed under nonstandard propositional
(or predicate) proofs in any countable recursively saturated model.

We end this chapter be listing the open questions stated in the thesis.

\closeoutputstream{que} \input{questions.out}

\bibliography{ref} \bibliographystyle{alpha}
\end{document}